\documentclass[11pt]{article}
\usepackage{amsmath}
\usepackage{amsthm}
\input{epsf.tex}

\newtheorem{theorem}{Theorem}[section]
\newtheorem{proposition}[theorem]{Proposition}
\newtheorem{lemma}[theorem]{Lemma}
\newtheorem{corollary}[theorem]{Corollary}
\newtheorem{define}[theorem]{Definition}

\def\Empty{}

\catcode`\@=11

\def\section{\@startsection {section}{1}{\z@}{-3.5ex plus -1ex minus 
-.2ex}{2.3ex plus .2ex}{\large\bf}}

\def\fnum@figure{{\small Figure \thefigure}}
\def\fakefigure{\def\@captype{figure}}

\long\def\@makecaption#1#2{
    \vskip 10pt 
    \def\FCap{#2} \def\NoCap{\ignorespaces}
    \ifx \FCap\NoCap
       \setbox\@tempboxa\hbox{#1}  
      \else
       \setbox\@tempboxa\hbox{#1: \small \it #2}
    \fi
    \ifdim \wd\@tempboxa >\hsize   
        \unhbox\@tempboxa\par      
      \else                        
        \hbox to\hsize{\hfil\box\@tempboxa\hfil}  
    \fi}

\def\@oddhead{\hbox{}\rightmark \hfil \rm\thepage}
\def\sectionmark#1{\markright {\sc{\ifnum \c@secnumdepth >\z@
      \S\thesection.\hskip 1em\relax \fi #1}}}

\catcode`\@=12

\def\oplabel#1{
  \def\OpArg{#1} \ifx \OpArg\Empty {} \else
  	\label{#1}
  \fi}

\newlength{\saveu}

\catcode`\@=11

\newcommand{\ccl}{\mbox{${\cal L}$}}

\newcommand{\otsf}{\mbox{${\overline \tau^s_F}$}}
\newcommand{\otuf}{\mbox{${\overline \tau^u_F}$}}

\newcommand{\hhu}{\mbox{${\cal H} ^u$}}
\newcommand{\hhs}{\mbox{${\cal H} ^s$}}

\newcommand{\hp}{\mbox{${\cal H}$}}

\newcommand{\pin}{\mbox{$\partial _{\infty}$}}
\newcommand{\cc}{\mbox{$\cal C$}}

\newcommand{\oo}{\mbox{$\cal O$}}
\newcommand{\oos}{\mbox{${\cal O}^s$}}
\newcommand{\oou}{\mbox{${\cal O}^u$}}

\newcommand{\rrrr}{\mbox{${\bf R}$}}

\newcommand{\mi}{\mbox{$\widetilde M$}}

\newcommand{\ls}{\mbox{$\Lambda^s$}}
\newcommand{\lu}{\mbox{$\Lambda^u$}}

\newcommand{\wls}{\mbox{$\widetilde \Lambda^s$}}
\newcommand{\wlu}{\mbox{$\widetilde \Lambda^u$}}
\newcommand{\wcl}{\mbox{$\widetilde {\cal L}$}}
\newcommand{\wlsf}{\mbox{$\widetilde \Lambda^s_F$}}
\newcommand{\wluf}{\mbox{$\widetilde \Lambda^u_F$}}
\newcommand{\wlse}{\mbox{$\widetilde \Lambda^s_E$}}
\newcommand{\wlsg}{\mbox{$\widetilde \Lambda^s_G$}}
\newcommand{\wlsfi}{\mbox{$\widetilde \Lambda^s_{F_i}$}}
\newcommand{\wlsl}{\mbox{$\widetilde \Lambda^s_L$}}

\newcommand{\wws}{\mbox{$\widetilde W^s$}}

\newcommand{\wwp}{\mbox{$\widetilde \Phi$}}
\newcommand{\wwr}{\mbox{$\widetilde \Phi_{\rrrr}$}}

\newcommand{\hh}{\mbox{${\bf H}^2$}}

\newcommand{\si}{\mbox{$S^2_{\infty}$}}
\newcommand{\su}{\mbox{$S^1_{\infty}$}}

\newcommand{\fol}{\mbox{$\cal F$}}
\newcommand{\gal}{\mbox{$\cal G$}}

\newcommand{\fn}{\mbox{$\widetilde {\cal F}$}}

\newcommand{\ws}{\mbox{$\widetilde  W^s$}}
\newcommand{\wu}{\mbox{$\widetilde  W^u$}}

\def\@evenhead{\rm \leftmark \hfil \thepage}
\def\chaptermark#1{\markboth {\sc {\ifnum \c@secnumdepth >\m@ne
      \@chapapp\ \thechapter. \ \fi #1}}{}}%

\catcode`\@=12

\def\centeredepsfbox#1{\centerline{\epsfbox{#1}}}

\begin{document}

\title{Geometry of foliations and flows I: 
Almost transverse pseudo-Anosov flows and 
asymptotic behavior 
of foliations}
\author{S\'{e}rgio R. Fenley
\thanks{Reseach partially supported by NSF grants 
DMS-0296139 and DMS-0305313.}
}
\maketitle

\vskip .2in

{\small{
\centerline{\bf {Abstract}}
Let $\fol$ be a foliation in a 
closed  $3$-manifold with negatively curved fundamental group
and suppose that $\fol$ is almost transverse to a
quasigeodesic pseudo-Anosov flow.
We show that the leaves of the foliation
in the universal cover extend continuously to the sphere
at infinity, therefore the limit sets of the leaves are continuous images
of the circle.
One important corollary is that  if
$\fol$ is  a
Reebless, finite depth foliation in a hyperbolic manifold, then
it has the continuous extension property.
Such finite depth foliations 
exist whenever
the second Betti number is non zero.
The result also applies to other classes of foliations, including
a large class of foliations where all leaves
are dense and infinitely many examples with one sided
branching.
One key tool is a detailed understanding
of asymptotic properties of almost pseudo-Anosov singular
$1$-dimensional foliations in the leaves of $\fol$ lifted
to the universal cover. 
}}

\section{Introduction}

A $2$-dimensional foliation in a $3$-manifold is called Reebless
if it does not have a Reeb component: a foliation of the solid
torus so that the boundary is a leaf and the interior is foliated
by plane leaves spiralling towards the boundary.
As such the boundary leaf does not inject in the fundamental
group level and is compressible. Novikov \cite{No} showed that 
Reebless foliations and the underlying manifolds have excellent
topological properties. This result was extended by Rosenberg
\cite{Ros}, Palmeira \cite{Pa}
and many others.

The goal of this article is to analyse geometric properties of
foliations. Let $\fol$ be a Reebless 
foliation in $M^3$ with negatively curved fundamental
group.
Reebless implies that $M$ is irreducible \cite{Ros}.
In this article we will not make use of Perelman's 
fantastic results \cite{Pe1,Pe2,Pe3}, which if confirmed
imply that the manifold is hyperbolic.
Reebless foliations exist for instance  whenever $M$ is irreducible,
orientable and
the second homology of $M$ is not finite \cite{Ga1,Ga3}.
They also exist in much more generality by work of Roberts 
\cite{Ro1,Ro2,Ro3}, Thurston \cite{Th5} and many others.

Let 
$M^3$ be closed, irreducible with negatively curved 
fundamental group.
The universal cover is canonically compactified with a sphere
at infinity (denoted by $\si$), with compactification
a closed ball \cite{Be-Me}. The covering transformations
act by homeomorphisms in the compactified space.
Let $\fn$ be the lifted foliation to the universal cover $\mi$.
The leaves of $\fn$ are topological planes \cite{No}
and they are properly embedded. 
Hence they only limit
in the sphere at infinity. For hyperbolic manifolds,
the relationship between objects
in hyperbolic $3$-space (isometric to $\mi$) and their
limit sets in the sphere at infinity is central to  
the theory of 
such manifolds \cite{Th1,Th2,Mar}. The same is true
if $\pi_1(M)$ is negatively curved \cite{Gr,Gh-Ha}.
There is
a metric in $M$ so that all leaves of $\fol$ are hyperbolic
(that is constant curvature $-1$) \cite{Ca} and so the universal
cover of each leaf of 
$\fol$ is isometric to the hyperbolic  plane ($\hh$).
The {\em continuous extension question} asks whether these
leaves extend continuously to the sphere at infinity,
that is: given the inclusion map from a leaf $F$ of $\fn$
to $\mi$ is there a continuous extension to a map
$F \cup \pin F$ to $\mi \cup \si$? 
Here $\pin F$ is the ideal boundary of $F$ which is homeomorphic
to a circle.
If this is true we say that $\fol$ has the 
{\em continuous extension property}.
In that case the 
restriction of the map to $\pin F$ expresses the limit
set of $F$ as the continuous image of a circle, showing
it is locally connected.

In this article we prove the continuous extension property
for a very large class of foliations.
A pseudo-Anosov flow is {\em almost transverse} to a foliation
if an appropriate blow up of the flow is transverse 
to the foliation.
A blow up is obtained by replacing a (possibly empty) collection
of singular orbits by a union of annuli.
Another property that is important is a metric property:
A flow is {\em quasigeodesic} if it is uniformly
efficient up to a bounded multiplicative
distortion in measuring distances in the universal cover.
This is extremely important for manifolds with
negatively curved fundamental group \cite{Th1,Th3,Gr}.
Our main result is the following:

\vskip .1in
\noindent
{\bf {Main theorem}} $-$ Let $\fol$ be a foliation in
$M^3$ closed, atoroidal.
Suppose that $\fol$ is almost
transverse to a quasigeodesic pseudo-Anosov flow $\Phi$,
which has some prong singularity (that is, not a topological
Anosov flow).
This implies that $M$ has
negatively curved
fundamental group. 
Then $\fol$ has the continuous extension property.
Therefore the limits sets of leaves of $\fn$ are locally 
connected.
The set of foliations almost transverse to a flow
is open in the Hirsch topology of foliations.
\vskip .1in

Notice that the hypothesis imply that $\fol$ is transversely
orientable.
Since $M$ has a singular pseudo-Anosov flow then
$M$ is irreducible and the stable/unstable foliations
of $\Phi$ split to genuine laminations in $M$.
A fundamental result of Gabai and Kazez \cite{Ga-Ka}
then implies that $M$ has negatively curved fundamental
group. For simplicity of statements we usually
use the group negative curvature hypothesis,
but in most places that could be substituted
by the atoroidal hypothesis.


Notice that it is not necessary to assume 
that $\fol$ is
Reebless $-$ we prove that 
the condition of being almost transverse to a pseudo-Anosov
flow implies that $\fol$ is Reebless.

As a first consequence we prove the continuous
extension property for all Reebless finite depth foliations
in hyperbolic $3$-manifolds. 
Roughly a foliation is {\em finite depth} if all leaves
are proper and the leaves are perfectly fitted along
the cutting surfaces of a hierarchy of the manifold.
In particular there are
compact leaves.
These foliations exist whenever the second
homology is not finite.

\vskip .1in
\noindent
{\bf {Corollary A}} $-$ Let $\fol$ be a Reebless finite depth 
foliation in $M^3$ closed hyperbolic. Then $\fol$ has
the continuous extension property. In particular the
limit sets of the leaves are all locally connected.
\vskip .1in

Hence any orientable, hyperbolic $3$-manifold with non finite
second homology has such a foliation with the continuous
extension property. 
Conjecturally any closed, hyperbolic $3$-manifold
has a finite cover with positive first Betti number.
This would
imply there is always a foliation with the continuous extension
property in a finite cover.
The proof of corollary A is simple given previous results:
If necessary take a double cover and assume
that $\fol$ is transversely oriented. Then
Mosher and Gabai proved that such $\fol$ is almost transverse
to a pseudo-Anosov flow $\Phi$ \cite{Mo5}. 
We proved, jointly with Mosher that these
flows 
are quasigeodesic \cite{Fe-Mo}.
The main theorem then implies corollary A.
The result in \cite{Fe-Mo} is only for finite depth foliations: 
the proof depends heavily on the existence of a compact
leaf, $M$ being hyperbolic and the direct association with
a hierarchy.
By Thurston's geometrization theorem 
\cite{Th1,Th2,Mor}
it suffices to assume that
$M$ is atoroidal.

\vskip .1in
\noindent
{\bf {Corollary B}} $-$ There are infinitely many 
foliations with all leaves dense which have the continuous
extension property. Many of these have one sided branching.
\vskip .1in

Foliations with all leaves dense can be obtained for example
starting with finite depth foliations and doing small perturbations 
$-$ keeping it still almost transverse to the same 
quasigeodesic pseudo-Anosov flow.
A construction is 
carefully explained by Gabai \cite{Ga3}, providing
infinitely many examples with dense leaves to which 
corollary B applies. 
The examples occur whenever the second Betti number
of $M$ is non zero.
In fact whenever a foliation $\fol$ satisfies the
hypothesis of the main theorem, then any $\fol'$ sufficiently
close to $\fol$ will also be transverse to the same
flow. By the main theorem again, 
$\fol'$  will have the continuous extension
property.
This perturbation feature of the main theorem is not
shared by any previous result
concerning the continuous extension property.

A foliation $\fol$ is ${\bf R}$-{\em covered} if 
the leaf space  of $\fn$ is homeomorphic to the
real numbers. Equivalently this leaf space is Hausdorff.
A foliation which is not $\rrrr$-covered has {\em branching},
that is there are non separated points in the leaf space.
This leaf space is oriented (being a simply connected, perhaps
non Hausdorff $1$-manifold) and there is a notion of branching
in the positive or negative directions. If it branches only
in one direction the foliation is said to have 
{\em one sided branching}.
Foliations with one
sided branching, where all leaves are dense and the
foliation is transverse to a suspension
pseudo-Anosov flow (which is quasigeodesic)  were constructed
by Meigniez \cite{Me}. This provides infinitely 
many examples with one sided 
branching to which corollary B applies.

A very important tool in the proof of the main theorem is
an analysis of the topological structure
of the pseudo-Anosov flow. 
Let $\Phi_1$ be the original pseudo-Anosov flow almost transverse
to $\fol$. To make the flow transverse to $\fol$ one
needs in general to blow up a collection of singular orbits
into a collection of flow saturated annuli so that
each boundary is a closed orbit of the new flow $\Phi$.
The blown up flow is called an {\em almost pseudo-Anosov flow}
(see section 3).
If $\wwp$ is the lifted flow to
the universal cover $\mi$ and $\oo$ is its orbit space, then
$\oo$ is homeomorphic to 
the plane $\rrrr^2$ \cite{Fe-Mo} $-$ this is true for
pseudo-Anosov and almost pseudo-Anosov flows. 
When one blows up some singular orbits into a collection of
joined annuli, the stable/unstable singular foliations also
blow up. The two new singular foliations $\Lambda^s, \Lambda^u$
are everywhere transverse
to each other
except at the singularities and the blown up annuli. The blown
up annuli are part of both singular foliations. 
Since
$\fol$ is transverse to the blown up foliations, then the
stable/unstable foliations $\Lambda^s, \Lambda^u$
 induce singular $1$-dimensional foliations
in the leaves of $\fol$ and hence also in the leaves
$\fn$. The behavior of this is
described in the following result, which is
of independent interest also: 

\vskip .1in
\noindent
{\bf {Theorem C}} $-$ Let $\fol$ be a foliation 
with hyperbolic leaves in $M^3$
closed. Let $\Phi_1$ be a pseudo-Anosov flow almost transverse
to $\fol$ and let $\Phi$ be a corresponding almost pseudo-Anosov flow
transverse to $\fol$. Let $\Lambda^s, \Lambda^u$ be
the stable/unstable $2$-dimensional
foliations of $\Phi$ and $\wls, \wlu$
the lifts to $\mi$. Given $F$ leaf of $\fn$,
let $\wlsf, \wluf$ be the induced singular 
$1$-dimensional foliations in $F$.
Then 

\begin{itemize}

\item
For every ray $l$ in a leaf of $\wlsf$ or $\wluf$,
then $l$ limits in a single point of $\pin F$. 

\item
If 
the stable/unstable foliations $\wls, \wlu$ of $\Phi$ have 
Hausdorff leaf space, then the leaves of $\wlsf, \wluf$
are uniform quasigeodesics in $F$, the bound is independent
of the leaf. In general the leaves of $\wlsf, \wluf$ are not quasigeodesic.

\item
Any non Hausdorffness (of say $\wlsf$) 
is associated to a Reeb annulus
in a leaf of $\fol$ and when projected to $M$ it either projects
to or spirals
to a Reeb annulus. 

\item
The set of ideal points of
leaves of $\wlsf$ is dense in $\pin F$ and similarly for $\wluf$.

\item
If two rays of the same leaf of $\wlsf$ limit to the same ideal
point in $\pin F$ then this leaf of 
$\wlsf$ does not contain a singularity and the region in $F$
bounded by the leaf projects in $M$ to a set in a leaf
of $\fol$ which is either contained in or asymptotic to
a Reeb annulus.
\end{itemize}

\vskip .1in
This is a completely general result:
one does not need 
negatively curved fundamental group of $M$ or
any metric properties of $\Phi$.
Theorem C is a key tool used in
the proof of the main theorem.

The article is organized as follows:
In the next section we discuss previous results
on the continuous extension property and the strategy
of the proof of the main theorem and theorem C.
In section 3 we present basic definitions and
results concerning pseudo-Anosov flows and almost
pseudo-Anosov flows. In section 4 we analyse 
projections to the orbit space.
In sections 5 and 6 we analyse the singular foliations
$\wlsf, \wluf$ and asymptotic properties of their
leaves, proving theorem C.
In section 7 we prove the continuous extension property
$-$ the main theorem.
In the final section we comment on general relationships
between foliations and Kleinian groups.

In a subsequent article we analyse other important
consequences of quasigeodesic behavior for flows
and foliations \cite{Fe10}.


\section{Historical remarks and strategy of proofs}

Here we review what is known about the continuous extension
property.
In a seminal work, Cannon and Thurston \cite{Ca-Th} proved
this property 
when $\fol$ is a fibration over
the circle. Previously Thurston showed that a fiberting 
manifold is hyperbolic when the monodromy of the fibration
is pseudo-Anosov \cite{Th1,Th3,Th4,Bl-Ca}. Since the fundamental group
of a leaf of $\fol$ is a 
normal subgroup of the fundamental group of $M$,
then every limit set of a leaf of $\fn$ is the whole
sphere. In this way Cannon and Thurston produced many examples of 
group invariant Peano curves.

Another extremely important class of foliations is the
following:
A foliation is {\em proper} if the  leaves never limit
on themselves $-$ this is in the foliation sense and
it means that a sufficiently small transversal to a given 
leaf only
meets the leaf in a single point.
In particular leaves are not dense.
In a proper foliation there are compact leaves which
are said to have {\em depth} $0$. The depth of a leaf is
inductively defined to be $i$ (for finite $i$) if $i-1$ is the maximum
of the depths of leaves in the (foliation) limit set of
the leaf. A foliation has {\em finite depth} if it is
proper and there is a finite upper bound to the depths
of all leaves.

Gabai proved that whenever a compact $3$-manifold $M$ is 
irreducible, orientable  and
the second homology group $H_2(M, \partial M, {\bf Z})$ is not 
finite,
then there is a Reebless finite depth, foliation associated
to each non torsion homology class \cite{Ga1,Ga3}. 
The foliation is directly associated to a hierarchy of the
manifold and as such is strongly connected with the topological
structure of the manifold.
These results had several fundamental consequences for the
topology of $3$-manifolds \cite{Ga1,Ga2,Ga3}.
If $M$ is hyperbolic (or atoroidal), then
one important question is whether these finite depth
foliations have the continuous extension property.


Subsequently 
Gabai and Mosher showed \cite{Mo5} that any Reebless finite
depth foliation in a closed, atoroidal $3$-manifold admits
a pseudo-Anosov flow $\Phi$ which is almost transverse
to it. Roughly a flow is {\em pseudo-Anosov} if it has transverse
hyperbolic dynamics $-$ even though it may have
finitely many singularities.
It has stable and unstable two dimensional foliations which in
general are singular.
The term {\em almost transverse} means that one may need to 
blow up one singular orbit (or more) into a finite
collection of joined annuli to make the flow transverse
to the foliation. 
See detailed definitions and comments in section 3.
Since $\fol$ has a compact leaf and $M$ is
atoroidal, then Thurston \cite{Th1,Th3} proved
that $M$ is in fact hyperbolic.

We proved, jointly with Mosher,
that these pseudo-Anosov flows almost transverse to finite depth 
foliations
in hyperbolic $3$-manifolds are {\em quasigeodesic} \cite{Fe-Mo}.
This means that flow lines are uniformly efficient in
measuring distance in relative homotopy classes, or 
equivalently, uniformly efficient in measuring distance
in the universal cover.
This was first proved by Mosher \cite{Mo4} for a class
of flows transverse to some examples of depth one foliations
obtained by handle constructions.
Another concept is that of quasi-isometric behavior:
a foliation (perhaps singular) is {\em quasi-isometric} if
its leaves are uniformly efficient in measuring distance
in the universal cover. There are no non singular $2$ dimensional
quasi-isometric foliations in closed $3$-manifolds
with negatively curved fundamental group
\cite{Fe1}. 
As for {\underline {singular}} foliations the situation is
quite different and there are examples.
The stable/unstable singular foliations of the quasigeodesic
flows above may be quasi-isometric \cite{Fe6} and may
not \cite{Mo5,Fe6}. 
In general quasi-isometric behavior of $\Lambda^s$ (or $\Lambda^u$)
implies that $\Phi$ is quasigeodesic.

If both the stable and unstable foliations are quasi-isometric
{\underline {and}}
 the flow is actually transverse (as opposed to being
almost transverse) to the finite depth foliation then
we proved in \cite{Fe6} that $\fol$ has the continuous extension
property.
We stress that this result only applies to finite depth foliations
$-$ the proof depends, amongst other things,  on induction in the depth.
To use this result we needed to check
the quasi-isometric behavior of $\Lambda^s, \Lambda^u$ and
transversality between $\fol$ and $\Phi$.
This was very tricky and we could only do that for some 
{\underline {depth one}} foliations.
More to the point, it is known that  these conditions do not
always hold for general finite depth foliations.
Corollary A proves the continuous extension property
for all finite depth foliations: there are no restrictions
on the depth of the foliation, or about transversality
of the flow with the foliation or quasi-isometric behavior
of $\Lambda^s, \Lambda^u$.

\vskip .11in
The continuous extension property was also proved for
another class of foliations:
A foliation $\fol$ is {\em uniform} if any two leaves
in the universal cover are a bounded distance apart,  the
bound depends on the individual leaves. 
Hence $\fol$ is $\rrrr$-covered.
Thurston \cite{Th5} proved
that there is a large class of
uniform foliations.
If in addition
$\pi_1(M)$ is negatively
curved, then 
Thurston \cite{Th5} proved 
there is a pseudo-Anosov flow
transverse to $\fol$. From this  it is easy to 
prove  that the flow has quasi-isometric stable/unstable
foliations. In this case it also easily implies that the
foliation $\fol$ has
the continuous extension property.
The arguments are a very clever generalization of the
fibering situation.

Notice that in all the previous
results, there is a pseudo-Anosov flow $\Phi$
{\underline {transverse}} to $\fol$ and so that the
stable/unstable foliations of $\Phi$ are 
quasi-isometric singular foliations.
Both of these properties were crucial in all proofs.
The main theorem implies all the previous results
about the continuous extension property and
it has the unique feature that it applies
to an open set of foliations.

\vskip .12in
The main theorem 
can potentially be widely applicable because of
the abundance of pseudo-Anosov flows almost transverse
to foliations: Thurston proved this for fibrations \cite{Th4}.
It is also true for all
$\rrrr$-covered foliations \cite{Fe7,Cal1}
and Calegari proved it for all foliations with 
one sided branching \cite{Cal2}, all minimal foliations \cite{Cal3}
and many other foliations \cite{Cal3}.
The missing ingredient is the quasigeodesic property of
these pseudo-Anosov flows which is needed to apply
the main theorem.
In general the quasigeodesic property for a pseudo-Anosov
flow (or an arbitrary flow) is very hard to obtain.
This property is only known when the pseudo-Anosov
flow is almost transverse to
a foliation of one of the following types: 
finite depth or uniform.
One main goal in the study of a pseudo-Anosov flow $\Phi$ in
$M$ with Gromov hyperbolic fundamental group is to decide
whether $\Phi$ is quasigeodesic.
There are many examples where it is not quasigeodesic \cite{Fe2}.

\vskip .1in
Why almost transversality and not just transversality? 
In many cases ($\rrrr$-covered, 
one sided branching)  the pseudo-Anosov flow is actually
transverse to $\fol$.
But for finite depth foliations (which have
two sided branching), there are many examples where
it is impossible to make the pseudo-Anosov flow transverse
to the foliation \cite{Mo5} and one can only get almost
transversality. We will have more comments about this
in section 3.
Finite depth foliations are extremely important as they
are strongly connected to the topology of the underlying
manifold.
Also, in some sense, foliations with two sided branching
are probably more common than foliations which are
either $\rrrr$-covered or with
one sided branching.

We also remark that 
in all previous results concerning the continuous extension
property, theorem C was 
a crucial property on which the whole analysis hinged.
In the previous situations, the analysis of
leaves of $\wlsf, \wluf$ was either trivial or 
substantially easier.
In particular in these situations the leaves of
$\wlsf, \wluf$ were always uniform quasigeodesics, which
simplified subsequent proofs considerably. Such is
not the case in general.
In particular if $\wls, \wlu$ do not have Hausdorff leaf
space, then a priori the leaves of $\wlsf, \wluf$ do not have
to be quasigeodesic.
The proof of theorem $C$ works for any pseudo-Anosov
flow almost transverse to a foliation with hyperbolic
leaves.
The proof uses the
denseness of contracting directions for foliations
as proved by Thurston \cite{Th6,Th7} when he introduced the
universal circle for foliations $-$  even though we do not
directly use the universal circle here.
The basic idea is: if any ray 
of (say) $\wlsf$ does not limit to a single
point  in $\pin F$
then it limits in a non trivial interval of $\pin F$.
Zoom into this interval and analyse the situation 
in the limit. This is actually the easiest statement to prove
in theorem C.
The facts about rays with same ideal point
and non Hausdorffness are much trickier, but they will be essential
in the analysis of the main theorem.
The results of theorem C are also used in other contexts,
for example to analyse rigidity of pseudo-Anosov flows
almost transverse to a given foliation. This will be
explored in a future article \cite{Fe10}.

The proof of the main theorem has 2 parts: given a leaf $F$ of $\fn$,
one first constructs an extension to the ideal boundary and
then show it is continuous. To define the extension, one
uses the foliations $\wlsf, \wluf$ as they hopefully
define a basis neighborhood of an ideal point $p$ of $F$.
The best situation is that the leaves
of $\wls, \wlu$ which contain these leaves
of $\wlsf, \wluf$ define basis neighborhoods of unique
points in $\si$, hence defining the image of $p$ in $\si$.
There are several difficulties here: first the leaves
of $\wlsf, \wluf$ are not quasigeodesics in $F$, so much
more care is needed. Another problem is that
the foliations $\wls, \wls$ in general do not have
Hausdorff leaf space. This keeps recurring
throughout the proof. A further difficulty is that if
intersections of leaves $L_i$ of say
$\wls$ with a leaf $F$ of $\fn$ escape
in $F$ as $i$ converges to infinity, it
does not follow that the $L_i$ themselves escape compact
sets in $\mi$.
Consequently there
are several cases to be analysed.

Another fact that is important for the analysis of the main
theorem and theorem C is the following:
Let $\Theta$ be the projection map from $\mi$ to $\oo$
(projects flow lines to points).
A leaf of $\fn$ 
intersects an orbit of $\wwp$ at most once defining an injective
projection of $F$ to $\Theta(F)$. The projection $\Theta(F)$ is equal
to $\oo$ if and only if the foliation is $\rrrr$-covered.
An important problem here is to determine the boundary 
$\Theta(F)$ as a subset of $\oo$. This turns out to be
a collection of subsets of stable/unstable leaves in $\oo$.
This result is different than what happens
for pseudo-Anosov flows transverse to foliations
and its proof is much more delicate.
This is analsyed in section 4.

\section{Preliminaries: Pseudo-Anosov flows and almost pseudo-Anosov flows}


Let $\Phi$ be a flow on a closed, oriented 3-manifold $M$. We say
that $\Phi$ is a {\em pseudo-Anosov flow} if the following are
satisfied:


- For each $x \in M$, the flow line $t \to \Phi(x,t)$ is $C^1$,
it is not a single point,
and the tangent vector bundle $D_t \Phi$ is $C^0$.

- There is a finite number of periodic orbits $\{ \gamma_i \}$,
called {\em singular orbits}, such that the flow is 
``topologically" smooth off of the
singular orbits (see below).

- The flowlines are tangent to two singular transverse
foliations $\ls, \lu$ which have smooth leaves off of $\gamma_i$
and intersect exactly in the flow lines of $\Phi$.
These are like Anosov foliations off of the singular orbits.
This is the topologically smooth behavior described above.
A leaf containing a singularity is homeomorphic 
to $P \times I/f$
where $P$ is a $p$-prong in the plane and $f$ is a homeomorphism
from $P \times \{ 1 \}$ to $P \times \{ 0 \}$.
In a stable leaf, $f$ contracts towards towards
the prongs and in an unstable leaf it expands away
from the prongs.
We restrict to $p$ at least $2$, that is, we do not allow
$1$-prongs.

- In a stable leaf all orbits are forward asymptotic,
in an unstable leaf all orbits are backwards asymptotic.

Basic references for pseudo-Anosov flows are \cite{Mo3,Mo5}
and for $3$-manifolds \cite{He}.

\vskip .05in
\noindent
\underline {Notation/definition:} \ 
The singular
foliations lifted to $\mi$ are
denoted by $\wls, \wlu$.
If $x \in M$ let $W^s(x)$ denote the leaf of $\ls$ containing
$x$.  Similarly one defines $W^u(x)$
and in the
universal cover $\ws(x), \wu(x)$.
Similarly if $\alpha$ is an orbit of $\Phi$ define
$W^s(\alpha)$, etc...
Let also $\wwp$ be the lifted flow to $\mi$.

\vskip .1in
We review the results about the topology of
$\wls, \wlu$ that we will need.
We refer to \cite{Fe4,Fe6} for detailed definitions, explanations and 
proofs.
The orbit space of $\wwp$ in
$\mi$ is homeomorphic to the plane $\rrrr^2$ \cite{Fe-Mo}
and is denoted by $\oo \cong \mi/\wwp$. 
Let $\Theta: \mi \rightarrow \oo \cong \rrrr^2$
be the projection map. 
If $L$ is a 
leaf of $\wls$ or $\wlu$,
then $\Theta(L) \subset \oo$ is a tree which is either homeomorphic
to $\rrrr$ if $L$ is regular,
or is a union of $p$-rays all with the same starting point
if $L$ has a singular $p$-prong orbit.
The foliations $\wls, \wlu$ induce $1$-dimensional foliations
$\oos, \oou$ in $\oo$. Its leaves are $\Theta(L)$ as
above.
If $L$ is a leaf of $\wls$ or $\wlu$, then 
a {\em sector} is a component of $\mi - L$.
Similarly for $\oos, \oou$. 
If $B$ is any subset of $\oo$, we denote by $B \times \rrrr$
the set $\Theta^{-1}(B)$.
The same notation $B \times \rrrr$ will be used for
any subset $B$ of $\mi$: it will just be the union
of all flow lines through points of $B$.

\begin{define}
Let $L$ be a leaf of $\wls$ or $\wlu$. A slice of $L$ is 
$l \times \rrrr$ where $l$ is a properly embedded
copy of the reals in $\Theta(L)$. For instance if $L$
is regular then $L$ is its only slice. If a slice
is the boundary of a sector of $L$ then it is called
a line leaf of $L$.
If $a$ is a ray in $\Theta(L)$ then $A = a \times \rrrr$
is called a half leaf of $L$.
If $\zeta$ is an open segment in $\Theta(L)$ 
it defines a {\em flow band} $L_1$ of $L$
by $L_1 = \zeta \times \rrrr$.
Same notation for the foliations $\oos, \oou$ of $\oo$.
\end{define}

If $F \in \wls$ and $G \in \wlu$ 
then $F$ and $G$ intersect in at most one
orbit. 
Also suppose that a leaf $F \in \wls$ intersects two leaves
$G, H \in \wlu$ and so does $L \in \wls$.
Then $F, L, G, H$ form a {\em rectangle} in $\mi$
and there is  no singularity 
in the interior of the rectangle \cite{Fe6}.
There will be two generalizations of rectangles: 1) perfect fits $=$ rectangle
with one corner removed 
and 2) lozenges $=$ rectangle with two opposite corners removed.
We will also denote by rectangles, perfect fits, lozenges
and product regions the projection of these regions to
$\oo \cong \rrrr^2$.

\begin{define}{(\cite{Fe2,Fe4,Fe6})}{}
Perfect fits -
Two leaves $F \in \wls$ and $G \in \wlu$, form
a perfect fit if $F \cap G = \emptyset$ and there
are half leaves $F_1$ of $F$ and $G_1$ of $G$ 
and also flow bands $L_1 \subset L \in \wls$ and
$H_1 \subset H \in \wlu$,
so that 
%
the set 

$$\overline F_1 \cup \overline H_1 \cup 
\overline L_1 \cup \overline G_1$$

\noindent
separates $M$ and forms an a rectangle $R$ with a corner removed:
The joint structure of $\wls, \wlu$ in $R$ is that of
a rectangle with a corner orbit removed. The removed corner
corresponds to the perfect of $F$ and $G$ which 
do not intersect.
\end{define}

\begin{figure}
\centeredepsfbox{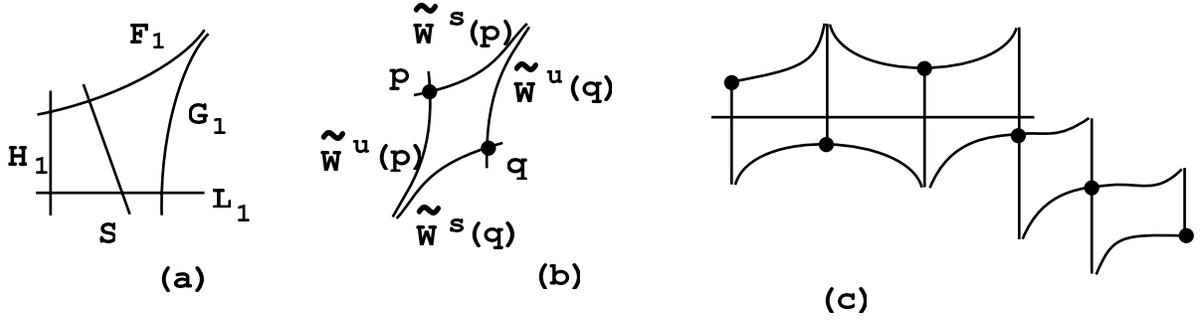}
\caption{a. Perfect fits in $\mi$,
b. A lozenge, c. A chain of lozenges.}
\label{loz}
\end{figure}

We refer to fig. \ref{loz}, a for perfect fits.
There is a product structure in the interior of $R$: there are
two stable boundary sides and two unstable one. An unstable
leaf intersects one stable boundary side (not in the corner) if
and only if it intersects the other stable boundary side
(not in the corner).
We also say that the leaves $F, G$ are {\em asymptotic}.

%
%
%
%

\begin{define}{(\cite{Fe2,Fe4,Fe6})}{}
Lozenges - A lozenge is a region of $\mi$ whose closure
is homeomorphic to a rectangle with two corners removed.
More specifically two points $p, q$ form the corners
of a lozenge if there are half leaves $A, B$ of
$\ws(p), \wu(p)$ defined by $p$
and  $C, D$ half leaves of $\ws(q), \wu(q)$ so
that $A$ and $D$ form a perfect fit and so do
$B$ and $C$. The region bounded by the lozenge
is $R$ and it does not have any singularities.
%
%
%
%
The sides are not contained in the lozenge,
but are in the boundary of the lozenge.
See fig. \ref{loz}, b.
\end{define}


There are no singularities in the lozenges,
which implies that
$R$ is an open region in $\mi$.
There may be singular orbits
on the sides of the lozenge and the corner orbits.


Two lozenges are {\em adjacent} if they share a corner and
there is a stable or unstable leaf
intersecting both of them, see fig. \ref{loz}, c.
Therefore they share a side.
A {\em chain of lozenges} is a collection $\{ \cc _i \}, 
i \in I$, where $I$ is an interval (finite or not) in ${\bf Z}$;
so that if $i, i+1 \in I$, then 
${\cal C}_i$ and ${\cal C}_{i+1}$ share
a corner, see fig. \ref{loz}, c.
Consecutive lozenges may be adjacent or not.
The chain is finite if $I$ is finite.

\begin{define}{}{}
Suppose $A$ is a flow band in a leaf of $\wls$.
Suppose that for each orbit $\gamma$ of $\wwp$ in $A$ there is a
half leaf $B_{\gamma}$ of $\wu(\gamma)$ defined by $\gamma$ so that: 
for any two orbits $\gamma, \beta$ in $A$ then
a stable leaf intersects $B_{\beta}$ if and only if 
it intersects $B_{\gamma}$.
%
%
This defines a stable product region which is the union
of the $B_{\gamma}$.
Similarly define unstable product regions.
\label{defsta}
\end{define}

There are no singular orbits of 
$\wwp$ in $A$.

%

We abuse convention and call 
a leaf $L$ of $\wls$ or $\wlu$ is called {\em periodic}
if there is a non trivial covering translation
$g$ of $\mi$ with $g(L) = L$. This is equivalent
to $\pi(L)$ containing a periodic orbit of $\Phi$.
In the same way an orbit 
$\gamma$ of $\wwp$
is {\em periodic} if $\pi(\gamma)$ is a periodic orbit
of $\Phi$.

We say that two orbits $\gamma, \alpha$ of $\wwp$ 
(or the leaves $\ws(\gamma), \ws(\alpha)$)
are connected by a 
chain of lozenges $\{ {\cal C}_i \}, 1 \leq i \leq n$,
if $\gamma$ is a corner of ${\cal C}_1$ and $\alpha$ 
is a corner of ${\cal C}_n$.

If ${\cal C}$ is a lozenge with corners $\beta, \gamma$ and
$g$ is a non trivial covering translation 
leaving $\beta, \gamma$ invariant (and so also the lozenge),
then $\pi(\beta), \pi(\gamma)$ are closed orbits
of $\wwp$ which are freely homotopic to the inverse of each
other.

\begin{theorem}{\cite{Fe4,Fe6}}{}
Let $\Phi$ be a pseudo-Anosov flow in $M^3$ closed and let 
$F_0 \not = F_1 \in \wls$.
Suppose that there is a non trivial covering translation $g$
with $g(F_i) = F_i, i = 0,1$.
Let $\alpha_i, i = 0,1$ be the periodic orbits of $\wwp$
in $F_i$ so that $g(\alpha_i) = \alpha_i$.
Then $\alpha_0$ and $\alpha_1$ are connected
by a finite chain of lozenges 
$\{ {\cal C}_i \}, 1 \leq i \leq n$ and $g$
leaves invariant each lozenge 
${\cal C}_i$ as well as their corners.
\label{chain}
\end{theorem}

A chain from $\alpha_0$ to $\alpha_1$ is called
{\em minimal} if all lozenges in the chain are distinct.
Exactly as proved in \cite{Fe4} for Anosov flows,
it follows that there
is a unique minimal chain from $\alpha_0$ to $\alpha_1$
and also all other chains have to contain all the lozenges
in the minimal chain.

\begin{figure}
\centeredepsfbox{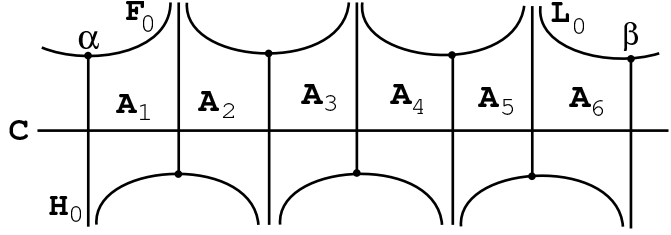}
\caption{
The correct picture between non separated
leaves of $\wls$.}
\label{pict}
\end{figure}

The main result concerning non Hausdorff behavior in the leaf spaces
of $\wls, \wlu$ is the following:

\begin{theorem}{\cite{Fe4,Fe6}}{}
Let $\Phi$ be a pseudo-Anosov flow in $M^3$. 
Suppose that $F \not = L$
are not separated in the leaf space of $\wls$.
Then the following facts happen:

\begin{itemize}

\item
$F$ is periodic and so is $L$.

\item
Let $F_0, L_0$ be the line leaves of $F, L$ which
are not separated from each other.
Let $V_0$ be the sector of $F$ bounded by
$F_0$ and containing $L$.
Let $\alpha$ be the periodic orbit in $F_0$ and
$H_0$ be the component of $(\wu(\alpha) - \alpha)$ 
contained in $V_0$.
Let $g$ be a non trivial covering translation
with $g(F_0) = F_0$, $g(H_0) = H_0$ and $g$ leaves
invariant the components of $(F_0 - \alpha)$.
Then $g(L_0) = L_0$.

\item
Let $\beta$ be the periodic orbit in $L_0$. As $g(\beta) = \beta$,
$g(\alpha) = \alpha$, then  $\pi(\alpha), \pi(\beta)$ are 
closed orbits of $\Phi$  which are  freely homotopic in $M$.
By theorem \ref{chain} \   $F_0$ and $L_0$ are connected by
a finite chain of lozenges 
$\{ A_i \}, 1 \leq i \leq n$.
Consecutive lozenges are adjacent.
They all intersect a common stable leaf $C$.
There is an even number of lozenges 
in the chain, see
fig. \ref{pict}.

\item
In addition 
let ${\cal B}_{F,L}$ be the set of leaves non separated
from $F$ and $L$.
Put an order in ${\cal B}_{F,L}$ as follows:
Put an orientation in the set of orbits of $C$ contained
in the union of the lozenges and their sides.
If $R_1, R_2 \in {\cal B}_{F,L}$ let $\alpha_1, \alpha_2$
be the respective periodic orbits in $R_1, R_2$. Then
$\wu(\alpha_i) \cap C \not = \emptyset$ and let 
$a_i = \wu(\alpha_i) \cap C$.
We define $R_1 < R_2$ in ${\cal B}_{F,L}$ 
if $a_1$ precedes
$a_2$ in the orientation of the set of orbits of 
$C$.
Then ${\cal B}_{F,L}$
is either order isomorphic to $\{ 1, ..., n \}$ for some
$n \in {\bf N}$; or ${\cal B}_{F,L}$ is order
isomorphic to the integers ${\bf Z}$.

\item
Also if there are $Z, S \in \wls$ so that
${\cal B}_{Z, S}$ is infinite, then there is 
an incompressible torus in $M$ transverse to 
$\Phi$. In particular $M$ cannot be atoroidal.
Also if there are $F, L$ as above, then there are
closed orbits $\alpha, \beta$ of $\Phi$ which
are freely homotopic to the inverse of each other.

\item
Up to covering translations,
there are only finitely many non Hausdorff
points in the leaf space of $\wls$.
\end{itemize}
\label{theb}
\end{theorem}

%

Notice that ${\cal B}_{F,L}$ is a discrete set in this order.
For detailed explanations and proofs, see
\cite{Fe4,Fe6}. 

\begin{theorem}{(\cite{Fe6})}{}
Let $\Phi$ be a pseudo-Anosov flow. If there is
a stable or unstable product region, then $\Phi$ is 
topologically conjugate to a suspension Anosov flow.
In particular $\Phi$ is non singular.
\label{prod}
\end{theorem}

%

\begin{proposition}{}{}
Let $\varphi$ be a (topological) Anosov flow so that every leaf of
its stable foliation $\wls$ intersects every leaf of its stable foliations
$\wlu$. Then $\varphi$ is topologically conjugate to a suspension 
Anosov flow. In particular $M$ fibers over the circle with
fiber a torus and Anosov monodromy.
\label{susp}
\end{proposition}

\begin{proof}{}
This result is proved by Barbot \cite{Ba1} when $\varphi$
is a {\underline {smooth}} Anosov flow. That means it is $C^1$ and
it has also {\underline {strong}} stable/unstable foliations and
contraction on the level of tangent vectors along the flow.
Here we only have the weak foliations and orbits being asymptotic
in their leaves. 
With proper understanding all the steps carry through to the
general situation.

Lift to a finite cover where $\Lambda^s, \Lambda^u$ are transversely orientable.
A cross section in the cover projects to a cross section in
the manifold (after cut and paste following Fried \cite{Fr}) and
so we can prove the result in the cover.

First, the flow $\varphi$
is expanding: there is $\epsilon >  0$ so that no distinct 
orbits are always less than $\epsilon$ away from each other. Inaba
and Matsumoto
then proved that this flow is a topological pseudo-Anosov flow \cite{In-Ma}.
The main thing is the existence of a Markov partition
for the flow. This implies that if $F$ is a leaf of $\wls$ which
is left invariant by $g$, then there is a closed orbit of $\varphi$
in $\pi(F)$ and all orbits are asymptotic to this closed orbit.
Similarly for $\wlu$.

What this means is the following: consider the action of $\pi_1(M)$ 
in the leaf space of $\wls$ which is the reals.
Hence we have a group action in $\rrrr$. Let $g$ in $\pi_1(M)$
which fixes a point. There is $L$ in $\wls$ with $g(L) = L$.
So there is orbit $\gamma$ of $\widetilde \varphi$ with $g(\gamma) = \gamma$.
Let $U$ be the unstable leaf of $\widetilde \varphi$ with $\gamma$ contained
in $U$. Then $g(U) = U$. If $g$ is associated to the positive direction
of $\gamma$ then $g$ acts as a contraction in the set of orbits of $U$
with $\gamma$ as the only fixed point. Since every leaf
of $\wlu$ intersects every leaf of $\wls$ then the set of orbits
in $U$ is equivalent to the set of leaves of $\wls$.
This implies the important fact:

\vskip .1in
\noindent
{\underline {Conclusion}} - If $g$ is in $\pi_1(M)$ has a fixed point
in the leaf space of $\wls$ then it is of hyperbolic type and
has a single fixed point.
\vskip .1in

Using this topological characterization Barbot \cite{Ba1} showed that
$G = \pi_1(M)$ is metabelian, in fact he showed that the commutator
subgroup $[G,G]$ is abelian. In particular $\pi_1(M)$ is solvable.
This used only an action by homeomorphisms in $\rrrr$ satisfying the
conclusion above.
Barbot \cite{Ba1} also proved 
that the leaves of $\Lambda^s, \Lambda^u$ are
dense in $M$.

Plante \cite{Pl1}, showed that if $\fol$ a minimal foliation
in $\pi_1(M)$ solvable then $\fol$ is transversely affine:
there is a collection of charts $f_i:U_i \rightarrow \rrrr ^2 \times \rrrr$,
so that the transition functions are affine in the second coordinate.
Using this Plante \cite{Pl1, Pl2} constructs a homomorphism

$$C: \ \pi_1(M) \ \ \rightarrow \ \ \rrrr$$

\noindent
which measures the logarithm of how much distortion there is
along an element of $\pi_1(M)$. This is a cohomology class.
Every closed orbit $\gamma$ of $\varphi$ has a transversal fence which
is expanding - this implies that $C(\gamma)$ is positive.
Plante then refers to a criterion of Fried \cite{Fr} to conclude that
$\varphi$ has a cross section and therefore it is easily
seen that $\varphi$ is topologically conjugate
to a suspension Anosov flow. 
This finishes the proof of the proposition.
\end{proof}

\vskip .1in
\noindent
{\bf {\underline {Blown up orbits and almost pseudo-Anosov flows}}}
\vskip .05in

We now describe almost pseudo-Anosov flows.
The description is taken from \cite{Mo5}.
First we describe the blow up of a singular orbit.
Let $\gamma$ be a singular orbit of a pseudo-Anosov flow
$\varphi$ and let $C$ be a small disk transverse to $\varphi$
through a point $p$ of $\gamma$. The set  \ $\ws(\gamma)
\cup \wu(\gamma)$ \ intersects $C$ in a tree $T$, which
is made of a single common vertex and $2n$ prongs,
if the stable/unstable leaves $\ws(\gamma), \wu(\gamma)$
have $n$ prongs.
The stable and unstable prongs alternate in $C$.
Orient the prongs so that stable prongs point
towards $p$ and unstable prongs point away from
$p$. Hence $T$ is an oriented tree.
Up to topological conjugacy we may assume in 
polar coordinates centered at $p$ that 
$T \ = \ \{ (r,\theta) \ | \ \theta = k\pi/n,
k = 0,1, ..., 2n-1 \}$.

Choose a Poincar\'{e} section for $\gamma$, that is,
a smaller disk $C'$ so that every point of $C'$ flows
forward to a point in $C$. This  defines a continuous 
map $f:C' \rightarrow C$ which fixes $p$ and is a
homeomorphism into its image.
This $f$ sends the corresponding pieces of
$T$ into $T$,
contracting stable direction and expanding unstable
directions.
There is a rotation or reflection $R$ in $C$ so
that $R$ prescribes where the sectors of $C$ defined
by $T$  go to under $f$.
For instance $R$ may be a rotation by $2\pi/n$.

We first blow up $f$.
Let $D$ be a small subdisk of $C'$ containing
$p$ in its interior.
Let $T^*$ be an oriented tree in $C$ which agrees
with $T$ outside $D$ and which is invariant
under $R$ as above and such that each vertex
$v$ of $T^*$ is ``pseudo-hyperbolic", meaning
as you go around the edges incident to $v$,
the orientations of the edges of $T^*$
alternate pointing toward and away from $v$.
Each vertex of $T^*$ must have $2i$ edges for
some $i \geq 2$.
The point here is that $T^*$ is created from
$T$ by replacing $p$ with a finite subtree.
There are finitely many ways to choose
$T^*$, up to compactly supported isotopy.
The new edges of $T^*$ created
by the blow up are called the finite edges of $T^*$.

Given this data there is a $C^0$ perturbation 
$f^*$ of $f$, and a continuous map $h: C \rightarrow C$
so that 

\begin{itemize}

\item
$f^*$ leaves $T^*$ invariant.

\item
For each edge $E$ of $T^*$, the first return map of 
$f^*$ acts as a translation on
$int(E)$ moving points in the direction of the
orientation on $E$.

\item 
$h$ collapses the finite edges of $T^*$  to the
point $p$ and $h$ is otherwise $1$ to $1$.

\item
$h$ is a semiconjugacy from $f^*$ to $f$,
that is, $f \circ h = h \circ f^*$.

\item
$h$ is close to the identity map in the sup norm
and $h$ equals the identity on 
$C - D$.

\end{itemize}

\noindent
We say that  $f^*$ is obtained from
$f$ by 
{\em dynamically blowing up}
the pseudo-hyperbolic fixed point $p$.
Each choice of $T^*$ determines a unique
$f^*$ up to conjugacy by isotopy.
There are therefore finitely many ways to dynamically
blow up a pseudo-hyperbolic fixed point,
up to conjugation by isotopy.

Now define a dynamic blow up of $\gamma$ by altering
the flow $\varphi$ near $\gamma$ as follows. Let $D'$ so that
$f(D') \subset D$.
Alter $\varphi$ so that the first return map
$f: D' \rightarrow D$ is replaced
by $f^*: D' \rightarrow D$.
This has the effect of alterating the generating vector
field of $\varphi$ inside the ``mapping torus"
$T_g = \{ \varphi(x,s) \ | \ x \in D', 0 \leq s \leq t(x) \}$
and leaving the vector field unaltered outside of
$T_g$.
Let $\varphi^*$ be the new flow.

The orbit $\gamma$ gets blown up into a collection
of flow invariant annuli. In each annulus $A$, the boundary
components are closed orbits of $\varphi^*$
which are isotopic to $\gamma$
as oriented orbits.
In the interior of $A$ orbits move from one boundary to the other,
as given by the orientation of the corresponding
edge of $T^*$.
There is a blow down map $\xi:M \rightarrow M$ which is
homotopic to the identity and isotopic to the
identity in the complement of the annuli.
It sends the collection of annuli into $\gamma$.
The map $\xi$ sends orbits of $\varphi^*$ to orbits
of $\varphi$ and preserves orientation.

\begin{define}{}{}
Let $\varphi$ be a pseudo-Anosov flow in $M^3$ closed.
Then $\varphi^*$ is an almost pseudo-Anosov flow associated to
$\varphi$ if $\varphi^*$ is obtained from $\varphi$ by dynamically blowing
some singular orbits of $\varphi$.
\end{define}


The reason for considering almost pseudo-Anosov flows is
as follows. All of the constructions of pseudo-Anosov flows
transverse to foliations are in fact constructions of
a pair of laminations $-$ stable and unstable $-$ which
are transverse to each other and to the foliation
\cite{Th4,Mo5,Fe7,Cal1,Cal2}.
The intersection of the laminations is oriented producing
a flow in this intersection.
One then collapses the complementary regions to the
laminations to produce transverse singular foliations and 
a pseudo-Anosov flow.
The transversality problem occurs in this last step, the
blow down of complementary regions. In certain situations,
for example for finite depth foliations, one cannot guarantee
transversality of the flow and foliation after the blow down. 
See extended explanations by Mosher in \cite{Mo5}.

The necessity of almost transversality as opposed to transversality
was first discovered by Mosher in \cite{Mo1} for positioning
surfaces with respect to certain ${\bf Z}$ invariant
flows (in non compact manifolds).
This was further explored in \cite{Mo2} where it was analysed
when it is necessary to 
blow up a pseudo-Anosov flow before
it became transverse to a given surface.

%
%
%
\begin{figure}
\centeredepsfbox{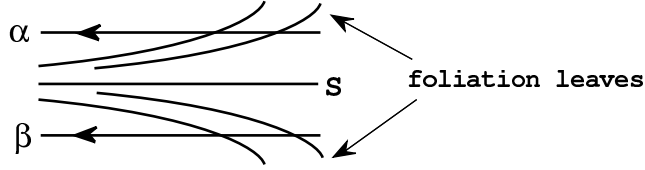}
\caption{
Obstruction to transversality.
Here $\alpha, \beta$ are isotopic closed orbits.
If the flow is to be pseudo-Anosov, then one has to
collapse $\alpha, \beta$ together. But then the collapsed
flow will not be transverse to the foliation.}
\label{prob}
\end{figure}

%
%

\vskip .1in
Given a dynamic blow up $\varphi^*$ of $\varphi$ with corresponding
blow down
map $\xi$, then $\xi^{-1}$ of $\Lambda^s, \Lambda^u$ are
singular foliations, which now have some flow invariant
annuli.
These are the annuli that come from blowing up some 
singular orbits $-$ they are called the
{\em blown up annuli}.
We still denote by $\ls, \lu$ the blown up stable/unstable foliations
of an almost pseudo-Anosov flows. They are transverse
to each other except at the blown up annuli. 
The same notation is used for $\wls, \wlu$, etc..

Because $\varphi^*$ is transverse to the foliation $\fol$,
then $\wls, \wlu$ are transverse to $\fn$ and induce
singular $1$-dimensional foliations $\wlsf, \wluf$
in any leaf $F$ of $\fn$.

The objects perfect fits, lozenges, product regions, etc.. all
make sense in the setting of almost pseudo-Anosov flows: they
are just the blow ups of the same objects for the corresponding
pseudo-Anosov flows. Since the interior of these objects does
not have singularities, the blow up operation does
not affect these interiors. 
There may be singular orbits in the boundary which get
blown into a collection of annuli. All the results in
this section still hold for almost pseudo-Anosov flows,
with the blow up operation.
For example if $F, L$ in $\wls$ are not separated from
each other, then they are connected by an even number
of lozenges all intersecting a common stable leaf.
Since parts of the boundary of these may have been
blown into annuli, there is not a product structure
in the closure of the union of the lozenges, but there
is still a product structure of $\wls, \wlu$ in the interior.

\section{Projections of leaves of $\fn$ to the orbit space}

\noindent
{\bf {Notation}} $-$ In some of the proofs and arguments that
follow we will be working with an almost pseudo-Anosov flow
transverse to a foliation $\fol$. 
In those cases, for notational simplicity 
we denote the almost pseudo-Anosov flow by
$\Phi$ and a corresponding blow down pseudo-Anosov
flow by $\Phi_1$. 
Hence $\Phi$ is some blow up of $\Phi_1$.
This is different than the
notation $\varphi, \varphi^*$ from the last section.

\vskip .1in
Let then $\Phi$ be an almost pseudo-Anosov flow which is transverse
to a foliation $\fol$. 
This implies that $\fol$ is Reebless $-$ we provide
a proof of this at the end of this section.
An orbit of $\wwp$ intersects a leaf of $\fn$ at most once $-$
because the leaves of $\fn$ are properly embedded and 
$\wwp$ is transverse to $\fn$. 
Hence the projection $\Theta: F \rightarrow \Theta(F)$ is injective.
We want to determine
the set of orbits a leaf of $\fn$ intersects $-$ 
in particular we want to determine the boundary $\partial \Theta(F)$.
As it turns out, $\partial \Theta(F)$ is composed of a disjoint
union of slice leaves in $\oos, \oou$.

We assume throughout this section that $\Phi$ is {\em blow down
minimal} with respect to being transverse to $\fol$. This means
that no blow down of some flow annuli of $\Phi$ produces a flow
transverse to $\fol$.

Since $\Phi$ is transverse to $\fol$, there is $\epsilon > 0$ so that
if a leaf $F$ of $\fn$ intersects an orbit of $\wwp$ at $p$ then 
it intersects every orbit of $\wwp$ which passes $\epsilon$ near
$p$ and the intersection is also very near $p$.
To understand $\partial \Theta(F)$ one main ingredient is that 
when considering pseudo-Anosov flows, then
flow lines in the same stable leaf are forward asymptotic. So
if $F$ intersects a given orbit in a very future time then it
also intersects a lot of other orbits in the same stable
since in future time
they converge. In the limit this produces a stable boundary leaf
of $\Theta(F)$. The blow up operation disturbs this: it is not 
true that orbits in the same stable leaf of an almost pseudo-Anosov
flow are forward asymptotic: when they pass arbitrarily near a blow
up annulus the orbits are distorted and their distance
can increase  enormously. This is the key difficulty in this section.
Hence we first analyse the blow up operation more
carefully.

\vskip .1in
\noindent
{\underline {Notation}} $-$ Given $\Phi$ an almost pseudo-Anosov flow,
let 
$\Phi_1$ be a corresponding pseudo-Anosov flow associated to $\Phi$. 
The term $\ws(x)$ will denote the stable leaf of $\wwp$ or 
$\wwp_1$, where the context will make clear which one it is.
\vskip .1in

Recall that $\pi: \mi \rightarrow M$ denotes the universal covering map.

We will start with $\Phi_1$ and understand the blow up procedure.
The blown up annuli come from singular orbits.
The {\em lift annuli} are the lifts of blown up annuli
to $\mi$. Their projections to $\oo$ are called {\em
blown segments}.
If $L$ is a blown up leaf of $\wls$ or $\wlu$ the components
of $L$ minus the lift annuli are called the {\em prongs}.
A {\em quarter} associated to an orbit $\gamma$ of $\wwp_1$
is the closure of a connected component of $\mi - (\wu(\gamma)
\cup \ws(\gamma))$. Its boundary is a union of $\gamma$ and
half leaves in the stable and unstable leaves of $\gamma$.
We will be interested in a neighborhood $V$ of $\gamma$ in
this quarter which projects to $M$ near the closed orbit
$\pi(\gamma)$. 
We will understand the blow up in the projection of a
quarter. Glueing up different quarter gives the overall picture
of the blow up operation.
In the projected quarter $\pi(V)$ in $M$ there
is a cross section to the flow $\Phi_1$. The orbits across
the cross section are determined by which stable and unstable
leaf they are in. The return map on the stable direction
is a contraction and an expansion in the unstable direction.
Any contraction is topologically conjugate to say
$x \rightarrow x/2$ and an expansion is conjugate to 
$x \rightarrow 2x$. Hence the local return map is 
topologically conjugate to

$$
\left[
\begin{array}{rr}
 1/2 & 0 \\ 0 & 2 \\ 
\end{array}
\right]$$

\noindent
a linear map. The whole discussion here is one of topological
conjugacy. The flow is conjugate to 
$(x,y,0) \ \rightarrow \ (2^{-t}x, 2^t y, t)$.
%
Think of the blow up annulus as the set of unit tangent vectors
to $\gamma$
associated to the quarter region. The flow in the
annulus is given by the action of $D V_t$ on the tangent vectors.
It has 2 closed orbits (the boundary ones corresponding to the stable
and unstable leaves). The other orbits are asymptotic to the
stable closed orbit in negative time and to the unstable
closed orbit in positive time. This makes it into a continuous flow
in this blown up part. 
For future reference recall this fact that in a blow up annulus the
boundary components are orbits of the flow and in the interior
the flow lines go from one boundary to the other without a Reeb annulus
picture (there is a cross section to the flow in the annulus).
Do this for each quarter region that is
blown up. One can then glue up the 2 sides of the appropriate annuli
because they are all of the same topological picture (using the
standard model above). This describes the blown
up operation in a quarter. There is clearly a blow down map which sends orbits
of the blown up flow $\Phi$
to orbits of $\Phi_1$ and collapses connected unions
of annuli into a single $p$-prong singular orbit.

We quantify these: let $\epsilon$ very small so that any
two orbits of $\Phi_1$ which are always less than $\epsilon$
apart in forward time, then they are in the same stable leaf.
Let ${\cal Z}$ the union of the singular orbits of 
$\Phi_1$ which are blown up.
Let $\epsilon' << \epsilon$ and let
$U$ be the $\epsilon'$ tubular neighborhood of ${\cal Z}$.
Let $U'$ (resp. $U$) be the $\epsilon'/2$ 
(resp. $\epsilon'$)
tubular neighborhood of ${\cal Z}$.
Choose the blow up map to be the identity in the complement
of $U'$, that is the blown up annuli are also contained in $U'$.
The blow down map is then 
an isometry of the Riemannian metric outside
$U'$. Choose the blow down to move points very little in $U'$.
Isotope $\fol$ so that it is transverse to the flow $\Phi$.
We are now ready to analyse $\partial \Theta(F)$.

\begin{proposition}{}{}
Let $F$ in $\fn$. Then $\Theta(F)$ is an open subset of $\oo$.
Any boundary component of $\Theta(F)$ is a slice of a leaf of 
$\oos$ or $\oou$.
If it is a slice of $\oos$, then as $\Theta(F)$ approaches $l$,
the corresponding points of $F$ escape in the positive direction.
Similarly for unstable boundary slices.
\label{bounda}
\end{proposition}

\begin{proof}{}
First notice that since $F$ is transverse to $\wwp$ then $\Theta(F)$ is
an open set. Hence $\partial \Theta(F)$ is disjoint from $\Theta(F)$.
The important thing is to notice that the metric is the same outside
the small neighborhood $U'$ of the blown up annuli. If two 
points are in the same stable leaf, then their orbits under
the blow down flow $\Phi_1$ are asymptotic in forward time.
The same is true
for $\Phi$, for big enough time if the point is outside $U$.
This is because the points of the
corresponding orbits of $\Phi_1$ will
be {\underline {both}} outside $U'$ $-$ this is the reason for
the construction of two neighborhoods $U', U$.
The following setup will be used in {\underline {all}} cases.

\vskip .1in
\noindent
{\underline {Setup}} $-$
Let $v$ in $\partial \Theta(F)$ and $v_i$ in $\Theta(F)$ with
$v_i$ converging to $v$.
Let $p_i$ in $F$ with $\Theta(p_i) = v_i$ and let $w$ in $\mi$
with $\Theta(w) = v$. Let $D$ be any small disk in $\mi$
transverse to $\wwp$ with $w$ in the interior of $D$. For $i$
big enough $v_i$ is in $\Theta(D)$ so there are 
$t_i$ real numbers with $p_i = \wwp_{t_i}(w_i)$
and $w_i$ are in $D$. As
$v$ is not in $\Theta(F)$, then $|t_i|$ grows without bound. Without loss
of generality assume up to subsequence that $t_i \rightarrow +\infty$.
Same proof if $t_i \rightarrow -\infty$.
We will prove that there is a slice leaf $L$ of $\ws(w)$ so that 
$\Theta(L) \subset 
\partial \Theta(F)$ and $F$ goes up as it ``approaches" $L$.
The stable/unstable leaves here are those of the almost pseudo-Anosov flow
and they may have blown up annuli.

One transversality fact used here is the following:
for each $\epsilon$ sufficiently small, there is $\epsilon' > 0$
so that if $x, y$ are $\epsilon$ close in $M$ then
the $\fol$ leaf through $x$ intersects the $\Phi$-orbit through
$y$ in a point $\epsilon'$ close to $x$. 
The $\epsilon'$ goes to $0$ as $\epsilon$ does.

\vskip .1in
\noindent
{\underline {Case 1}} $-$ Suppose that $w$ is not in a blown up leaf.

First we show that we can assume no $w_i$ is in 
$\ws(w)$. Otherwise up to subsequence assume all $w_i$ are in
$\ws(w)$. The orbits through $w_i$ and $w$ start out very close
and aside from the time they stay in $\pi^{-1}(U)$ they are always very close.
Let $B$ be the component
of the intersection 
of $F$ with the flow band from $
\wwr(w_i)$ to $\wwr(w)$ in the stable leaf $\ws(w)$,
which contains $p_i$.
Then $B$ does not intersect $\wwr(w)$ so it has
to either escape up or down.
If it escapes down it will have to intersect a small segment
from $w_i$ to $w$ and hence so does $F$. For $i$ big enough
$w_i$ is arbitrarily near $w$, so
transversality of $\fol$ and $\Phi$ then
implies that $F$ will intersect $\wwr(w)$ near $w$, contradiction
see fig. \ref{neck}, a.

We now consider the case that $B$ escapes up. If the forward orbit
through $w$ is not always in $\pi^{-1}(U)$ then at those times outside
of $\pi^{-1}(U)$ it will be arbitrarily
close to $\wwr(w_i)$ and transversality
implies again that $F$ intersects $\wwr(w)$. If the forward orbit
of $w$ always stays in $\pi^{-1}(U)$ the same happens after the blow down
so the blow down orbit is in the stable leaf of the singular
orbit which is being blown up. 
This does not happen in case 1.

We can now assume that all $v_i$ are in a sector of $\oos(v)$
with $l$ the boundary of this sector and 
$L = l \times \rrrr$,
the line leaf of
$\ws(w)$ which is the boundary of this sector.

Let now $q$ in $l$. We will show that $q$ is in $\partial \Theta(F)$
so $l \subset \partial \Theta(F)$. There is a segment $[q,v]$ contained
in $l$. 
Choose $x$ in $L$ with $\Theta(x) = q$. Let $\alpha$ be a segment
in $\ws(w)$ transverse to the flow lines and going from $x$ to $w$. 
Notice that $x, w$ are in $\mi$ and $q, v$ are in the orbit
space $\oo$.
Let $x_i$ converging to $x$ and
$x_i$ in $\ws(w_i)$. We can do that since all $w_i$
are in the same sector of $\ws(w)$.
 Choose segments $\alpha_i$ from $x_i$ to $w_i$ in 
$\ws(w_i)$ and transverse to the flowlines of $\wwp$ in $\ws(w_i)$.

\vskip .1in
\noindent
{\underline {Claim}} $-$ For every orbit $\gamma$ 
of $\wwp$ intersecting $\alpha_i$
in $y$ then $\gamma$ intersects $F$ in $\wwp_s(y)$ where
$s$ converges to $\infty$ as $i \rightarrow \infty$.


Suppose there is $a_0 > 0$ so that for some $i_0$ then

$$\wwp_{[a_0,t_i]}(w_i) \ \subset \ \pi^{-1}(U)
\ \ \ {\rm for \ all} \ \ i \geq i_0$$

\noindent
Then $\wwp_{[a_0,\infty)}(w)$ is contained in the closure
of $\pi^{-1}(U)$.
As seen before this implies that $w$ is in a blown up
stable leaf, which is not the hypothesis of case 1.
Therefore up to subsequence, 
there are arbitrary
big times $s_i$ between $0$ and $t_i$ so that $\wwp_{s_i}(w_i)$ is
not in $\pi^{-1}(U)$. 
Hence $\wwr(x_i)$ is very close to $\wwp_{s_i}(w_i)$
and since $F$ cannot escape up or down then $F$ intersects 
$\wwr(x_i)$. Hence the segment
$[\Theta(x_i), v_i]$ of $\oos(v)$
is contained in $\Theta(F)$ and
so $[x,v]$ is contained in the closure of $\Theta(F)$. Also
the time $s$ so that $\wwp_s(y)$ hits $F$ goes to $\infty$, hence
$[x,v]$ cannot intersect $\Theta(F)$ $-$ else there would
be bounded times where it intersects $F$, by transversality.
We conclude that 
$[x,v] \subset \partial \Theta(F)$, hence $l \subset \partial \Theta(F)$
as desired.
If there is a sequence $z_i$ in $F$
escaping down with $\Theta(z_i)$ converging to 
a point in $l$, then by connectedness
there is
one intersecting a compact middle region $-$ this would force
an intersection of $F$ with $l \times \rrrr$ which is impossible.

This finishes the proof of case 1.
In this case we proved there is a {\underline {line leaf}} $l$
of $\Theta(L)$ with $l \subset \partial \Theta(F)$ and $F$ escapes
up
as $\Theta(F)$ approaches $l$.

\begin{figure}
\centeredepsfbox{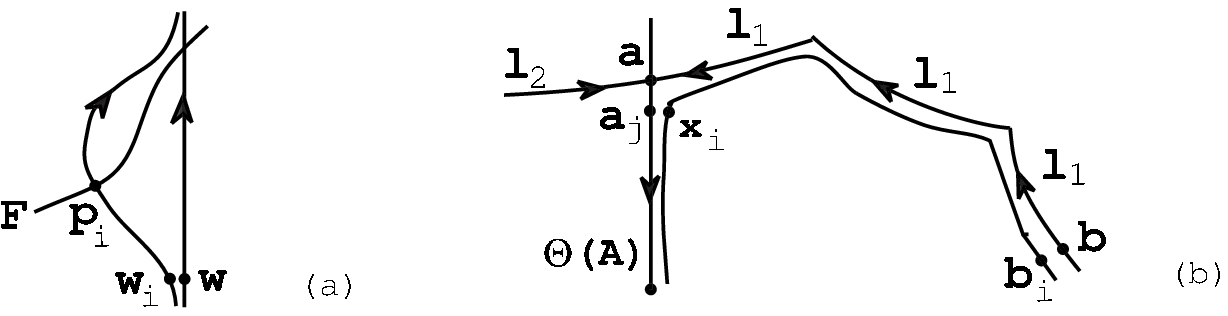}
\caption{
a. A strangling neck is being forced,
b. A slice in a leaf of $\oos$ or $\oou$.
$x_i = \Theta(z_i)$.}
\label{neck}
\end{figure}

\vskip .1in
\noindent
{\underline {Case 2}} $-$ $w$ is in a blown up leaf, but
$F$ does not intersect a lift annulus in $\ws(w)$.


Refer to the setup above.
As before we first show we can assume $w_i$ are not in $\ws(w)$.
Otherwise, up to subsequence assume all $w_i$ are in
$\ws(w)$.
Since $F$ does not intersect lift annuli in $\ws(w)$, then
$w_i$ are all in prongs of $\ws(w)$. Up to subsequence
we can assume they
are all 
in the same prong $C$ of $\ws(w)$ which has boundary an orbit
$\gamma$ of $\wwp$.
It follows that $w$ is in $\gamma$.
All the orbits in $C$ are forward asymptotic to $\gamma$,
even in the blown up situation.
The strangling necks analysis of case 1 shows that $F$
will be forced to intersect $\wwr(w)$.
This cannot occur.

Hence assume all $v_i$ are in a sector of $\oos(v)$ bounded
by a line leaf $l$. Let $L$ be $l \times \rrrr$.
Let $q$ be a point in $l$ and choose $x, \alpha,
x_i$ and $\alpha_i$ as in the proof of case 1.
Choose a small disc $D$ which is transverse to $\wwp$ and
has $\alpha$ in its interior. For $i$ big enough  then
$D$ intersects lift annuli  only in $\ws(x)$.
This is because the union of the blown annuli forms a compact set
in $M$, so  either $\alpha$ intersects a lift annulus,
in which case there is no other lift annulus nearby
or $D$ is entirely disjoint from lift annuli.
From now on the arguments of case 1 apply perfectly.
This shows that 
$\Theta(L)$ is contained in $\overline {\Theta(F)}$,
it is disjoint from $\Theta(F)$ and so it is $\partial \Theta(F)$
and $F$ escapes up as it approaches $L$.
This finishes the proof of case 2.


Now we need to understand what happens when $F$ intersects
a lift annulus in general. We separate that in a special case.
We need the following facts before addressing this case.
A lift annulus $W$ through $b$ is contained in $\ws(b)$ and
$\wu(b)$ so there is not stable/unstable flow directions
in $W$. However there are still such directions
in $\partial W$, because one attracts nearby orbits
of $\wwp$ in $W$ and the other one repels nearby orbits
in $W$. In
this generalized sense the first one is stable and the
second one is unstable. In this sense
if $a$ is in an endpoint of a blown segment,
then all local components of $\oos(a) - a, \oou(a) - a$ near $a$
are either generalized stable or unstable. With
this understanding there is an even number of such
components and they alternate between generalized stable and unstable.
Some local components of $\oos(a) - a$ are also local
components of $\oou(a) -a$ if they are blown segments.
One key thing to remember is that generalized 
stable and unstable alternate.

\vskip .1in
\noindent
{\underline {Case 3}} $-$
Suppose that $F$ intersects some lift annulus $A$ contained
in $\ws(w)$.

Then 
$F$ does not intersect both boundary orbits of $A$.
Otherwise collapse $\pi(A)$ to a single orbit, still
keeping $\Phi$ transverse to $\fol$. 
This contradicts that $\Phi$ is blow down minimal with
respect to $\fol$.
Hence either $F \cap A$ is contained in
the interior of $A$
or it intersects only one boundary leaf.

Assume without loss of generality that $F$ escapes 
{\underline {up}} in one
direction.
This defines an orbit $\gamma$ of $\wwp$ with
$a = \Theta(\gamma)$ in $\partial \Theta(F)$.
The orbit $\gamma$ has to be in the boundary of the lift
annulus $A$. This is because an interior orbit is asymptotic
to both boundary orbits and hence would intersect $F$.
We now look at the picture in $\oo$. Consider the stable leaf
$\oos(a)$. Notice that $\Theta(F)$ intersects $\Theta(A)$. 
From the point
of view of $\gamma$, orbits in $A$ move away from $\gamma$ in
future time, that is $A$ is an unstable direction
from $\gamma$. This means that $\Theta(A)$ is
generalized unstable as seen from $a$.
It follows that 
there are two generalized stable sides of $\oos(a)$
one on each side of $\Theta(A)$
which are the closest to $\Theta(A)$.
Choose one side, start at $a$ and follow along $\oos(a)$ if needed
through blown segments and eventually into a prong in 
$\oos(a)$ so as to
produce a piece of a line leaf of $\oos(a)$ in that direction.
This path is regular on the side associated
to $\Theta(A)$ and defines a half leaf
$l_1$ of $\oos(a)$. Similarly define
$l_2$ in the other direction, see fig. \ref{neck}, b.
Let $l$ be the union of $l_1$ and $l_2$. Then $l$ is a slice leaf
of $\oos(a)$ but is not a line leaf since $\Theta(A)$ is in $\oos(a)$
and is not in $l$.

\vskip .1in
\noindent
{\underline {Claim}} $-$ $l$ is contained in $\partial \Theta(F)$
and $F$ escapes positively as $\Theta(F)$ approaches $l$.

Let $b$ in $l_1$ with $b$ not in blown segment, that is,
$b$ in a prong.
Choose  $b_i$ in $\oou(b)$, with 
$b_i \rightarrow b$ and in that component of $\oo - l$.
Let $D$ be an embedded disc in $\mi$ which is transverse to $\wwp$ and
projects to $\oo$ to a neighborhood of the arc $\xi$ in $l_1$ 
from $a$ to $b$.
Let $y_i$ in $D$ with 
$\Theta(y_i) = b_i$, $y_i \rightarrow y$ with $\Theta(y) = b$.
Assume that $y$ is not in $\pi^{-1}(U)$.
Choose $b$ so that it is not in the unstable leaf of one singular orbit,
hence $\wu(y)$ does not contain lift annuli.
In addition choose $y_i$ so that $\ws(y_i)$ does not contain lift annuli either.

Choose points $u_j$ in $F \cap A$ so that 
$\Theta(u_j) = a_j$ converges to $a$.
For each $j$ the set $\Theta(F)$ contains a small neighborhood
$V_j$ of $\Theta(u_j)$ with $V_j$ converging to $a$ when
$j$ converges to infinity. 
Notice that $a$ is not in $V_j$ as $a$ is not in $\Theta(F)$.
The leaves $\oos(b_i)$ are getting closer 
and closer to $l_1$ and $\Theta(A)$.
For $j$ fixed there is $i$ big enough so that
$\oos(b_i)$ intersects $V_j$. Let

$$z_i \ \in \ F \cap \ws(y_i) \ \ {\rm with} \ \ 
\Theta(z_i) \in V_j$$

\noindent
here $i$ depends on $j$.
Let $z_i = \wwp_{t_i}(r_i)$ with $r_i$ in $D$. 
By choosing $j$ and $i$ converging
to infinity we get that 
$\Theta(z_i)$ converges to $a$ and we can ensure that the
arc of  $D \cap \ws(y_i)$  between $r_i$ and $y_i$ is converging
to the arc $\eta$ of $\ws(\gamma) \cap D$ with $\Theta(\eta) = \xi$.
We can also choose $V_j$ small enough so that $t_i$ converges to
$+\infty$.

The orbits $\wwr(y_i), \wwr(r_i)$ are very close 
in the forward direction as long
as they are outside $\pi^{-1}(U)$. Since $\ws(y_i)$ does not
contain lift annuli then for times $s$
converging to infinity $\wwp_s(y_i)$ is not in $\pi^{-1}(U)$.
Consider the flow band $C$ in $\ws(y_i)$
between $\wwr(r_i)$ and $\wwr(y_i)$. 
The leaf $F$ intersects $\wwr(r_i)$ in $\wwp_{t_i}(r_i)$ with 
$t_i$ converging to infinity. 
Then an analysis exactly as in case 1 considering strangling necks
and the arcs $B$ in that proof, shows that $F \cap \ws(y_i)$ 
cannot escape up 
before intersecting $\wwr(y_i)$. 

Suppose that $F$ escapes down before intersecting $\wwr(y_i)$.
We show that this is impossible.
Since $F \cap \ws(y_i)$ has points $z_i$ in the forward
direction from $D$ and points in the backwards direction
from $D$ it follows that
$F \cap \ws(y_i)$ must intersect  $D$ in at least a
point $q_i$. Up to subsequence we may assume that $q_i$
converges to $q$ in $\ws(y)$.
This will be an iterative process. 
Let $u_1 = q$.
It is crucial
to notice that in the flow band of $\ws(y)$ between $\wwr(y)$ and
$\gamma$ the flow lines tend to go closer to $\gamma$,
that is, either they project to closed orbits freely
homotopic to $\pi(\gamma)$ or they are asymptotic
to one of these orbits moving closer to $\gamma$.
We now consider the component
of $F \cap \ws(y)$ containing $u_1$ and follow
it towards $\gamma$. 
This component does not intersect $\gamma$ and by
the above it can only escape down in $\ws(y)$.
As it escapes down it produces points $c_i$ in $\ws(r_i)$ 
and as before produces points $c'_i$ in $D$, which
up to subsequence converge to $c$ in $D \cap F$.
By construction $c$ is not $u_1$ and its
orbit is closer to $\gamma$.
Let $u_2 = c$. 
We can iterate this process. 
Notice the $u_i$ cannot accumulate in $D$, or else
all the corresponding points of $F$ are in a compact
set of $\mi$.
On the other hand the process does not terminate.
This produces a contradiction.

The contradiction shows that in fact  
the arc
$\Theta(C)$ is in $\Theta(F)$ which implies that 
$\xi = \Theta(\eta)$ is
contained in $\overline {\Theta(F)}$. As the time to hit $F$
from $D$ grows with $i$, this shows that $\Theta(F)$ does
not intersect $\xi$ and hence $\xi$ is contained in $\partial \Theta(F)$.
As $b$ is arbitrary this shows that $l \subset \partial \Theta(F)$ and
$F$ escapes up as $\Theta(F)$ approaches $l$.
This finishes the analysis of case 3.

\vskip .1in
\noindent
{\underline {Case 4}} $-$ $w$ is in a blown up stable leaf and
$F$ intersects some lift annulus $A$ in $\ws(w)$.

The difference from case 3 is that in case 3 we obtained a slice
boundary $l$ of $\Theta(F)$ $-$ but in our situation 
we do not yet know
if it contains $\Theta(w)$ and  whether it is a stable or
unstable. Here we prove it is a stable slice and it contains
$\Theta(w)$.

Recall the setup: $v = \Theta(w)$ is in $\partial \Theta(F)$
and there are $v_i$ in $\Theta(F)$ with $v_i$ converging  to 
$v$ and with $p_i$ in $(v_i \times \rrrr) \cap F$.
Also $p_i = \wwp_{t_i}(w_i)$ with $w_i$ converging to $w$
in $\mi$ and $t_i$ converging to infinity.
Let $\xi$ be the blown segment $\Theta(A)$.

The analysis of case 3 shows that $\Theta(F)$ contains 
the interior of $\Theta(A)$.
Suppose first that 
$v$ is in $\xi$. Then $v$ is in the boundary of 
$\xi$ and by case 3 again $F$ escapes up
or down when $\Theta(F)$ approaches a slice
which contains $v$.
If it escapes up, then the slice is a stable slice
and we obtain the desired result in this case.
We now show that $F$ does not escape down. 
Let $l$ be the unstable slice in $\partial \Theta(F)$
associated to this. Then $l$ cuts in half a small
disk neighborhood of $v$ in $\oo$.
The set $\Theta(F)$ intersects only one component
of the complement, the one which intersects $\xi$.
As $F$ escapes down when $\Theta(F)$ approaches $l$, then
for all points in $\Theta(F)$ near enough $v$ the corresponding
point in $F$ is flow backwards from $D$.
This contradicts the
fact that $t_i$ is converging to infinity.
Therefore $F$ cannot escape down as it approaches $l$.

We can now assume that $v$ is not in $\xi$.
By changing $\xi$ if necessary assume that $\xi$ is the
blown segment in $\oos(v)$  intersected by $\Theta(F)$ which is
closest to $v$.
Let $z$ be the endpoint of $\xi$ separating the rest of $\xi$
from $v$ in $\oos(v)$.

We first show that $z$ is not in $\Theta(F)$. Suppose
that is not the case and let $b$ the intersection point
of 
$z \times \rrrr$ and $F$.
Since $\xi$ is the last blown segment of
$\oos(v)$ between $\xi$ and $v$ intersected by
$\Theta(F)$ and $\Theta(F)$ contains an open neighborhood
of $z$, it follows that $v$ is in a prong $B$ of
$\oos(v)$ with endpoint $z$. 
Let $\tau$ be the component of $F \cap \ws(b)$ containing
$b$. Since $F$ does not intersect
$v \times \rrrr$ then it escapes. As the region between
$\wwr(b)$  and $z \times \rrrr$ (should it be
$v \times \rrrr$ instead of $z \times \rrrr$? previous
report???)
is a prong,
then $F$ cannot escape up. As seen in the arguments
for case 3, $F$ cannot escape down either.
This shows that $z$ cannot be in $\Theta(F)$.

%
%

It follows that $F$ escapes either up or down as $\Theta(F)$
approaches $z$. 
Suppose first that it escapes up.
Then we are in the situation
of case 3 and we produce a stable slice $l$ in
$\partial \Theta(F)$ with $F$ going up as $\Theta(F)$ approaches
$l$. 
If $v$ is not in $l$ then $l$ separates $v$ from
$\Theta(F)$. This contradicts $v_i$ in $\Theta(F)$
with $v_i$ converging to $v$. Hence $v$ is in $l$
with $F$ escaping up as $\Theta(F)$ approaches $l$.
This is exactly what we want finishing the analysis
in this case.

The last situation is $F$ escaping down in $A$ as $\Theta(F)$
approaches $z$. By case 3 there is a slice leaf $l$ in
$\oou(z)$ with $l$ contained in $\partial \Theta(F)$ and
$F$ escaping down as $\Theta(F)$ approaches $l$.
We want to show that this case cannot happen.
Notice that the blown segments 
of $\oos(z)$ are exactly the
same as the blown segments of $\oou(z)$. The sets $\oos(z), \oou(z)$
differ exactly in the prongs and as they go around the collection
of blown segments. 
The collection of all prongs in $\oos(z), \oou(z)$
also alternates between stable and unstable as it goes around
the union of the blown segments.

Suppose first that $v$ is in $l$.
This contradicts $F$ escaping down and $t_i \rightarrow \infty$.
Finally suppose that $v$ is not in $l$.
We claim that in this case $l$ separates $v$ from $\Theta(F)$.
Let $\alpha$ be the path in $\oos(v)$ from $z$ to $v$.
If $\alpha$ only intersects $l$ in $z$, then the separation
property follows because
$l_1$ and $l_2$ contain the local components
of $\oos(z) \cup \oou(z) - z$ which are closest to $\Theta(A)$.
This was part of the construction
of $l$ in case 3.
Here the $\xi$ is generalized stable at $z$
and $l_1, l_2$ are generalized unstable at $p$. The path from $z$ to 
$v$ in $\oos(v)$ cannot start in $\xi$ or $l_1$ or $l_2$, hence
$l$ separates $\Theta(F)$ from $v$.


If on the other hand $\alpha \cap l = \delta$ is not
a single point, then it is a union of blown segments.
Let $u$ be the other endpoint of $\delta$. By regularity
of $l_1$ and $l_2$ on the $\Theta(F)$ side it follows that each
blown up segment in $\delta$ has flow direction away from
$z$. Hence $\delta$ is generalized stable at $u$.
Therefore the closest component of
$\oos(u) \cup \oou(u) - u$ on the
$\Theta(F)$ side is generalized unstable and that is contained
in $l$.
In this case it also follows that $l$ separates $v$ from
$\Theta(F)$. As seen before this is a contradiction.

This finishes the proof of proposition \ref{bounda}
\end{proof}

This has an important consequence that will be used extensively
in this article.

\begin{proposition}{}{}
Let $F$ in $\fn$ and $L$ in $\wls$ or $\wlu$.  Then the intersection
$F \cap L$ is connected.
\label{conn}
\end{proposition}

\begin{proof}{}
By transversality of $\fol$ and
$\Phi$, the intersection $C = \Theta(F) \cap \Theta(L)$
is open in $\Theta(L)$.
Suppose there are 2 disjoint components $A, B$ of $C$.
Then there is $v$ in $\partial A$ with $v$ separating
$A$ from $B$. There are $v_i$ in $A$ with $v_i$ converging
to $v$. By the previous proposition $F$ escapes up or down
in $A \times \rrrr$ as $\Theta(F)$ approaches $v$.
Assume wlog that $F$ escapes up. Then there is a slice
leaf $l$ of $\oos(v)$ with $l \subset \partial \Theta(F)$
and $F$ escapes up as $\Theta(F)$ approaches $l$.
Since $l$ and $\Theta(F)$ are disjoint then $B$ is
disjoint from $l$. In addition $v$ separates
$B$ from $A$ in $\Theta(L)$. It follows from the construction
of the slice $l$ as being the closest to $A$, that
$l$ separates $A$ from $B$. Hence $\Theta(F)$ cannot
intersect $B$, contrary to assumption.
This finishes the proof.
\end{proof}

As promised, we now prove that $\fol$ being almost 
transverse to a pseudo-Anosov flow implies that $\fol$
is Reebless.

\begin{proposition}{}{}
Let $\fol$ be a foliation almost transverse to a pseudo-Anosov
flow  $\Phi_1$ and transverse to a corresponding almost pseudo-Anosov
flow $\Phi$. Then $\fol$ is Reebless.
\end{proposition}

\begin{proof}{}
Suppose that $\fol$ is not Reebless and consider a Reeb component
which is a solid torus $V$ bounded by a torus $T$.
Assume that the flow $\Phi$ is incoming along $T$.

Recall that  there are some singular orbits of $\Phi_1$
which blow up into a collection of flow annuli of
$\Phi$. Suppose that $V$ intersects one of these annuli
$A$. Then since $\Phi$ is incoming along $T$, the torus
$T$ cannot intersect the closed orbits in $\partial A$.
Hence it intersects the interior of $A$, say in
a point $p$ and the forward orbit of $p$ will
limit in a closed orbit which is contained in
the interior of $V$.

If on the other hand $V$ does not intersect these blown
annuli then the blow down operation does not affect
the flow in $V$. That means we can assume that $\Phi_1$ is
equal to $\Phi$ in $V$. Since orbits of $\Phi_1$
are trapped inside $V$
once they enter $V$, the shadow lemma for pseudo-Anosov
flows \cite{Han,Man,Mo1}, shows that 
there is also a periodic orbit of $\Phi_1$ (and hence
also of $\Phi$) in $V - T$.
Notice that the shadow lemma is for pseudo-Anosov flows
and not for almost pseudo-Anosov flows and that is why
we split the analysis into 2 cases.

In any case there is a closed orbit $\gamma$ of $\Phi$ contained
in the interior of $V$.
Consider the generalized stable/unstable local
leaves at $\gamma$. Since $\Phi$ is incoming along $T$,
the generalized unstable leaves have to be contained
in $V$. We eventually obtain that a whole half leaf
of $W^u(\gamma)$ is contained in $V$. 
A lift $\widetilde V$ to $\mi$ is homeomorphic to $D^2 \times
\rrrr$, because closed orbits of $\Phi$ are not null homotopic.
The procedure above produces a half leaf of $\wu(\widetilde \gamma)$
contained in $\widetilde V$. This contradicts the
fact that $\wu(\widetilde \gamma)$ is properly 
embedded \cite{Ga-Oe}.
This shows that $\fol$ is Reebless.
\end{proof}

\section{Asymptotic properties in leaves of the foliation}
\label{asym}

Let $\Phi$ be an almost pseudo-Anosov flow transverse
to a foliation $\fol$ with hyperbolic leaves. Let
$\Lambda^s, \Lambda^u$ be the singular foliations
of $\Phi$. Given leaf $F$ of $\fn$ let $\wlsf, \wluf$
be the induced one dimensional singular foliations in $F$.
In this section we study asymptotic properties of 
rays in $\wlsf$.
First we mention a result of Thurston \cite{Th5,Th7} concerning
contracting directions, which for convenience we
state for 3-manifolds:

\begin{theorem}{}{(Thurston)} 
Let $\fol$ be a codimension one foliation with hyperbolic
leaves in $M^3$ closed. Then for every $x$ in any leaf
$F$ of $\fn$ 
and every $\epsilon > 0$ there is a dense set of
geodesic rays 
of $F$ starting at $x$ such that:
for any such ray $r$
there is a transversal $\beta$ to $\fn$ starting at $x$ so
that any leaf $L$ intersecting $\beta$ and any
$y$ in $r$, then the distance between $y$ and $L$ is
less than $\epsilon$.
If there is not a holonomy invariant transverse measure
whose support contains $\pi(F)$ then one can show
that the directions are actually contracting, that is:
if $y$ escapes in $r$ then the distance between $y$ and
$L$ converges to $0$.
Finally if $\pi(F)$ is not closed one can choose the 
$\beta$ above to have $x$ in the interior.
\end{theorem}

There is a carefully written published version of this
result in \cite{Ca-Du}.
The directions above where distance to nearby $L$ goes to
$0$ are called contracting directions. The other ones
where distance is bounded by $\epsilon$
are called $\epsilon$ non expanding directions.

The goal of this section is to show that given a leaf
$L$ of $\fn$ and a ray $l$ of ${\widetilde \Lambda}^s_L$, 
then $l$ converges to a single point in $\pin L$.
We first analyse 
the non $\rrrr$-covered case.
The proof is very involved and is done by way 
of contradiction.
Later we deal
with the $\rrrr$-covered situation.
This result is a natural extension of a result by
Levitt  \cite{Le} who proved that if
$\gal$ is a foliation with prong singularities
in a closed hyperbolic surface, then in the universal
cover, an arbitrary ray converges to a single point in
the circle at infinity.
The situation for non compact leaves of foliations
is much more delicate.

\begin{proposition}{}{}
Let $\Phi$ be an almost pseudo-Anosov flow transverse to
a foliation $\fol$ in $M^3$ closed and $\fol$ with hyperbolic
leaves. 
Suppose that $\fol$ is not $\rrrr$-covered.
Given a leaf $L$ of $\fn$ and an arbitrary ray $l$ 
in a leaf of ${\widetilde \Lambda}^s_L$ or
${\widetilde \Lambda}^u_L$ 
then $l$ limits to a single
point in $\pin L$. The limit depends on the ray $l$.
\label{as1}
\end{proposition}

\begin{proof}{}
We do the proof for 
${\widetilde \Lambda}^s_L$.
Let $\epsilon$ positive so that if $p$ in $\mi$ is
less than $\epsilon$ from a leaf $F$ of $\fn$, then
the flow line through $p$ intersects $F$ less than
$2 \epsilon$ away from $p$.
Let $l$ be a ray in ${\widetilde \Lambda}^s_L$.
Because $\fol$ and $\Phi$ are transverse,
$L$ is properly embedded in $\mi$ and leaves of $\wls$ are
properly embedded, it follows that $l$ is a properly embedded
ray in $L$. Therefore it can only limit in 
$\pin L$. 

Suppose by way of contradiction that $l$ 
limits on 2 distinct points $a_0, b_0$ in $\pin L$.
Fix $p$ a basepoint in $L$. 
Since
$l$ limits in $a_0, b_0$, there are compact arcs $l_i$ of
$l$ with endpoints which converge to $a_0, b_0$ respectively
in
$L \cup \pin L$ and so that the distance from $l_i$ to $p$ in
$L$ converges to infinity. Also we can assume that the
$l_i$ converges to a segment $v$ in $\pin L$, where $v$ connects
$a_0, b_0$. This is in the Hausdorff topology of closed sets
in $L \cup \pin L$,  which is a closed disk.

The key idea is to bring this situation to a compact part of $\mi$.
Choose a sequence $p_i$ at bounded distance
from points in $l_{k_i}$ so that that $p_i$ converges
to a point $a$ in the interior of $v$. 
The bound depends on the sequence.
Up to subquence
assume that there are convering translations $g_i$ in 
$\pi_1(M)$ and a point $p_0$ in $\mi$ 
so that $g_i(p_i)$ converges to 
$p_0$ in $\mi$.

We claim that the set of possible limits $p_0$ 
obtained as above projects
to a sublamination of $\fol$.
Clearly if $g_i(p_i)$ converges to $p_0$ and $q$ is in the same 
leaf $L_0$ of $\fn$ as $p_0$, then the distance 
from $p_0$ to $q$ is finite 
and there are $q_i$ in
$L$ with $d_L(q_i, p_i)$ bounded and
$g_i(q_i)$ converging to $q$. Also $q_i$ converges to $a$ in $\pin L$.
In addition if a sequence of such limits $c_j$ converges
to $c_0$ then a diagonal process shows that $c_0$ is also
obtained as a single limit.
This proves the claim.
Choose a minimal sublamination $\ccl$.

A leaf $F$ of $\fn$ is isometric to the hyperbolic plane. A 
{\em wedge} $W$ in $F$ with corner $b$ and ideal set 
an interval $B \subset \pin F$
is the union of the rays in $F$ from $b$ with ideal point in $B$.
The angle of the wedge is the angle that the boundary rays of
$W$ make at $b$.
For any such sequence $p_i$ as above, then
the visual angle at $p_i$ subintended by the arc $v$ in $\pin L$ grows
to $2 \pi$. Therefore the angle of wedge with corner
$p_i$ and ideal set 
$\pin L - v$ converges to $0$.
This is called the {\em bad wedge}.

Assume up to subsequence that $g_i(p_i)$ is converging
to $p_0$ in a leaf $L_0$ of $\fn$ and that the 
directions of the
bad wedges with corners $g_i(p_i)$ in $g_i(L)$ are converging 
to the direction $r_0$  of $L_0$. Let $c$ be the ideal
point of $r_0$ in $\pin L_0$.

Suppose first that $\pi(L_0)$ is not compact $-$ we shall
see briefly that this is in fact always the case.
Thurston's theorem shows that the set of two sided contracting
directions (or $\epsilon$ non
expanding directions) in $L_0$ is dense
in $\pin L_0$. 
We will use these to transport a lot of the structure
of ${\widetilde \Lambda}^s_{L_0}$ to nearby leaves.
Choose $s_0, s_1$ to be rays in $L_0$ defining
contracting directions
(or $\epsilon$ non expanding directions)
very near $r_0$ so that together they form a small wedge
$W$ in $L_0$ with corner $p_0$.
There is an interval of leaves near $L_0$ so that
any such leaf $V$ is less than $\epsilon$ away from 
$s_0, s_1$. Then a flow line  of $\wwp$ through any point
in $s_0$ or $s_1$ intersects $V$ less than $2 \epsilon$
away. So $s_0$ flows to a curve in $V$, where we can
assume it has geodesic curvature very close to $0$,
if $\epsilon$ is sufficiently small. It is therefore 
a quasigeodesic with a well defined ideal point.
The same happens for $s_1$ and the flow images $u_0, u_1$
of $s_0, s_1$ in $V$ define a generalized wedge $W'$
in $V$.
The ideal points $e_0, e_1$ of $u_0, u_1$  are
close and bound an interval $I$ which is 
almost all of $\pin g_i(L)$.

By construction $g_i(l)$ is a
ray which limits in an interval of $\pin g_i(L)$ which
contains $I$ in its interior if
$i$ is big enough. 
There are then  subarcs $\tau_j$ of $g_i(l)$ with endpoints
$a_j, b_j$ in $u_0, u_1$ respectively so that $a_j$ converges
to $e_0$ and $b_j$ converges to $e_1$ and $\tau_j$ converges
to $I$, see fig. \ref{inde}.
Here $i$ is fixed and $j$ varies. 
Since $a_j, b_j$ are in $u_0, u_1$ then they flow 
(by $\wwp$) to points
in $L_0$. The images in $L_0$ are in the same leaf of  $\wls$.
By proposition \ref{conn} these images 
are in the same leaf of
$\widetilde \Lambda^s_{L_0}$.
Hence the whole segment $\tau_j$   flows into
$L_0$.

\begin{figure}
\centeredepsfbox{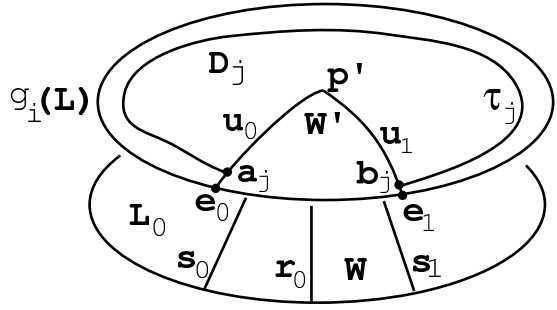}
\caption{
Transporting the
structure between leaves $g_i(L)$ and $L_0$.}
\label{inde}
\end{figure}

The point $p_0$ flows into $p'$ in $g_i(L)$ under the flow.
The arc $\tau_j$ together with 
subarcs or $u_0, u_1$ from $a_j, b_j$ to $p'$ 
bound a disc $D_j$ in $g_i(L)$. The arguments above show
that the boundary of $D_j$ flows into $L_0$ producing
a curve in $L_0$ bounding 
a disc $B_j$. The segments of $\wwp$ connecting points
in $\partial D_j$ to points in
$\partial B_j$ produce an annulus $C_j$. Then $D_j \cup C_j
\cup B_j$ is an embedded sphere in $\mi$ and hence
bounds an embedded ball. 
Since orbits of $\wwp$ are properly embedded in $\mi$,
it follows that all orbits of $\wwp$ intersecting
$D_j$ will also intersect $B_j$. Hence there is product
flow in this ball. Since this is true for all $j$ then
the union of the $D_j$ flows into $L_0$.
The union of the $D_j$ is the closure of $g_i(L) - W'$.
The image is contained in the closure of $L_0 - W$
in $L_0$ $-$ call the closure $J$. 

We claim
that the image is in fact $J$.
All the $\tau_j$ are in the same leaf of $\wls$ and hence
all their flow images  in $L_0$ also are. 
Since rays of $\widetilde \Lambda^s_{L_0}$ are
properly embedded in $L_0$ then 
when $j$ converges to infinity the images of $\tau_j$
in $L_0$ escape
compact sets. This shows the claim.
Therefore the flow produces a homeomorphism between
the closure of $L_0 - W$ and the closure of 
$g_i(L) - W'$. Clearly the same is true for any leaf
in the interval
associated to the contracting (non expanding) directions
$s_0, s_1$. 
In particular we have the following conclusions:

\vskip .1in
\noindent
{\bf {Conclusion}} $-$ In any limit leaf $L_0$ with 
a limit direction $r_0$ of bad wedges the following
happens: Let $c$ be the ideal point of $r_0$ and
$A$ a closed interval of $\pin L_0 - \{ c \}$.
Then there is a leaf $l$ of $\widetilde {\Lambda}^s_{L_0}$
with compact subsegments $l_i$ so that the endpoints
of $l_i$ converge to the endpoints $a, b$ of
$A$ and $l_i$ converges to $A$.
In particular $l_i$ escapes compact sets.
There are also subsegments $v_i$ with both endpoints
converging to $a$ and so that $v_i$ converges to 
sets in $\pin L_0$ which contain $A$.
Finally for sufficiently near leaves 
there is a wedge in $L_0$ which forms
a product flow region with these nearby leaves.
\vskip .1in

To get the second assertion above just follow $l$ beyond the 
endpoint of $l_i$ near $b$ until it returns near $a$
again.
As a preliminary step to obtain proposition \ref{as1} we prove
the following:

\begin{lemma}{}{}
For any limit $g_i(p_i)$ converging to $p_0$,
the distinguished direction of the bad wedge associated
to $g_i(p_i)$ converges to a single direction at $p_0$.
In the second case
this direction varies continuously with the leaves in $
\widetilde {\ccl}$.
\label{single}
\end{lemma}

\begin{proof}{}
Suppose there are subsequences $q_i, p_i$ converging to 
points in (interior) $v$ with $g_i(p_i), h_i(q_i) \rightarrow p_0
\in L_0 \in \fn$,
but the directions of the wedges converge to $r_0, r_1$
distinct geodesic rays in $L_0$.
We will first show that there is an interval of leaves
of $\fn$  so that the flow $\wwp$ is a product flow 
in this region.

Using the limit direction $r_0$ we produce a wedge
$W$ in $L_0$ so that the closure of $L_0 - W$
is part of a product flow region with nearby leaves
of $\fn$.
Using the other limit direction $r_1$
we produce a flow product region associated to
another wedge region $W_*$ disjoint from $W - p_0$.
Together they produce a global product structure
of the flow in a neighborhood of $L_0$.






This shows that there is a neighborhood $N$ of $L_0$ in the
leaf space of $\fn$ so that
the flow is a product flow in $N$.
In particular there is no non Hausdorffness of $\fn$ in
this neighborhood.
This is a very strong property as we shall see below. It 
implies a global product structure of the flow.

Notice that the structure of $\widetilde \Lambda^s_{g_i(L)}$ in
$g_i(L) - W'$ flows over to $L_0$. In particular there are 
many rays of $\widetilde \Lambda^s_{L_0}$ which do not
have a single limit in $\pin L_0$. This implies that
$\pi(L_0)$ is not compact. This is because Levitt \cite{Le}
proved that given any singular foliation with prong singularities
in a closed hyperbolic surface $R$, then 
the rays of the lift to $\widetilde R$ all have unique
limit points in the ideal boundary.
This shows that the minimal lamination $\ccl$ is not a compact leaf
and hence it has no compact leaves.

\vskip .1in
Consider the neighborhood $N$ as above.
Consider the translates $g(N)$ where $g$ runs through all elements
of the fundamental group. 
Let $P$ be the component of the union containing $N$. 
It is easy to see that the
set $P$ is precisely invariant: if $g$ is in $\pi_1(M)$
and $g(P)$ intersects $P$ then $g(P)$ is equal to $P$.
In addition $\fol$ restricted to $P$ has leaf space
homeomorphic to $\rrrr$ because of the product flow
property.
We are assuming that $N$ is open.

Suppose first that $P$ is not all  of $\mi$, hence
$\partial P$ is a non empty collection of leaves of $\fn$.
Let $C$ be the projection of $P$ to $M$.
Then $C$ is open, saturated by leaves of $\fol$.
Notice that $g(P)$ does not intersect $\partial P$ for
any $g$ in $\pi_1(M)$ for otherwise $g(P)$ intersects
$P$ and so $g(P) = P$. It follows that $\pi(\partial P)$
is disjoint from $C$ hence $C$ is a proper open, foliated
subset of $M$.

Dippolito \cite{Di} developed a theory of such open, saturated
subsets. Let $\overline C$ be the metric completion of $C$.
There is an induced foliation in $\overline C$, which
we will also denote by $\fol$. Then

$$\overline C \ \ = \ \ V \ \cup \ \bigcup_1^n V_i$$

\noindent
where $V$ is compact and may be all of $\overline C$.
Each nonempty $V_i$ is an $I$-bundle over a non compact
surface with boundary, so that $\fol$ is a foliation
transverse to the $I$-fibers. Each component of
the intersection $\partial V_i \cap V$ is an annulus
(or Moebius band) with induced foliation transverse
to the $I$ fibers.
In our situation with $\Phi$ transverse to the flow,
if $V$ is not $\overline C$, we can choose $V$
big enough so that the flow is transverse to $\fol$ 
in each $V_i$ and induces an $I$-fibration there.

Let $\widetilde R$ be a component of $\partial P$ and
$R$ the projection to $M$, so $R$ is component
of $\partial C$.
Parametrize the leaves of $\fn$ in $P$ as
$F_t, 0 < t < 1$ with $t$ increasing with
flow direction. A leaf in the boundary of $P$ which
is the limit of leaves in $P$ which are limiting from
the positive side above has to be the
limit of $F_t$ as $t$ goes to $0$: 
Suppose that $S$ is in the boundary of $P$ and there
are $x_i$ in $F_{t_i}$ with $t_i$ converging
to $t_0 > 0$ 
and $x_i$ converging to $x$ in $S$.
Then $S$ and $F_{t_0}$ are not separated from
each other.
For $i$ big enough the flow line
through $x_i$ will intersect $S$ 
and therefore this flow line will not intersect
$F_{t_0}$. This contradicts the fact
that  $F_{t_o}$ and $F_{t_i}$ 
have a flow product structure.

Suppose then that $\widetilde R$ is a limit of $F_t$
where $t$ converges to $0$.
Suppose first that $R$ is compact. 
Suppose there are $t_i$ converging to $0$ so that
$F_{t_i}$ are in $\widetilde {\ccl}$. Then since
$\ccl$ is a closed subset of $M$ it follows that
$\widetilde R$ is in $\widetilde {\ccl}$ and so
$R$ is in $\ccl$. But $R$ is closed, contradicting
the fact that $\ccl$ has no closed leaves.
There is then $a > 0$ which is the smallest $a$ 
so that $F_a$ is in $\widetilde {\ccl}$ $-$ notice
that $\widetilde {\ccl}$ has leaves in $P$.
For any $g$ in $\pi_1(R)$ then $g(N) \cap N$ is
not empty hence $g(N) = N$. It follows that
$g(F_a) = F_b$ for
some $b$. If $b$ is not $a$ then by taking
$g^{-1}$ if necessary we may assume that $b < a$.
But as $F_b$ is in $\widetilde {\ccl}$, this contradicts
the definition of $a$. Hence $g(F_a) = F_a$ for
any $g$ in $\pi_1(R)$. 
This implies that $\pi(F_a)$ is a closed surface,
again contradiction.

We conclude that $R$ is not compact, hence it eventually
enters some $V_i$ (the point here is that $V$ is
not $\overline C$). 
The flow restricted to any component of $\partial V_i
\cap \overline C$ goes from one component to the
other in the annulus. This implies that all $\pi(F_t)$ 
intersect this annulus. 
There is then a leaf $B$ of $\ccl$ which enters
$V_i$. Going deeper and deeper in this non compact $I$-bundle
will produce a limit point which is not in $C$.
This shows the very important fact that $\ccl$ is
not contained in $C$ and therefore

$${\cal E} \ = \ \ccl  \cap   (M - C) \ \not  =  \ \emptyset$$

\noindent
In addition ${\cal E}$ is not equal to $\ccl$ 
since 
$\ccl$ has leaves in $C$
and
$(M - C)$ is closed.
Hence ${\cal E}$ is a non trivial, proper sublamination
of  $\ccl$. This contradicts the fact that $\ccl$ is
a minimal lamination.

This shows that the assumption $P \not = \mi$ is
impossible. Hence $P = \mi$, which implies the
flow $\wwp$ produces a global product picture
of $\fn$ and in particular $\fol$ is $\rrrr$-covered,
contrary to assumption.

This shows 
the limits of the bad wedges are unique directions
in the limit leaves. It also shows that they vary
continuously from leaf to leaf, for otherwise one
obtains bad wedges in very near leaves which have
definitely separated directions. The same proof above
then applies.
This finishes the proof of lemma \ref{single}.
\end{proof}

\vskip .1in
\noindent
{\underline {Continuation of the proof of proposition \ref{as1}}}


By the previous lemma
 we know that limit directions of
bad wedges are unique and they vary continuously in leaves
of $\widetilde {\ccl}$.
These unique directions are distinguished in their
respective leaves.

We first show that any complementary region of $\ccl$
(if any)
is an $I$-bundle with a product flow.

Lift to a double cover if necessary to assume that $M$
is orientable. Assume this is the original foliation
$\fol$, flow $\Phi$, etc..
Let $Z$ be a leaf of $\widetilde {\cal L}$. Since $Z$
has a distinguished ideal point, then the fundamental
group of $\pi(Z)$ can be at most ${\bf Z}$. Since there
is a transverse flow and $M$ is orientable this implies
that $\pi(Z)$ is either a plane or an annulus.

Let $U$ be a complementary region of $\ccl$ with
boundary leaves $R_1, R_2, R_3$, etc..
As explained before the completion of $U$ has
a compact thick part and the non compact 
arms which are in thin, $I$-bundle
regions. 
Suppose first that $R_1$ is a plane.
There is a big disk $D$ so that $R_1 - D$ is contained
in the thin arms and flows across $U$ to another boundary
components of $U$. By connectedness it flows into a single
boundary component $R_2$ of $U$. Then $\partial D$ flows
into a curve $\gamma$ in $R_2$ which is null homotopic
in $M$. 
The flow segments in $M$ produce an annulus
$C$ in the completion of $U$. Since $\fol$ is Reebless
then $\gamma$ bounds a disk $D'$ in $R_2$ and so $R_2$
is a plane. The union $D \cup C \cup D'$ is an
embedded sphere in $M$ which bounds a ball $B$.
Since orbits of $\wwp$ are properly embedded in $\mi$,
it follows that 
the flow has to a product flow
in $B$ as well. 
This shows that flow is a product
in the completion of $U$.

Suppose now that each $R_i$ is an annulus. Let
$F$ be a lift of $R_1$ to $\mi$ with $F$ in the boundary 
of a component $\widetilde U$ of $\pi^{-1}(U)$.
In $R_1$ there are two disjoint open annuli $A_1, A_2$ contained
in the thin arms so that $B = R_1 - (A_1 \cup A_2)$
is a closed annulus in the core.
Then $A_1, A_2$ flow into two annuli leaves
$R_2, R_3$ in the boundary of $U$.
Lifting to $F = \widetilde R_1$ we see leaves
of $\wlsf$ limiting in an interval of $\pin F$ with very
small complement (near the distinguished ideal
point of $F$). This implies they will have
points in the lifts $\widetilde A_1, \widetilde A_2$
of $A_1, A_2$ to $F$.
This shows that $\widetilde A_1, \widetilde A_2$
are in the same leaf of $\fn$.
This implies that $R_2 = R_3$. 
In the same way a half of the infinite strip $\widetilde B$
flows into $\widetilde R_2$. Since $B$ is compact,
then all of $B$ flows into $R_2$.
This implies that the region $U$ is an $I$-bundle.
It is also easy to show that the flow is a product
in this $I$-bundle.

This implies that we can collapse this complementary region
along flow lines to completely eliminate it.
This is because even in the universal cover we are 
eliminating product regions of the flow and the
asymptotic behavior is still preserved in the remaining
regions.
This can be done to all complementary regions and
therefore we can assume there are no complementary
regions, that is $\ccl = \fol$ or that $\fol$ is minimal.

Let $F_1, F_2$ be leaves of $\fn$ which are not
separated from each other. Consider leaves $F$ of 
$\fn$ which are very close to points in both $F_1$
and $F_2$. As stated in the conclusion in the beginning
of the proof of this theorem,
there is a wedge of $F$ which flows into $F_1$
and similarly for $F_2$. Hence there are half planes
$E_1, E_2$ of $F$ which flow into $F_1, F_2$.
As $F_1, F_2$ are not separated this implies that
$E_1, E_2$ are disjoint. Fix a point $w$ in $F$ and
a big enough radius $r$ so that the disk $D$ of radius
$r$ around $w$ intersects both $E_1, E_2$.
Again as seen in the conclusion above
there is an arc $l$ in a leaf of $\wlsf$
so that both endpoints of $l$ are outside $D$ and in $E_1$
and so that $l$ is entirely outside $D$ and as seen
from $p$ the visual measure of $l$ is almost $2 \pi$.
This implies that $l$ intersects $E_2$. Since the
endpoints of $l$ are in $E_1$, which flows to $F_1$,
then proposition \ref{conn} implies that the whole arc $l$ flows
into $F_1$. The points of $l$ in $E_2$ will also flow
to $F_2$. This is a contradiction.

This contradiction 
finishes the proof of proposition \ref{as1}
\end{proof}

Next we analyse the $\rrrr$-covered situation which has interest
on its own:

\begin{theorem}{}{}
Let $\fol$ be an $\rrrr$-covered foliation and $\Phi$ be a pseudo-Anosov flow
almost transverse to $\fol$. Then $\Phi$ is actually transverse to $\fol$.
In addition for any leaf $F$ of $\fn$ and for any ray $l$ in $\wlsf$
it converges to a unique ideal point in $\pin F$.
The limit usually depends on $l$.
\label{as2}
\end{theorem}

\begin{proof}{}
If $\Phi$ is not transverse to $\fol$, let 
$\Phi^*$ be an almost 
pseudo-Anosov flow which is transverse to $\fol$ and
is a blow up of $\Phi$. 
Notice this is not the same notation as in proposition \ref{as1}
$-$ here we prove $\Phi$ is equal to $\Phi^*$.
There is flow annulus $A$ of $\Phi^*$
with closed orbits $\gamma_1, \gamma_2$ in
the boundary, so that $A$ blows down to a single orbit of $\Phi$.

The foliation induced by $\fol$ in $A$ 
has leaves which spiral to at least one
boundary component $-$ which they do not intersect.
Lifting this picture to the universal cover one
obtains an orbit of $\wwp^*$ which does not intersect
every leaf of $\fn$. This means that the flow $\wwp^*$
is not {\em regulating} for $\fn$ \
\cite{Th6, Th7}. 
We also say that $\Phi^*$ does not regulate $\fol$.
In \cite{Fe9} we analysed a similar situation and proved 
the following: if $\Upsilon$ is a pseudo-Anosov flow transverse
to an $\rrrr$-covered foliation and $\Upsilon$ is not regulating,
then $\Upsilon$ is an $\rrrr$-covered
Anosov flow.
The same arguments work with an almost pseudo-Anosov flow
transverse to an $\rrrr$-covered foliation.
This shows that $\Phi^*$ is an $\rrrr$-covered Anosov flow
and has no (topological) singularities.
In particular $\Phi^*$ is equal to $\Phi$, that is
the original flow is already transverse to $\fol$.
This proves the first assertion of the theorem.


Assume by way of contradiction that there
is $L'$ in $\widetilde {\Lambda}^s$ and $l$ in 
$\widetilde {\Lambda}^s_{L'}$ which 
does not converge to a single point in $\pin L'$.
As in the proof of theorem \ref{as1} we construct a minimal
sublamination $\ccl$ of $\fol$ such that:
for every $L$ in $\widetilde {\ccl}$ there is
an ideal point $u$ in $\pin L$ so that for every closed
segment $J$ in $\pin L - \{ u \}$ there is a ray
$l$ of $\widetilde {\Lambda}^s_L$ which has subsegments
limiting to $J$..
As shown in the proof of theorem \ref{as1}, 
$\ccl$ cannot be
a compact leaf.

Suppose first that every leaf of $\fol$ is a plane.
Then Rosenberg \cite{Ros} proved that $M$ is the $3$-dimensional
torus $T^3$. 
This manifold is a Seifert fibered space. In this case
Brittenham \cite{Br1} proved that an essential lamination
is isotopic to one which is either vertical (a union
of Seifert fibers) or horizontal (transverse to the 
fibers). So after isotopy assume $\ccl$ has one
of these types. If $\ccl$ has a vertical leaf $B$, then
geometrically it is a product of the reals with the circle.
Hence it is an Euclidean  leaf and in the universal
cover it has polynomial growth of area.
If $\ccl$ has
a horizontal leaf $B$, then
because the fibration is a product,
there is a projection to a $T^2$ fiber, which distorts
distances by a bounded amount. Again the same growth properties
hold.  But the leaves of $\fol$ are hyperbolic, which is 
a contradiction. We conclude that $M$ cannot be $T^3$.



Let then $F$ in $\widetilde {\ccl}$ with $\pi(F)$ not simply 
connected. Let $g$ in $\pi_1(M)$ non trivial with $g(F) = F$
and $\xi$ be the axis of $g$ in $F$. At least one ideal point
of $\xi$, call it $u$, is not the direction of a fixed 
limit of bad wedges. Then as explained before there is a 
ray $l$ of $\wlsf$ and segments $l_i$ 
of $l$, bounded by $a_i, b_i$ both points
in $\xi$, so that $l_i$ escapes compact sets
and converges to a non trivial segment in $\pin F$.
We may assume that $l_i \cap \xi = \{ a_i, b_i \}$ and also
that all $l_i$ are in the same side of $\xi$.
Let $e_0$ be the translation length of $g$ in $F$.

If the distance from $a_i$ to $b_i$ along $\xi$ is bigger than $e_0$ then
this produces a contradiction as follows:
There is an integer $n$ so that $g^n(a_i)$ is  in the 
open segment $(a_i,b_i)$ of $\xi$ and
and $g^n(b_i)$ is outside of the closed segment $[a_i,b_i]$.
Since the arc $l_i$ only intersects $\xi$ in $a_i, b_i$, then 
$l_i$, together with $[a_i,b_i]$
bounds a closed disk in $F$ and $g^n(a_i)$ is in $(a_i,b_i)$.
But $g^n(b_i)$ is outside and $g^n(l_i)$ is also on this side
of $\xi$, so this produces
a transverse self intersection of $\wlsf$.
If $g^n(l_i)$ is contained in the leaf $v$ which contains $l_i$,
then $g^n(v) = v$ and this produces infinitely many
singularities in $v$, which is impossible. Hence $g^n(l_i)$ is
not in $v$ and the transverse intersection is impossible.
The same arguments deal with the case that $l_i$ intersects
$\xi$ in other points besides $a_i, b_i$.


We conclude that the distance in $\xi$ from $a_i$ to $b_i$ is bounded.
Up to subsequence we may assume there are integers $n_i$
so that 
$g^{n_i}(a_i)$ converges to $a_0$ and $g^{n_i}(b_i)$ converges to $b_0$,
both limits in $\xi$ of course.
Since the lengths of $g^{n_i}(l_i)$ are converging to infinity,
it follows that $a_0, b_0$ are not in the same leaf of $\wlsf$.
By proposition \ref{conn} it follows that $a_0, b_0$ are not
in the same leaf of $\wls$. But for each $i$, the pair
of points $g^{n_i}(a_i), g^{n_i}(b_i)$ is in the same leaf
of $\wls$. This implies that the leaf space of $\wls$ is
not Hausdorff.

First of all this implies that $\Phi$ is regulating for $\fol$,
for otherwise the aforementioned result from \cite{Fe9}
shows that $\Phi$ is an $\rrrr$-covered Anosov flow $-$ in
particular $\wls$ has Hausdorff leaf space.
Also by theorem \ref{theb} the fact that $\wls$ has non Hausdorff
leaf space implies that there are closed orbits $\alpha, \beta$
of $\Phi$ so that $\alpha$ is freely homotopic to the inverse
of $\beta$.
Let $h$ be a covering translation associated to $\alpha$ and
$\widetilde \alpha$, $\widetilde \beta$ lifts of $\alpha, \beta$
to $\mi$ which are left invariant by $h$. Without loss of generality
assume that $h$ acts in $\widetilde \alpha$ sending points forwards.
As $\alpha \cong \beta^{-1}$ this implies that $h$ acts
on $\widetilde \beta$ taking points backwards. But since
both of them intersects all leaves of $\fn$ (by the regulating
property) then
as seen from $\widetilde \alpha$ the translation $h$ acts
increasingly in the leaf space of $\fn$, with opposite
behavior when considering $\widetilde \beta$.
This is a contradiction, which shows that this cannot happen.
This finishes the proof of theorem \ref{as2}.
\end{proof}


Proposition \ref{as1} and theorem \ref{as2} imply the following:

\begin{corollary}{}{}
Let $\Phi$ be an almost pseudo-Anosov flow transverse to
a foliation $\fol$ with hyperbolic leaves in $M^3$ closed.
For any leaf $L$ of $\fn$ and any ray $l$ of $\wlsf$
or $\wluf$, then $l$ converges to a single point in $\pin L$.
\end{corollary}

\vskip .1in
\noindent
{\bf {Remark}} - Group invariance and compactness of $M$ are both
essential here.  For example start with a nicely behaved singular foliation
of $\hh$, so that all rays converge. It could be a foliation by geodesics
or for instance the lift of the stable singular foliation associated
to a suspension. 
Fix a base point $p$. Now rotate the leaves at a distance $d$ of $p$ by
an angle $d$. In this situation all rays limit in all points of $\pin L$,
in fact they spiral indefinitely into it.
Another operation is to fix a ray through $p$ and then distort the rest more
and more one way and the other way. Here we have the leaves getting
closer and closer to segments in $\pin F$ which 
are complementary to the ideal point associated to the
ray.

\section{Properties of leaves of $\wlsf, \wluf$ and their ideal points}

In this section $\Phi$ is an almost pseudo-Anosov flow
transverse to a foliation $\fol$. 
As in the previous section
there is no restriction on $M$ here.
In the previous section we proved that
for any ray $r$ of a leaf of $\wlsf$ or $\wluf$, then
it has a unique ideal point in $\pin F$. The notation for this ideal
point will be $r_{\infty}$.
We now analyse further properties of leaves of $\wlsf$ and their
ideal points.
Analogous results hold for $\wluf$.

First we want to show that if $E$ is a fixed
leaf of $\wls$ (or $\wlu$)
then the ideal points in $\pin F$ of rays of $E \cap F$
vary continuously with $F$.
In order to do that we first
put a topology on the union of ideal boundaries of an interval
of leaves. Let $p$ in $F$ leaf of $\fn$ and $\tau$ a transversal
to $\fn$ with $p$ in the interior.
For any $L$ in $\fn$ intersecting $\tau$, the ideal
boundary is in 1-1 correspondence with the unit
tangent bundle to $L$ at $\tau \cap L$:
ideal points correspond to rays in $L$ starting at $L \cap \tau$.
This is a homeomorphism.
This puts a topology in 

$$ {\cal A}  \ \ = \ \ 
\cup \ \{ \pin L \ \ | \ \ L \cap \tau \not = \emptyset \}$$

\noindent
making it into an annulus homeomorphic to 
$\cup  \ \{ T^1_q \fn, q \in \tau \}$ as a subspace
of the unit tangent bundle of $M$.
This topology in ${\cal A} $ is independent of the choice of transversal
$\tau$. The following definition/result is proved
in \cite{Fe7} or \cite{Cal1}.

\begin{define}{(markers)}{}
Given a foliation $\fol$ by hyperbolic leaves
of $M^3$ closed, then  there is $\epsilon > 0$ so
that:
Let $v$ be a geodesic ray in a leaf $F$ so that 
it is associated to a contracting (or $\epsilon$ non
expanding) direction of $F$.
For any leaf $L$ sufficiently near $F$, then
all the points of $v$ flow into $L$ and define
a curve denoted by $v_L$. 
Then $v_L$ has a unique ideal point denoted by $a_L$. The
union $m$ of the $a_L$ is called a marker and is
a subset of ${\cal A}  = \cup \ \{ \pin L \}$. 
Then $m$ is an embedded curve in ${\cal A} $ in the topology
defined above.
\end{define}

In addition the markers are dense in ${\cal A} $ in the following sense:
Let $z$ in $\pin F$ and $a_i, b_i$ in $\pin F$ which are 
in markers associated to contracting (non expanding) directions
on a fixed side of $F$. Suppose that the sequence of open 
intervals
$(a_i, b_i)$ in $\pin F$ contains $z$ and converges
to $z$ as $i$ converges to infinity.
Let $\alpha_i, \beta_i$ be the markers in that
side of $\pin F$ containing $a_i, b_i$ respectively.
Let $L_i$ in $\fn$ be a sequence of leaves converging
to $F$ and on that side of $F$ so that $\pin L_i$ intersects
both $\alpha_i$ and $\beta_i$.
In 
the annulus ${\cal A} $ of circles at infinity, consider the
rectangle $R_i$ bounded by $(a_i, b_i)$ in $\pin F$,
the parts of $\alpha_i, \beta_i$ between $\pin F$ and
$\pin L_i$ and the small segment in $\pin L_i$ bounded
by \ $\pin L_i \cap \alpha_i$ \ and \ $\pin L_i \cap \beta_i$.
Then the sets $R_i$ converge to $z$  as $i$ converges
to infinity. This is proved in \cite{Fe7}.

From now on the $\epsilon$ is chosen small enough
to also satisfy
the conclusions of the definition above and also that
any set in $\mi$ of diameter less than $10 \epsilon$
is in a product box of $\fn$ and $\wwp$.
Given a curve $\zeta$ in a leaf $F$ with starting
point $p$ and limiting on a unique point $q$ in $\pin F$,
let $\zeta^*$ denote the geodesic ray of $F$ with
same starting and ideal points.

\begin{lemma}{}{}
Let $E$ be a leaf of $\wls$ and $p$ the starting point of the ray $r$ of 
$E \cap F$.
Assume that $r$ does not have any singularity.
For any $L$ near $F$, then $E \cap L$ has a ray
$r_L$ which is near $r$. The ideal points
of $r_L$ in $\pin L$ vary continuously with $L$
in the topology of ${\cal A}$ defined above.
\label{cont}
\end{lemma}

\begin{proof}{}
We do the proof for say the positive side of $F$.
We consider $r$ without singularity or else we would
have to check the 2 exterior rays in $\wlsf$ emanating from
$p$. We can always get a subray of $r$ which has
no singularities.

Let $u = r_{\infty}$.
 Choose 
contracting (or $\epsilon$ non expanding) directions  
in both sides of $u$, with ideal points very close to $u$.
Let them be defined by geodesic rays $r_0, r_1$ starting at $p$.
There is $\tau$ a small flow segment starting at $p$ and
 in that side of $F$  so that for any $L$ intersecting
$\tau$, then $L$ is asymptotic to $F$ along the
$r_0, r_1$ rays, or at least always $\leq \epsilon$ from $F$.
Hence $r_0, r_1$ flow along $\wwp$ to $L$. Let 
$s_0, s_1$ be the flow images in $L$. 
The $\epsilon$ is also chosen small enough so that 
$s_0, s_1$ have geodesic curvature very small
(this $\epsilon$ depends only on $M$ and $\fol$).
In particular the curves $s_0, s_1$ 
are a small bounded distance (depending only on $\epsilon$) 
from
the corresponding geodesic arcs $s^*_0, s^*_1$.
Let the 
ideal points of $s_0, s_1$  in $\pin L$
be denoted by  $v_0, v_1$
and let $J_L$ be the small closed interval
in $\pin L$ bounded by $v_0, v_1$.
Then $v_0, v_1$ 
are in the markers
associated to $r_0, r_1$ respectively and so they vary continuously
with $L$.

\begin{figure}
\centeredepsfbox{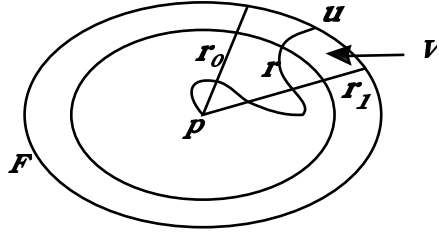}
\caption{
Leaf in wedge defined by markers.}
\label{bet}
\end{figure}

Consider $\xi = E \cap L$ and the rays $l$ of
$\xi$ starting at $\tau \cap L$ and
containing some points which flow back to points
in $r$. It may be that $\xi$ has singularities $-$
even if $r$ does not $-$ but there are only
finitely many such rays.
We want to prove that the ideal point of any such 
is in $J_L$.
As the rectangles $R_i$ defined above converge to 
$u$ in 
${\cal A}$ this will prove the continuity property of the lemma.

Choose $d > 0$ so that outside of a disk $D$
of radius $d$ in $F$, then
$r$ is in the small wedge $W$ of $F$ defined by $r_0, r_1$,
see fig. \ref{bet}.
Choose $\tau$ small enough so that if $L$ intersects $\tau$, then
the entire disk $D$ is
$\epsilon$ near $L$.
Let $V$ be the closure in $F$ of $W - D$.
The boundary $\partial V$ consists of subrays 
of $r_0, r_1$ and an arc in $\partial D$. Therefore
all points in $\partial V$ are less than $\epsilon$
from $L$ and flow to $L$ under $\wwp$ with image
a curve $\gamma$. This curve contains subrays of 
$s_0, s_1$ and it is properly embedded in $L$.
Points of $F$ near $\partial V$ also flow to $L$ so there
is a unique component $U$ of $L - \gamma$ which has
some points flowing back to points in $V$.
We want to show that the ray $l$ is eventually
contained in $U$.

Let $r_{init}$ be the subarc of $r$ between $p$ and
the last point $c_0$ of $r$ in $D$.
As $p$ and $c_0$ flow into $L$, then proposition \ref{conn}
shows that the entire arc $r_{init}$ flows into $L$
and let $\delta$ be its image in $L$.
As $r$ is singularity free, then so is $\delta$ and
hence $\delta$ is contained in any ray $l$ of $E \cap L$
in that direction.
After $c_0$ the curve $r$ enters $V$ and so $l$ must
enter $U$ after $\delta$.
If after that the ray $l$ exits $U$ then it must
cross $\partial U = \gamma$ in some point, call it $c_1$.
But $c_1$ flows back to $F$ and one can apply
proposition \ref{conn} again in the backwards direction
to show that $c_1$ has to flow to a point in $r$.
This contradicts the choice of $c_0$.

This shows that $l$ is eventually entirely contained
in $U$ and therefore $l_{\infty}$ is a point
in $J_L$. This shows the continuity property 
as desired and finishes the proof of the lemma.
\end{proof}

Now we have a property which will be crucial to a lot of our analysis.

\begin{proposition}{}{}
Suppose that $\fol$ is not topologically conjugate to the
stable foliation of a suspension Anosov flow.
Then the set of ideal points of rays of $\wlsf$ is dense in
$\pin F$.
\end{proposition}

\begin{proof}{}
Suppose that there is $F$ in $\fn$ so that the set of
ideal points in $\wlsf$ is not dense in $\pin F$.
Let $J$ be an open interval in $\pin F$ free of such
ideal points. Choose $p_i$ in $F$, $p_i$ converging
to a point in $J$. 
The visual angle of $J$ as seen from $p_i$ converges
to $2 \pi$, so the complementary wedge $W_i$
with corner  $p_i$ has angle which converges to
zero. Up to subsequence assume that 
$g_i(p_i)$ converges to $p_0$ in a leaf $L$ of $\fn$
and the small wedges $g_i(W_i)$ converge to a geodesic
ray $s$ in $L$ with ideal point $z$.

\vskip .15in
\noindent
{\underline {Claim }} $-$ In $L$ all the rays of $\wlsl$ converge
to $z$.

Suppose there is
$x$ different from $z$ which is an ideal point of a ray
$r$ in $\wlsl$. Then $r$ is contained in $\ws(c_0)$ for
some $c_0$ in $\mi$  and
for $g_i(F)$ sufficiently near $L$ then $\ws(c_0)$ intersects
$g_i(F)$. Any 
ray of $\ws(c_0) \cap g_i(F)$ which is near $r$ will have ideal point
near $x$ in the topology of corresponding annulus ${\cal A}$
of ideal circles near $\pin L$. This is a consequence of the previous
lemma.
But $g_i(W_i)$ converges to $r$ in this topology of ${\cal A}$,
so the sets $g_i(\pin F - J)$ converge to $z$ in ${\cal A}$.
There are no ideal points of leaves of 
$\widetilde {\Lambda}^s_{g_i(F)}$ in $g_i(J)$.
This contradicts the fact that the ideal points above are
very near $x$ and proves the claim.
\vskip .1in

The proof of the proposition
is similar to that of theorem \ref{as1}.
As in that theorem consider the set of possible limits
$g_i(p_i)$ as above. This projects to a lamination in $M$
and let ${\cal L}$ be a minimal sublamination.
The claim shows that each leaf of $\wcl$ has a distinguished
ideal point towards which all rays of $\wlsl$ converge.
The arguments in the claim also prove that if $\tau$ is
a transversal to $\fn$, then the ideal points of
leaves of $\wcl$ intersecting $\tau$ vary continuously
in the corresponding ideal annulus. Because of the
distinguished ideal point property, then each leaf of
${\cal L}$ has fundamental group at most ${\bf Z}$.
If needed lift to a double cover so that all leaves
of $\fol$ are orientable. Hence a leaf of ${\cal L}$ is
either a plane or an annulus.

Consider a complementary component $U$ of ${\cal L}$
and a boundary leaf $A$ of $U$. If $A$ is a plane then
as in the proof of theorem \ref{as1},
the region $U$ is an $I$-bundle over $A$ and the flow
$\Phi$ is a product in $U$. This region can be
collapsed away.

Suppose now that $A$ is an annulus.
Assume that flow lines through $A$ flow into $U$.
Again we want to show that $U$ is a product region.
As in the proof of theorem \ref{as1} let $A_1, A_2$ be
two noncompact, disjoint annuli in $A$ with $A - (A_1 \cup A_2)$
a compact annulus and $A_1, A_2$ contained in the 
thin, $I$-bundle region. Then $A_1, A_2$ flow entirely
into leaves $B$ and $C$ in $\partial U$. 
Suppose first that $B, C$ are different. Lift to the
universal cover to produce lifts $\widetilde U, \widetilde
A, \widetilde A_1, \widetilde A_2, \widetilde B, \widetilde C$.
Then $\widetilde A_1, \widetilde A_2$ are disjoint half
planes of $\widetilde A$ which flow positively respectively
into $\widetilde B$ and $\widetilde C$. 
Let $g$ be the generator of the isotropy group of $\widetilde A$,
which has fixed points in $z, x$ where $z$ is the
distinguished ideal point in $\widetilde A$.
The argument will show there is a leaf in
$\widetilde {\Lambda}^s_{\widetilde A}$ which also has
ideal point in $x$, contradiction.

From a point in 
$\widetilde A_1$ draw a geodesic segment of 
$\widetilde A$ to a point
in $\widetilde A_2$. Let $p$ be the first point
of this segment which does not flow positively into 
$\widetilde B$.
Then $\Theta(p)$ is in the boundary of $\Theta(\widetilde B)$.
Also  points in the segment near $p$ flow to $\widetilde B$
in positive time, hence there is a slice leaf $l$ of
$\oos(\Theta(p))$  which is in the boundary of $\Theta(\widetilde B)$.
Notice that every point in $l$ is a limit of points in $\Theta(B)$ on
that side.
The set $(l \times \rrrr)$ intersects $\widetilde A$ in at least
$p$: if
$l$ is contained in $\Theta(\widetilde A)$ then it 
generates a properly embedded copy of the reals
in a leaf $s$ of $\widetilde {\Lambda}^s_{\widetilde A}$
otherwise the part that is contained  in $\Theta(\widetilde A)$
also does. 
Every point of $s$ is a limit of points  that flow
positively into $\widetilde B$. Therefore no point in $s$
can flow positively in $\widetilde C$ or else we would have
points flowing both in $\widetilde B$ and $\widetilde C$.

This shows that the leaf $s$ of $\widetilde {\Lambda}^s_{\widetilde A}$
is a bounded distance from the axis $r$ of $g$.
Iterate $s$ by powers of $g$ acting with $z$ as an expanding
fixed point.
The iterates $g^n(s)$ with $n > 0$ are all distinct. 
Either they are all nested or they are disjoint.
If they are not nested since they all have to be in a bounded
distance neighborhood of the axis of $g$ and have both
endpoints in $z$, then eventually they will have two 
points which are far along the leaf, but close in $\widetilde A$.
By Euler characteristic reasons, this would force
a center or one prong singularity, which is impossible.
Hence they are nested, increasing and 
they limit to a leaf of $\widetilde {\Lambda}^s_{\widetilde A}$
which has ideal limit points in $z$ and $x$. This is a contradiction.
This shows that $B = C$.
In fact the same arguments show that all of the points
in $A$ flow into $B$, since that happens for the complement
of a compact annulus in $A$ and then the arguments above 
apply here.
Hence $U$ is a product region.
Therefore
we can collapse $\fol$
to a minimal foliation. 

%
%

\vskip .1in
As in theorem \ref{as1}
we can then show that $\fol$ is $\rrrr$-covered. Suppose this
is not the case and let $F_1, F_2$ be non separated leaves.
Let $L_i$ in $\fn$ leaves converging to both $F_1, F_2$.
Let $u_1, u_2$ be the distinguished ideal points
in $\pin F_1, \ \pin F_2$  respectively.
Let $a_1, b_1$ be points in $\pin F_1$ very near 
$u_1$ and on opposite sides of $u_1$ and which are
in markers associated to contracting or $\epsilon$
non expanding directions in $F_1$ associated to the $L_i$
side. 
Let $r_1$ be the geodesic in $F_1$ with ideal
points $a_1, b_1$.
Similarly for $F_2$ producing $a_2, b_2, r_2$.
For $i$ big enough $L_i$ is at most $\epsilon$ far
from all points in $r_1, r_2$. Therefore $r_1$ flows (by $\wwp$)
into a curve $s_1$ in $L_i$ and $r_2$ flows into
$s_2$.  This implies that $s_1, s_2$ are disjoint
in $L_i$.
Also $s_1$ has ideal points $a'_1, b'_1$
which are in markers containing $a_1, b_1$ respectively
(this is using a transversal to $\fn$ through a point
in $F_1$). Similarly $s_2$ has ideal points $a'_2, b'_2$
in markers containing $a_2, b_2$
(using transversal to $\fn$ through a point in $F_2$).
As $s_1, s_2$ are disjoint then $a'_1, b'_1$ do not
link $a'_2, b'_2$ in $\pin L_i$, see fig. \ref{link}.

\begin{figure}
\centeredepsfbox{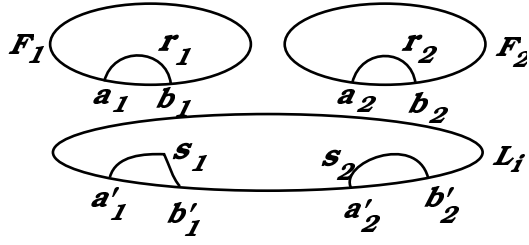}
\caption{
Pushing ideal points near.}
\label{link}
\end{figure}

The ideal point $a'_1$ cannot be in a marker to $\pin F_1$
and to $\pin F_2$ at the same time since they are non
separated leaves. Hence the points $a'_1, b'_1,
a'_2, b'_2$ are all distinct.
Let $J_1$ be the interval of $\pin L_i$ bounded by
$a'_1, b'_1$ and not containing the other points
and similarly define $J_2$.
For simplicity we are ommitting the dependence of $J_1, J_2$ 
on $L_i$ (or on $i$).
Now consider $E$ a leaf of $\wls$ intersecting $F_1$.
Then $E \cap F_1$ has a ray with ideal point $u_1$,
which is in the interval $(a_1,b_1)$ of $\pin F_1$.
The proof of lemma \ref{cont} shows that if
$L_i$ is close enough to $F_1$ then the ideal
points of the corresponding rays of $(E \cap L_i)$ 
have to be in $J_1$. In the same way using $F_2$ one
shows that the distinguished ideal point has to be
in $J_2$. Since $J_1, J_2$ are disjoint, this is
a contradiction. This shows that $\fol$ is
$\rrrr$-covered.

Since $\fol$ is $\rrrr$-covered then theorem \ref{as2} implies
that $\Phi$ can be chosen to be a pseudo-Anosov flow.

Also as $\fol$ is $\rrrr$-covered we can choose a transversal
$\tau$ intersecting all the leaves of $\fn$. This shows that
the union of all the circles at infinity has a natural topology
making it into a cylinder ${\cal A}$. This situation of
$\rrrr$-covered foliations is carefully analysed in \cite{Fe7}.
The fundamental group of $M$ acts in ${\cal A}$ by
homeomorphisms.
The union of the distinguished
ideal points of leaves of the distinct leaves of $\fn$
is a continuous curve $\zeta$ in ${\cal A}$ which
is group invariant.

Suppose first that $\fol$ admits a holonomy invariant
transverse measure. Since $\fol$ is minimal then
the transverse measure has full support.
Under these conditions Imanishi \cite{Im} proved
that $M$ fibers over the circle
with fiber a closed surface.
In addition $\fol$ is approximated arbitrarily near
by a a fibration.
The pseudo-Anosov flow is also transverse to these
nearby fibrations and so the same situation occurs
for the fibrations: there is a global invariant
curve in the cylinder at infinity. Since now
there are compact leaves, this is impossible.

We conclude that there is no holonomy invariant
transverse measure. Therefore Thurston's theorem
shows the existence of contracting directions
and not just $\epsilon$ non expanding directions.
So the markers are associated to contracting directions.
If $\zeta$ intersects a marker $m$, that corresponds
to a direction in a leaf of $\fn$ which is contracting.
Under the flow $\wwp$ this
gets reflected in the contracted leaves nearby, that
is the marker is contained in $\zeta$.
Since $\fol$ is minimal and $\zeta$ is $\pi_1(M)$ invariant,
this shows that the entire curve $\zeta$ is
a marker associated to  contracting directions.
The results from \cite{Fe7} apply here, in particular
lemma 3.17  through proposition 3.21 of \cite{Fe7}:
they show that no other direction in $\fn$ (outside
of $\zeta$) is a contracting direction.
By Thurston's theorem again, there would be a
holonomy invariant transverse measure, contradiction.

Therefore $\zeta$ has no contracting directions.
The same analysis of \cite{Fe7} now shows that 
for any leaf $F$ in $\fn$
and every direction other than the distinguished direction,
then it is a contracting direction.
In fact it is a contracting direction with any other leaf
of the foliation.

This is a very interesting situation. Let $a_F$ be the
distinguished ideal point of $F$ leaf of $\fn$.
Consider a one dimensional 
foliation in $\mi$ whose leaves are geodesics in leaves $F$ of $\fn$ 
and which have one ideal point $a_F$.  Let $\widetilde \xi$ be the
flow which is unit speed tangent to this foliation and moves
towards the ideal point $a_F$. 

This is a flow in $\mi$. Clearly in each leaf of $\fn$, it is a
smooth flow. If $q_i$ in $L_i$ of $\fn$
converge to $q$ in $L$, then the
geodesics of $L_i$ with ideal point $a_{L_i}$ converge to the 
geodesic through $q$ in $L$ 
with ideal point $a_L$. This is 
because the ideal points $a_F$ vary continuously with $F$
and $q_i$ converges to $q$ $-$ this is the local
trivialization of the union of the circles at infinity
using the tangent bundles to a transversal.
Hence $\widetilde \xi$ varies continuously.

Since $\zeta$ is group invariant, this induces a flow 
in $M$, which is tangent
to the foliation $\fol$. Clearly it is smooth along the leaves of
$\fol$ and usually just continuous in the transverse direction.

This flow is a topological Anosov flow: the stable
foliation is just the original foliation $\fol$. 
The unstable foliation: Let $p$ in leaf $L$ of $\fn$, let
$\gamma$ be the flow line of $\widetilde \xi$ through $p$.
Then $\gamma$
has positive ideal point  $a_L$ and negative ideal point
$v$.  As explained above $v$ is in a marker $m$ which is associated
to a contracting direction and so that $m$
intersects
all ideal circles. For each $F$ in $\fn$, let $m_F$ be
the intersection of $m$ and $\pin F$. Let $\gamma_F$ be
the geodesic in $F$ with ideal points $a_F$ and $m_F$.
Let $E_p$ be the union of these $\gamma_F$. Then all
orbits of $\widetilde \xi$ in $E_p$ are backwards
asymptotic by construction. By construction the $E_p$ are
either disjoint or equal as $p$ varies in $\mi$ and they
form a group invariant foliation in $\mi$. This is the
unstable foliation.
Hence $\xi$ is a topologically Anosov flow. 
Notice that in the universal cover 
every stable leaf intersects every unstable leaf
and vice versa.

By proposition \ref{susp} it follows that $\xi$ is topologically
conjugate to a suspension Anosov flow.
The foliation $\fol$ is then topologically conjugate
to the stable foliation of this flow.
This finishes the proof of this proposition.
\end{proof}

%

\vskip .1in
\noindent
{\bf {Remark}} $-$ The hypothesis is necessary. Suppose that
$\fol$ is the stable foliation of a suspension Anosov flow,
$\xi$ so that it is transversely orientable. Perturb the flow
slightly so that flow lines are still tangent to the original
{\underline {unstable}} foliation of  $\xi$. The new flow, 
call it
$\Phi$ is transverse to $\fol$, it 
has the same unstable foliation as $\xi$ but different
stable foliation. 
The flow $\Phi$ is not regulating for $\fol$. The intersections
of leaves of $\wls$ with leaves $F$ of $\fn$ are all horocycles
with the same ideal point which is the positive ideal point
of flow lines in $F$. So the ideal points of rays of leaves of
$\wlsf$ are not dense in $\pin F$. Notice  these leaves
are not quasigeodesics in $F$ either.
This example is studied in detail in section 7 of \cite{Fe9}.
\vskip .1in

Now we want to study metric properties of  slices of 
leaves of $\wlsf$.
The best metric property such leaves could have is that they
are {\em quasigeodesic}: this means that length along the
curve is at most a bounded multiplicative distortion
of length in the leaf $F$ of $\fn$ \cite{Th1,Gr,Gh-Ha,CDP}. 
If the bound is $k$ then we say the curve is a $k$-quasigeodesic.
Since $F$ is hyperbolic this would imply that such leaves
(the non singular ones) are a bounded distance from
true geodesics. 
Very unfortunate for us, this is not true in general. 
But there are still some good properties. 

Let $\hhs$ be the leaf space of $\wls$ and $\hhu$ be
the leaf space of $\wlu$.
Clearly since $\hhs$ may be non Hausdorff, it could
be that some $\wlsf$ does not have Hausdorff leaf space. 
This easily would imply that the slices of $\wlsf$ are
not uniformly quasigeodesic \cite{Fe1}.
This in fact occurs, see Mosher \cite{Mo1,Mo5}. Still it could be
that given a ray in $\wlsf$, it is a quasigeodesic $-$
with the quasigeodesic constant depending on the particular
ray. We are not able to prove this and we cannot conjecture
what happens in generality. But we are able to
prove a weaker property, which will be enough for our purposes.
If $r$ is a ray in a leaf of $\wlsf$, recall that 
$r^*$ is the unique
geodesic ray in $F$ with same starting point as $r$ and same
ideal point. We would like to prove that $r, r^*$ are a bounded
distance apart, but we do not know if that is true. But 
we can prove the following important property:

\begin{lemma}{}{}
There is $\delta_0 > 0$ so that for any $F$ in $\fn$ and any ray $r$ in
a leaf of $\wlsf$, then given any segment of length $\delta_0$ in
$r^*$, there is a point in this segment which is less than
$\delta_0$ from $r$ in $F$. That implies that $r^*$ is in the neighborhood
of radius $2 \delta_0$ of $r$ in $F$.
\label{boun}
\end{lemma}

\begin{proof}{}
This means that 
$r^* \subset N_{2 \delta_0}(r)$ in $F$. We do not know if the converse
holds. Suppose the lemma is not true. Then there are $F_i$ leaves
of $\fn$, $r_i$ rays of $\wlsfi$
and $p_i$ in $r^*_i$ so that $B_{2i}(p_i)$ (in $F_i$) does not
intersect $r_i$. There is one 
 side of $r^*_i$ in $F_i$ so that $r_i$ goes
around that side, see fig. \ref{aroun}, a.
Let $q_i$ inside a half
disk of $B_{2i}(p_i)$ 
with $B_i(q_i)$ tangent to $r^*_i$ and
$\partial B_{2i}(p_i)$, see fig. \ref{aroun}, a.

\begin{figure}
\centeredepsfbox{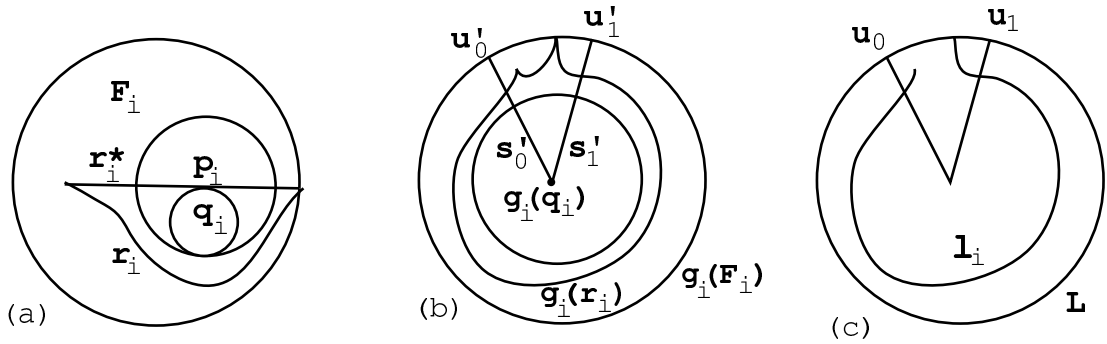}
\caption{
a. Limits of points, b. Going around disks
in $F_i$,
c The picture in $L$.}
\label{aroun}
\end{figure}

As usual up to subsequence there are $g_i$ in $\pi_1(M)$ with
$g_i(q_i)$ converging to $q_0$ in $L$ leaf of $\fn$
and so that the geodesic
segments $\zeta_i$ from $g_i(q_i)$ to $g_i(p_i)$ in $F_i$
converge to a geodesic
ray $s$ in $L$. Choose two markers with 
points $u_0, u_1$ in $\pin L$ very close to $s_{\infty}$
and on opposite sides of it.
The markers are associated to the side of $L$
where the $g_i(F_i)$ are limiting to.
Let $s_0, s_1$ be the geodesic rays of $L$ starting at $q_0$ and
with ideal points $u_0, u_1$.
For $i$ big enough $g_i(F_i)$ is $\epsilon$ close to both
$s_0$ and $s_1$ and so these two rays flow (under $\wwp$)
to curves $s'_0, s'_1$ in $g_i(F_i)$. 
The ideal points $u'_0, u'_1$ of $s'_0, s'_1$ are in the markers
above.

For $i$ big enough the ray $g_i(r_i)$ has a subray which goes 
around $g_i(B_i(q_i))$ in $g_i(F_i)$ and has ideal
point in the small segment of $\pin g_i(F_i)$ defined
by $u'_0, u'_1$, see fig. \ref{aroun}, b.
Since $s'_0, s'_1$ flows back to $L$ this figure flows
back to $L$ producing a ray $l_i$ of $\wlsl$
which goes around a big disk in $L$ centered at $q_0$
and has ideal point in the small segment bounded by
$u_0, u_1$, see fig. \ref{aroun}, c.
As $i$ goes to infinity, these $l_i$ escape
to infinity in $L$ because bigger and bigger disks
in $g_i(F_i)$ flow to $L$. 
This implies that there is no ideal point of a ray
of $\wlsl$ outside the small segment of $\pin L$
bounded by $u_0, u_1$. This contradicts
the previous proposition that such ideal points are
dense in $\pin L$.

This finishes the proof of the lemma.
\end{proof}

\begin{lemma}{}{}
The limit points 
of rays of $\wlsf$ 
vary continuously in $\pin F$
except for the
non Hausdorffness in the leaf space of $\wlsf$.
\label{encon}
\end{lemma}

\begin{proof}{}
Suppose that $p_i$ converges to $p$ in $F$, with respective rays
$r_i$ converging to the ray $r$ of $\wlsf$. 
Let $l$ be the leaf of $\wlsf$ through $p$.
Up to subsequence
assume the $r_i$ are all in the same sector of $l$ defined by $p$
and that they form a nested sequence of rays.
Then the ideal ponts 
$(r_i)_{\infty}$ form a monotone sequence
in $\pin F$. Perhaps some
ideal points are the same.
If $(r_i)_{\infty}$ does not converge to $r_{\infty}$ there
is an interval $v$ in $\pin F$, between the limit and $r_{\infty}$.
Since the ideal points are dense in $\pin F$, there is $w$ leaf
of $\wlsf$ with $w_{\infty}$ in $v$. Therefore there is $l'$ not
separated from $l$ with $r_i$ converging to $l'$ as well.
In this fashion we can go from $l$ to $l'$. 
This shows that if there is no leaf of $\wlsf$ non separated from
$l$ in that side and in the direction the rays $r_i$ go, then
the limit points vary continuously.

We analyse a bit further the non Hausdorffness.
In the setup above there are subrays of $r_i$ with points
converging to a point in $l'$ and we can restart the
analysis with $l'$ instead of $l$.
If there are finitely many leaves non separated from $l$ and
$l'$ we can assume that $l, l'$ are consecutive. Then
they have subrays which share an ideal point. 
If $m$ is the last leaf non separated from $l, l'$ 
in the direction the rays $r_i$ go to, then
there is a ray $\zeta$ of $m$ so that there are subrays
of $r_i$ with points converging to a point in $\zeta$
and $(r_i)_{\infty}$ converges to $\zeta_{\infty}$.
If there are infinitely many such leaves non separated
from $l$, then we can order them as $\{ l_j \}, j \in {\bf N}$ 
all in the
direction the rays $r_i$ go to. The ideal points of
\ $l_j$ \ form a monotone sequence in $\pin F$ which converge
to a point $u$ in $\pin F$. The arguments above show
that $(r_i)_{\infty}$ converges to $u$.
\end{proof}

Our next goal is to analyse the 
non Hausdorffness in the leaf space of $\wlsf$.
We also want to understand when can the ideal points of two
different rays of 
the same leaf of $\wlsf$ be the same.
A {\em Reeb annulus} is an annulus $A$ with a foliation so that
the boundary components are leaves and every leaf in the
interior is a topological line which spirals towards 
the two boundary components in the same direction.
In the universal cover the lifted foliation does
not have Hausdorff leaf space.
The lifted foliation to the universal cover is 
called a {\em Reeb band}.

\begin{define}{(spike region)}{}
A stable spike region in a leaf $F$ of $\fn$
is a closed $\wlsf$ saturated set ${\cal E}$
satisfying:

\begin{itemize}

\item
There are
finitely many  boundary leaves of ${\cal E}$
which are line leaves of $\wlsf$.
The ideal points of consecutive rays in the boundary of ${\cal E}$
are the same, otherwise they are distinct (like an ideal
polygon).

\item
The region ${\cal E}$ is a bounded distance
from the ideal polygon with these vertices. The bound
is not universal in $\fn$. 

\item
There is an ideal
point  $z$ of  ${\cal E}$
so that every leaf in the interior of ${\cal E}$ has both ideal
points equal to $z$. In addition the leaves in the
interior are nested.
The finitely many leaves in the boundary
are all non separated from each other and they are limits
of the interior leaves.

\item
There is no singularity of $\wlsf$ in the interior of
${\cal E}$. If there is a singularity of $\wlsf$ 
in a boundary leaf $\tau$ of ${\cal E}$ then 
the interior of ${\cal E}$ is contained in the sector
defined by the line leaf $\tau$.
\end{itemize}

\noindent
Similarly define an unstable spike region. A spike region is either
a stable or unstable one.
\end{define}

\begin{proposition}{}{}
Let $E$ be a leaf in $\fn$ and $\upsilon$ a slice of a leaf
$\upsilon_0$ of $\wlse$.
Suppose that both ideal points of $\upsilon$ are the same.
Then $\upsilon$ is contained in 
the interior of a stable spike region $B$ of $E$.
In addition either $B$ projects to a Reeb annulus in a leaf
of $\fol$ or for any two consecutive rays in $\partial B$,
the region between them projects to a set asymptotic to
a Reeb annulus in a leaf of $\fol$.
Similarly for $\wluf$.
\label{band}
\end{proposition}

\begin{proof}{}
We do the proof for $\wlsf$.
Let $\upsilon$ be a slice as above with ideal point $x$ in 
$\pin E$. 
Let $C$ be the region bounded by $\upsilon$ in $E$ which
only limits in $x$.
First we assume 
assume that $\upsilon$ is a line leaf of some leaf of
$\wlsf$.
We will show that 
the region $C$ as it approaches $x$, projects to a set
in $M$ which limits to a Reeb annulus in  a leaf of $\fol$.
The process will be done in a series of steps.
The proof of this proposition is very long with
several intermediate results and lemmas.

Choose $z_0$ in $\upsilon$ and let $e_1, e_2$ be the rays
of $\upsilon$ defined by $z_0$.
Let $\zeta^*$ be the geodesic ray of $E$ starting at
$z_0$ and with ideal point $x$. Then $\zeta^*$ is 
contained in the $2 \delta_0$ neighborhood of $e_1$ or $e_2$,
where $\delta_0$ is the constant of lemma
\ref{boun}. It follows that 
we can choose $p_i, q_i$ in $e_1, e_2$ respectively 
with $p_i, q_i$ converging to $x$ in $E \cup \pin E$
and also $d_E(p_i, q_i) < 4 \delta_0$. 
Let $e^i_1$ be the subray of $e_1$ starting at $p_i$
and $e^i_2$ the subray of $e_2$ starting at $q_i$.
Up to subsequence  there are $p_0, q_0$ in $\mi$ and
are $g_i$ in $\pi_1(M)$ with $g_i(p_i)$, $g_i(q_i)$
converging to $p_0$, $q_0$ respectively. 
The distance condition implies 
$p_0, q_0$  are in the same leaf of $\fn$, let $F$ be this leaf.
Then $g_i(E)$ converges to $F$
and perhaps other leaves as well.

For $i$ big enough
the flowlines of $\wwp$ through $g_i(p_i), g_i(q_i)$
go through to $u_i$ and $v_i$ in $F$. Also $u_i \rightarrow p_0,
v_i \rightarrow q_0$.
If the leaf of $\wlsf$ through $p_0$ contains $q_0$ then 
for $i$ big enough the arcs 
in leaves of $\wlsf$ from $u_i$ to $v_i$
will have bounded length and
bounded diameter. The same will happen for the arcs of
of $g_i(\upsilon)$ between
$g_i(p_i)$ and $g_i(q_i)$, contradiction.
Hence $p_0, q_0$ are not in the same leaf of $\wlsf$. Let
$l$ be the leaf of $\wlsf$ through $p_0$ and $r$ be the one
through $q_0$. Let $L, R$ be leaves of $\wls$ containing $l$ and
$r$ respectively.
Since the intersection of
a leaf of $\wls$ with $F$ is connected, then $L$ and $R$ are
distinct and also are not separated from each other
in the leaf space of $\wls$.

For simplicity assume that the leaves of $\wls$ through
$u_i$ form a nested collection with  $i$.

\vskip .1in
The first goal is to show that we can choose $l, r$ line leaves
of $\wlsf$ as above so that they also share an ideal
point.
Let $\beta_i$ be a ray in the leaf of $\wlsf$ through
$u_i$ starting at $u_i$ and containing points in the 
flowlines which go to the ray
$g_i(e^i_1)$. Similarly let $\gamma_i$ be a subray
in the same leaf starting at $v_i$ and associated with 
the ray $g_i(e^i_2)$.
Let ${\cal C}_1$ \ (resp. ${\cal C}_2$) be the collection of line leaves of 
$\wlsf$ that $\beta_i$ \ (resp. $\gamma_i$) converges to, including the
ray of $l$ \ (resp. $r$).
Let ${\cal C}$ be the collection of all line leaves of $\wlsf$ which
are non separated from $l, r$. 
Then ${\cal C}$ contains ${\cal C}_1$ and ${\cal C}_2$. 
For any element $\tau$ in ${\cal C}$, let
$B(\tau)$ be the leaf of $\wls$ containing it.
 All of the $B(\tau)$ are not separated
from each other, and they are in the set of leaves ${\cal B}$ of $\wls$
non separated from both $L, R$.
By theorem \ref{theb}, the set ${\cal B}$ has a linear order,
making it order isomorphic to either ${\bf Z}$ or a finite
set. This induces an order in ${\cal C}$ where we can choose
this so that an arbitrary element of ${\cal C}_1$ is
bigger than any element in ${\cal C}_2$.

\begin{figure}
\centeredepsfbox{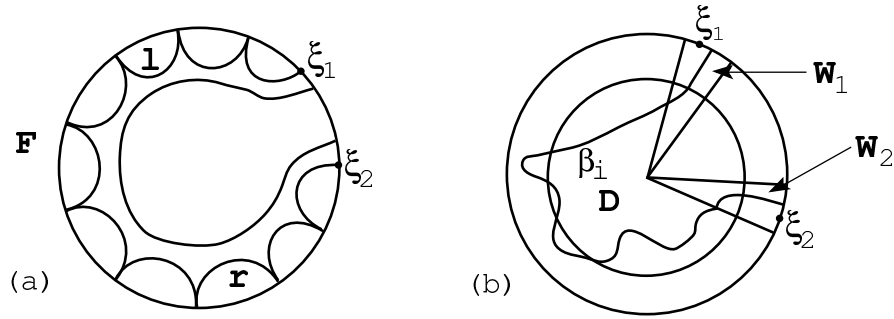}
\caption{
a. Non Hausdorffness in the limit,
b. Showing $\xi_1 = \xi_2$.}
\label{spre}
\end{figure}

If there are finitely many elements in ${\cal C}_1$ let 
$l'$ be the last one and let $\xi_1$ be the ideal point of
the ray of $l'$ corresponding to the direction of the
rays $\beta_i$.
Otherwise the ideal points of the leaves
in ${\cal C}_1$ form a weakly monotone sequence in $\pin F$
and let $\xi_1$ be the limit of this sequence.
Similarly define $\xi_2$ associated to $r$, see fig. \ref{spre},
a.

Fix a basepoint $x_0$ in $F$.
The first thing to prove is the following:

\begin{lemma}{}{}
$\xi_1 = \xi_2$.
\end{lemma}

\begin{proof}{}
Suppose by way of contradiction that this is not true.
Choose 2 markers very near $\xi_1$ bounding
an interval $J_1$ in $\pin F$ with $\xi_1$ in the
interior and similarly choose markers near $\xi_2$ and
interval $J_2$ so that $J_1, J_2$ are disjoint.
Let $W_1$ be the wedge of $F$ centered at the point
$x_0$ with
ideal set $J_1$ and $W_2$ the wedge of $F$ centered
also at $x_0$ with ideal set $J_2$. For i big enough
both boundaries of $W_1$ and $W_2$ flow into
$g_i(E)$. 

Suppose first that there is a last leaf $l'$ in ${\cal C}_1$.
Then $l'$ has a ray which is eventually contained in
a strictly smaller wedge $W'_1$, since 
the ideal point of $l'$ is  $\xi_1$.
Now choose a big disk $D$ of $F$ centered in $x_0$ .
Let $N_1$ be the closure of $W_1 - D$.
Choose $D$ big
enough so that $l'$ enters $N_1$ through $\partial D$
and is then entirely in $W'_1$. 
For $i$ big enough $\beta_i$ will be close to $l'$ for 
a long distance. By lemma \ref{encon} the ideal 
points of $\beta_i$ converge to $\xi_1$ as $i$ converges
to infinity, since $l'$ is the last leaf non separated
from $l$ in that side. The ideal point is in the limit
set of the subwedge $W'_1$. If the rays $\beta_i$ keep
exiting $W_1$ then since they are trapped by $l'$ and
$\beta_{i_0}$ (for some $i_0$), it follows that  the
sequence \ $\beta_i$ \ 
has additional limits besides the leaves in ${\cal C}_1$,
contradiction.
Therefore for big enough $i$, the $\beta_i$ enters $N_1$ through
$\partial D$ and stays in $N_1$ from then on.

We want to get the same result when ${\cal C}_1$ is infinite.
In that case let \ $\{ \nu_j, \ j \in {\bf N} \}$ 
\ be the leaves in
${\cal C}_1$ ordered with same ordering as in ${\cal C}_1$
and $\nu_1 = l$.
Since the leaves $v_i$ are non separated from each other
then they cannot accumulate anywhere in $F$ as $i \rightarrow \infty$
and the leaves $\nu_j$ escape compact sets as $j$ grows.
The ideal points of $\nu_j$
are also converging to $\xi_1$. By 
density of ideal points of $\wlsf$ in $\pin F$ the
leaves 
$\nu_j$
cannot be getting closer to non trivial intervals in $\pin F$.
This implies that there is $j_0$ so that for 

$$j \geq j_0, \ \ 
\nu_j \ \ {\rm  is \ very \ close  \ to } \ \ 
\xi_1 \ \ {\rm in} \   F \cup \pin F$$

\noindent
and so contained in $W_1$.
Now an argument entirely similar as in the case ${\cal C}_1$
finite implies that for $i$ big enough then $\beta_i$ has
subrays entirely contained in $N_1$.
The same holds for $\gamma_i$ producing subrays entirely
contained in the corresponding set $N_2$ $-$ the disk 
$D$ may need to be bigger to satisfy all these conditions.

There is $a_1 > 0$ and 
$i_0$ so that for $i \geq  i_0$ then except for the initial
segment of length $a_1$ then  $\beta_i$ is entirely
contained in $N_1$ and similarly for $\gamma_i$ and $N_2$.
Choose $k_0$ big enough so that $D$ is $\epsilon$ close
to $g_k(E)$ for any $k \geq k_0$.
Then $D$ flows in $g_k(E)$ under $\wwp$
and so do $\partial W_1, \partial W_2$.
For $i$ bigger than both $i_0, k_0$
the ray $\beta_i$ flows into the ray $g_i(e^i_1)$
(notice these do not
have singularities). 
The ray $g^i(e^i_1)$ has to be in the generalized wedge
which is bounded by the image of $\partial W_1$ in $g_i(E)$.
Similarly for $\gamma_i$. 
This argument is done in lemma \ref{cont}.
These two generalized wedges
have disjoint ideal sets in $\pin g_i(E)$. Therefore
$g_i(e^i_1)$ and $g_i(e^i_2)$ do not have the same
ideal points. This is a contradiction
because $e_1, e_2$ have the same ideal point in $\pin E$. 

This proves that $\xi_1 = \xi_2$.
\end{proof}

\noindent
{\underline {Continuation of the proof of proposition \ref{band}}}

The first part of the proof was this: in $E$ zoom in
towards an ideal point $x$ of  
$\pin E$ and use covering translations $g_i$ of $\mi$
to map back these points near a point in $\mi$ which
is in a leaf $F$.
We will redo this process starting with $F$.
By taking translates of $F$ we will limit to a leaf $F^*$.
The difference is that now we have leaves
of $\wlsf$ which are non separated from each other.
These non separated leaves are much better suited
to perturbation arguments as seen below.

The lemma shows that $\xi_1 = \xi_2$ and this
implies that the ideal points
of $\beta_i, \gamma_i$ are all the same and equal to $\xi_1$.
Let $\xi = \xi_1$.
The $\beta_i, \gamma_i$ are rays in leaves of
$\wlsf$ and contained in $F$.
Let $\mu$ be the geodesic ray in $F$ starting at $x_0$ 
(the basepoint in $F$) with
ideal point $\xi$.
Since $(\beta_i)_{\infty} = (\gamma_i)_{\infty} = \xi$,
then lemma \ref{boun} implies that for $z$ in
$\mu$ far enough from $p_0$, we can choose
a point in $\beta_i$ which is less than 
$2\delta_0$ away from $z$ in $F$.
Call this point $b_i(z)$. Similarly define $c_i(z)$ in $\gamma_i$.
This is for any $i$ in ${\bf N}$. 

For each $z$ we may take a subsequence
of the $b_i(z)$ which converges in $F$ and 
the limit is denoted by $b(z)$. Similarly define $c(z)$.
By definition of ${\cal C}_1$ the point $b(z)$ has
to be in one of the leaves of ${\cal C}_1$ and
similarly for $c(z)$.
The $b(z), c(z)$ are not uniquely defined and most likely
do not vary continuously with $z$.

\begin{lemma}{}{}
There is at least one element  $\zeta$ of ${\cal C}_1$ which
has ideal point $\xi$.
Similarly for ${\cal C}_2$.
\end{lemma}

\begin{proof}{}
If there are finitely elements in ${\cal C}_1$ then the
last one satisfies this property.
Suppose then there are infinitely many elements
in ${\cal C}_1$.
As $z$ varies in $\mu$, then so does $b(z)$. If 
there are $z$ escaping in $\mu$ so that $b(z)$ is
in the same element $\zeta$ of ${\cal C}_1$ then
$\zeta$ has an appropriate ray with ideal
point $\xi$. In this case we are done.

Otherwise we can find $z_k$ in $\mu$ converging to $\xi$
so that $b(z_k)$ are in leaves $\nu_{m(k)}$ of ${\cal C}_1$
which are all distinct. We can choose $z_k$ so that
the $m(k)$ increases with $k$.
In the same way we have $c(z_k)$ in 
distinct elements of ${\cal C}_2$.
Let 

$$B_k \ = \ \ws(b(z_k)), \ \ C_k \ = \ \ws(c(z_k)),
\ \ \ {\rm both \ in} \  \ {\cal B}$$

\noindent
Recall that ${\cal B}$ is the set of leaves of $\wls$
non separated from both $L, R$.
As the length from $b(z_k)$ to $c(z_k)$ in $F$
is bounded by $4\delta_0$, then up to subsequence assume
$\pi(b(z_k)), \pi(c(z_k))$ converge in $M$.
For $n, k$  big enough there is $h_{nk}$ covering
translation of $\mi$ so that $h_{nk}(b(z_n))$ is very close
to $b(z_k)$ and $h_{nk}(c(z_n))$ is very close to 
$c(z_k)$. Suppose $n >> k$, let $h = h_{nk}$ for simplicity.
Then $B_k$ has a point $b(z_k)$  very close to 
$h(b(z_n)) \in  h(B_n)$ and similarly $c(z_k)$ in 
$C_k$ very close to $h(c(z_n)) \in h(C_n)$.
But $B_k$ is non separated from $C_k$ and
similarly for $h(B_n), h(C_n)$, so the only
way this can happen is that 

$$h(B_n) \ = \ B_k, \ \ \ \ h(C_n) \ = \ C_k$$

\noindent
This implies that $h$ sends the set of leaves non
separated from $B_i, C_i$ to itself, that is $h$ acts
on the set ${\cal C}$ and therefore acts on ${\cal B}$ as
well. 
Notice that $B_k < B_n$ in the order of ${\cal B}$ 
because $n > k$ 
and $C_k \geq C_n$ (the $C_k$ could be all the same, but
if they are not then they decrease in the order).
Since $h(B_n) = B_k$ then
$h$ acts as a decreasing translation in the ordered
set ${\cal B}$.
But since $h(C_n) = C_k$ then $h$ acts as a non
decreasing translation. 
These two facts are incompatible.

This implies that we have to have at least one element
in ${\cal C}_1$ with ideal point $\xi$.
The same happens for ${\cal C}_2$.
This finishes the proof of the lemma.
\end{proof}


Since the sequence $\{ \beta_i \}$ 
also converges to the leaf $\zeta$ of $\wlsf$ we can rename
the objects and assume that $l = \zeta$ and $p_0$ is
a point in $l$. This can be accomplished by
choosing different points $p_i$ in the
ray $e_1$. Similarly do the same thing in the
other direction. We state this conclusion:

\vskip .1in
\noindent
{\bf {Conclusion}} $-$  There are $p_i, q_i$ in $e_1, e_2$
respectively, escaping these rays, so that 
$d_E(p_i, q_i) < 4 \delta_0$
and there are covering translations $g_i$ so that:
$g_i(p_i)$ converges to $p_0$, 
\ $g_i(q_i)$ converges to $q_0$,
both in $F$ and in rays $l, r$ of $\wlsf$.
Also $l, r$ converge to the same ideal point
$\xi$ in $\pin F$ and $l, r$ are not separated
in the leaf space of $\wlsf$.
\vskip .1in

%

We will continue this perturbation approach. 
We want to show that the region in $F$ ``between"
$l$ and $r$  projects to a Reeb annulus of $\fol$ in $M$.
Let then $z_i$ in $l$ converging to $\xi$ and $w_i$ in
$r$ converging to $\xi$, so that $d_F(z_i,w_i)$ is always
less than $4 \delta_0$. Up to subsequence assume there
are $h_i$ covering translations
with 

$$h_i(z_i) \rightarrow z_0, \ \ \ h_i(w_i) \rightarrow w_0$$

\noindent
Notice that $h_i(L), h_i(R)$ are 
non separated from each other and $h_i(L) \rightarrow \ws(z_0), 
\ h_i(R) \rightarrow \ws(w_0)$.
The argument in the previous lemma then implies
that $h_i(L) = h_j(L), \ h_i(R) = h_j(R)$ for all $i, j$
at least equal to some $i_0$.
Discard the first $i_0$ terms
and postcompose $h_i$ with $(h_{i_0})^{-1}$ 
(that is $(h_{i_0}^{-1} \circ h_i)$)
and assume that this is the original sequence $h_i$.
This implies that
$h_i(L) = L,
h_i(R) = R$ for all $i$.
So the $h_i$ are all in the intersection of
the isotropy groups of $L$ and $R$. This group is
generated by a covering translation
$h$.
Therefore there are $n_i$ with $h_i = h^{n_i}$. 
Since 
$h_i(z_i) \rightarrow z_0$ and the $\{ z_i, \ | \ i 
\in {\bf N} \}$ do not
accumulate in $\mi$ then 
$|n_i| \rightarrow \infty$.
In addition since $L, R$ are not separated from each other,
then $h$ preserves each individual line leaf, slice and
possible lift annulus of $L$.

Let $F^*$ be the leaf of $\fn$ containing $z_0, w_0$. Then
$h_i(F)$ converges to $F^*$.
Since $z_i$ is in $L$ and $h_i(L) = L$ then $h_i(z_i)$ is
in $L$ and so $z_0$ is in $L$. It follows that $L$ intersects
$F^*$.

Up to subsequence and perhaps taking the inverse of $h$,
assume that $n_i$ converges to $+\infty$.
If $h(F) = F$, then since $h(L) = L$ this produces
a closed leaf of $\Lambda^s \cap \fol$
in $\pi(F)$. Similarly 
$h(R \cap L) = R \cap L$ so produces another closed leaf in
$\pi(F)$ and together bound an annulus in $\pi(F)$ with a sequence of
leaves of $\Lambda^s \cap \pi(F)$
converging to the boundary leaves.
By Euler characteristic reasons, there can be no singularities
inside the annulus, so we conclude that the annulus in $\pi(F)$
has a Reeb foliation.

Now assume that $h(F)$ is not equal to $F$.
Let $\hp$ be the leaf space of $\fn$. This is a one dimensional
manifold, which is simply connected, but usually not
Hausdorff \cite{Ba2}. The element $h$ acts on $\hp$. An analysis of
group actions on simply connected non Hausdorff spaces
was iniatilly done by Barbot in
\cite{Ba2} and subsequently in \cite{Ro-St,Fe8}.
One possibility
is that $h$ acts freely in $\hp$. Then $h$ has an axis $\tau$
in $\hp$ which is invariant under $h$.  In general this
axis is not properly embedded, see \cite{Fe8}.
Since all the $h^{n_i}(F)$ intersect a common transversal,
then the analysis in \cite{Ba2} shows
that $F$ has to be in the axis of $h$ and 
$h^n(F)$ converges to a collection of non separated
leaves. In this case we get that $F^*$ and $h(F^*)$
are non separated from each other.

The other situation is that $h$ has fixed points in $\hp$.
In general the set of fixed points of $\hp$ is not a closed
set, but the set of points $z$ in $\hp$ so that $z$ and
$h(z)$ are not separated in $\hp$ is a closed subset
$Z$ of $\hp$ \cite{Ba2,Ro-St}.
None of the images of $F$ under $h$ can be in $Z$,
so $F$ is in a component of $\hp - Z$.
Then $h$ permutes these components. In addition
$h$ preserves an orientation in $\hp$ $-$ since
$\fol$ is transversely orientable.
Since
$h^{n_i}(F)$ all intersect a common transversal then
they have all to be in the same component $U$ of $\hp - Z$.
Let $i_0$ be the smallest positive integer
so that $h^{i_0}(U) = U$.
It follows that all $n_i$ are multiples of $i_0$.
Since $h^{n_i}(F)$ converges to $F^*$ then
the leaf $F^*$ is in the boundary of the component
$U$ and $h^{i_0}(F^*) = F^*$. 

The only remaining case to be analysed
is that $h$ acts freely and $h^n(F)$ converges
to $F^*$ with $h(F^*)$ non separated from $F^*$.
In this particular case we prove this is not
possible, that is:

%
%

\vskip .1in
\noindent
{\bf {Claim}} $-$ $h(F^*) = F^*$.

Suppose this is not true. 
The leaves $h(F^*), F^*$ are 
not separated in $\hp$. 
This implies that $\Theta(F^*)$ and $\Theta(h(F^*))$ are
disjoint subsets of $\oo$, see fig. \ref{jump}.
Therefore there are boundary leaves separating them.
But $L$ intersects both $F^*$ and $h(F^*)$
as $L$ intersects $F^*$ and is invariant under $h$.
Therefore both 
$\Theta(F^*)$ and $\Theta(h(F^*))$ intersect
the same stable leaf $\Theta(L)$.

Suppose that there is a stable boundary component of
$\Theta(F^*)$ separating it from $\Theta(h(F^*))$.
Then it has to be  a slice of $\Theta(L)$ as this
set intersects both of them. It would not be a line
leaf of $\Theta(L)$.
But as remarked before, $h$ leaves invariant all the
slices, line leaves and lift annuli of $L$ and
this contradicts $\Theta(h(F^*))$
being disjoint from $\Theta(F^*)$.
This implies there is
an {\underline {unstable}} boundary component
of $\Theta(F^*)$ separating it from $\Theta(h(F^*))$,
see fig. \ref{jump}.

\begin{figure}
\centeredepsfbox{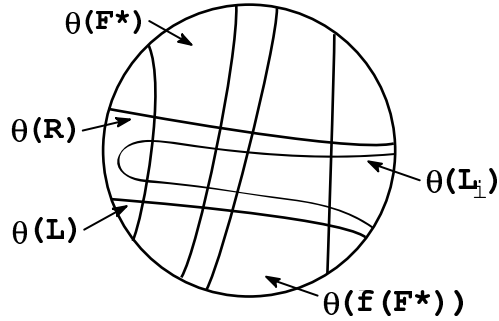}
\caption{
Contradiction in the orbit space $\oo$.}
\label{jump}
\end{figure}

%

In the same way $\Theta(R)$ intersects both $\Theta(F^*)$ and
$\Theta(h(F^*))$.
Let $L_i = \ws(u_i)$. Recall from the
beginning of the proof of  proposition \ref{band}
that $u_i, v_i$ are  points in $F$ with $u_i$ converging
to $p_0$ in $L$ and $v_i$ converging to $q_0$ in $R$.
Then 
$\Theta(L_i)$ converges to 
$\Theta(L) \cup \Theta(R)$ (maybe other leaves as well).
So $\Theta(L_i)$ intersects $\Theta(F^*)$ and $\Theta(h(F^*))$ 
for $i$ big enough.
The intersection of $\Theta(L_i)$ with 
at least one of $\Theta(F^*)$ or $\Theta(h(F^*))$ 
cannot be connected,
see fig. \ref{jump}.
This contradicts  propostion \ref{conn}.
This contradiction implies that
$h(F^*) = F^*$ and proves the claim.

\vskip .1in
So far we have  proved the following:
in any case there is $i_0$ a positive integer
so that if $f = h^{i_0}$ then
$f(F^*) = F^*$. As $f(L) = L$ then
$f(F^* \cap L) = F^* \cap L$
and similarly $f(F^* \cap R) = F^* \cap R$.
This produces an annulus $B$ in 
$\pi(F^*)$ with a Reeb foliation.
The region of $F^*$ bounded by $F^* \cap R$ and
$F^* \cap L$ bounds a band $B$
which is a bounded distance from a geodesic in
$F^*$ and projects to a Reeb annulus in a leaf of $\fol$.

But to prove proposition \ref{band},
we really want these facts for $F$ and not just $F^*$.
That is, we want a region in $\pi(E)$ which spirals
towards a Reeb annulus.
This turns out to be true: $\pi(E)$ has points converging
to $\pi(F)$ and $\pi(F)$ has points converging to a Reeb
annulus in $\pi(F^*)$. Since the annulus is compact,
it turns out the second step is unnecessary.
This depends on an analysis of holonomy of the foliation
$\fol$ near the annulus in $\pi(F^*)$ as explained below.

\vskip .1in
\noindent {\bf {Claim}} $-$ The point
$\pi(p_0)$ of $\pi(F)$ is in the boundary of
a Reeb annulus of $\fol$ contained in $\pi(F)$.
This implies that $F = F^*$.

%

The point $z_0$ is in $F^* \cap L$. Then $\pi(z_0)$ is
in $\pi(F^* \cap L) = \alpha$ which is a closed curve
since $h^{i_0}$ leaves invariant both $F^*$ and $L$
and their intersection is connected.
Previous arguments in the proof imply that for
$i$ big enough $h_i(z_i)$ is in the same local
sheet of $\wls$ as $z_0$.
Hence the points $\pi(z_i)$ are in $W^s(\pi(z_0)) = \pi(L)$
and converge to $\pi(z_0)$. This shows that $\pi(F \cap L)$
is asymptotic to $\alpha$ in the direction corresponding
to the projection of the direction of escaping $z_i$ in
the ray of $F \cap L$.
Namely $\alpha$ has contracting holonomy (of $\fol$) in
the side the $\pi(z_i)$ are converging to and eventually
$\pi(z_i)$ is in the domain of contraction of $\alpha$.

This means that the direction of $F$ associated to the ideal
point $\xi$ is a contracting direction towards
$F^*$. The rays in the leaves
$F^* \cap L, \ F^* \cap R$ in $F^*$ are
a bounded distance from a geodesic ray in $F^*$ with
same ideal point. The contraction above implies
that the corresponding rays $F \cap L, \ F \cap R$ of $F$ 
are also a bounded distance from a ray in $F$ with ideal
point $\xi$.

Now recall the points $p_i$ in $E$. We have $g_i(p_i)$
very close to $p_0$ in the leaf $l$ of $\wlsf$.
Also $\pi(l)$ is eventually in a region contracting
towards a Reeb annulus of $\fol$.
Hence if $i$ is big enough the $g_i(p_i)$ will also
be in this region. The leaf of $\Lambda^s \cap \fol$
through $\pi(p_i)$ will
be contracted towards the Reeb annulus in that direction.
This implies that the limit of the $\pi(p_i)$ is
already in a Reeb annulus, consequently the limit
of the $g_i(p_i)$ is
already in a Reeb band. 
%
%
%
%
%
%
%
%
%
%
%
%

It now follows that 
$\pi(F) = \pi(F^*)$.
That means that the second perturbation procedure (from points
in $F$ to points in $F^*$) in fact does not
produce any new leaf. 
This implies that up to covering translations then
the leaf $E$ is asymptotic to $F$ in the direction of
the ideal point $x$ in $\pin E$.
Let $V$ be the region of $E$ bounded by $\upsilon$ with
ideal point $x$. Then outside of a compact part
it projects very near a Reeb annulus in $\pi(F)$ and so 
this tail of $V$
has no singularity of the foliation $\wlse$.
By Euler characteristic reasons it follows that
the interior of 
$V$ has no singularities in the compact part of $V$ also.
In fact the arguments show that the tail of $V$ flows
into the {\underline {interior}} of a Reeb band 
in a nearby leaf $U$ of $\fn$.
Then leaves of $\wlse$ near $\upsilon$ on the outside of
$V$ will also flow to interior of the Reeb band  in $U$.
Therefore there are appropriate 2 rays on the outside 
so that they will be in the same leaf of $\widetilde \Lambda^s_U$
and hence in the same leaf of $\wls$.
It follows that in $E$ the leaf $\upsilon$ is also approximated
on the outside by a leaf which has a line leaf with 
both ideal points the same. This implies that $\upsilon$
has no singularities.

So far we proved the following:

\vskip .1in
\noindent
{\bf {Conclusion}} $-$ Let $\upsilon$ be a slice of 
$\wlse$ with 
two rays converging to the same ideal point $x$ of
$\pin E$  and $V$ is the region of $E$ bounded by $\upsilon$.
Then $\upsilon$ has no singularity of $\wlse$ and neither does $V$.
Also $\pi(V)$ is either contained in or asymptotic
to a Reeb annulus in a leaf of $\fol$ and so $E$
is asymptotic to a Reeb band in a leaf $F$ in
the direction $x$.
\vskip .1in

\begin{figure}
\centeredepsfbox{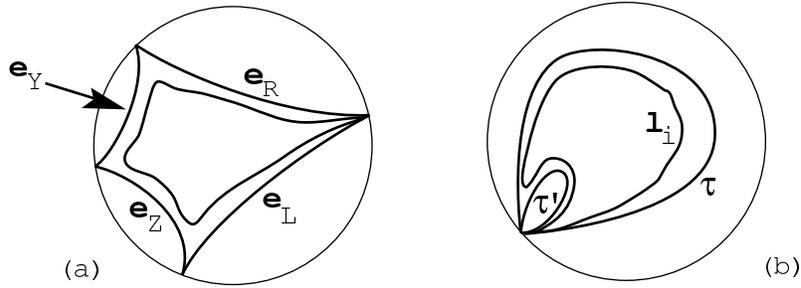}
\caption{
a. $l_i$ converging to non separated
leaves $e_L, e_Z, e_Y, e_R$ of $\wlse$, b. Nested families
and identifications of ideal points.}
\label{nonsep}
\end{figure}

\noindent
{\underline {Continuation of the proof of proposition \ref{band}}}.

What we want to prove is that in $E$ itself the region
$V$ is contained in the interior of a spike region.
Notice it is not true in general that $\pi(V)$ is
contained in a Reeb annulus, only that it is
asymptotic to a Reeb annulus. For instance start
with a leaf of $\fol$ having a Reeb annulus and blow that
into an I-bundle. Then produce holonomy associated
to the core of the Reeb annulus. Then one produces
Reeb bands asymptotic to but not contained in Reeb annuli.

Since $V$ is asymptotic to the Reeb band in $F$, it turns
out that (after rearranging by covering translations) that
$E$ intersects both $L$ and $R$ leaves of $\wls$.
Their intersection produces two leaves $e_L, e_R$ of 
$\wlse$ which are not separated from each other and which
have the same ideal point $x$. There are then leaves
$l_i$ of $\wlse$ all with ideal point $x$ and which 
converge to $e_L \cup e_R$. This follows from
the fact that in $F$ the same is true and $E$ is
asymptotic to $F$ in that direction, plus the connectivity
of the intersection of $E$ with leaves of $\wls$.

Now the sequence $l_i$ can converge to other
leaves as well, all of which will be non separated
from $e_L, e_R$. 
The set of limits is an ordered set and the any other 
leaf is 
between $e_L$ and $e_R$.
By theorem \ref{theb}
there are only finitely many of them.
We refer to
fig. \ref{nonsep}, a, where for simplicity we consider there
are 4 leaves in the limit: $e_L, e_Y, e_Z, e_R$ contained
in leaves $L, Z, Y, R$ of $\wls$. These leaves of $\wls$
are non separated from each other and form an ordered set.
Let $\xi$ be the region of
$E$ which is the  union of the region bounded by all the $l_i$
plus the boundary leaves, which are non separated from $e_L, e_R$.
Clearly every leaf in the interior has ideal point $x$ and
has no singularity. We want to show that $\xi$ is
a spike region.

Any two consecutive leaves 
of $\partial \xi$ in this ordering will have rays
with same ideal point and leaves $l_i$ converging to them.
This situation is important on its own and is analysed in
the following proposition:

\begin{proposition}{}{}
Suppose $v_1, v_2$ are non separated leaves in $\wlsg$ for some $G$ 
leaf of $\fn$.
Suppose there are no leaves non separated from $v_1, v_2$ in
between them. Then
the corresponding rays 
of $v_1, v_2$ 
have the same ideal point in $\pin G$. In addition
they are a bounded distance from a geodesic ray of $G$ with same ideal point.
In $M$ this region either projects to or is asymptotic to a Reeb
annulus.
\label{nh}
\end{proposition}

\begin{proof}{}
We do the essentially the
same proof as in the case of leaves of $\wlsg$
with same ideal points, except that we
go in the direction of the non Hausdorfness. 
Because there are no non separated leaves in between
$v_1, v_2$, then the corresponding rays have the same ideal
point.
Choose $w_i, y_i$ in these rays of $v_1, v_2$ and escaping
towards the ideal point and so that $d_G(w_i,y_i)$ is
less than $4 \delta_0$.
We do the limit analysis 
using $f_i(w_i), f_i(y_i)$ converging in $\mi$. Because 
$v_1, v_2$ are non separated it follows that $f_i(w_i), f_j(w_j)$
are in the same stable leaf (of $\wls$) for $i, j$ big enough.
Hence we can readjust so that they are all in the same
stable leaf and similarly for $f_i(y_i)$. The same arguments
as before show that that region of $G$ between $v_1, v_2$ projects
in $M$ to set in a leaf of $\fol$
which is either contained in or asymptotic to
a Reeb annulus. The results follow.
In general nothing can be said about the other direction
in the leaves
$v_1, v_2$: in particular it does not follow at all 
that the other rays of $v_1, v_2$ have to have
the same ideal point.
\end{proof}

Given this last proposition then for any two consecutive
rays in $\partial {\cal E}$ it follows that they are a bounded
distance from a geodesic ray in $E$. 
All that is needed to show that ${\cal E}$ is
a spike region is to prove that the ideal
points of the rays in the boundary are distinct except
for consecutive rays.

Suppose there are other identifications of ideal points of
leaves in the boundary of ${\cal E}$. Then there is at least one
line leaf $\tau$ in the boundary of ${\cal E}$ so that
$\tau$ has identified ideal points.
Our analysis so far shows that $\tau$ is in the interior
of another region similar to the one constructed above
so that all leaves have just one common ideal point.
Since the $l_i$ limit on $\tau$, then the ideal
point of $\tau$ has to be $x$. In addition the leaves in
this new region have to be nested. But if the $l_i$ together
with $\tau$ are a nested family of leaves of $\wlsf$, 
then the $\tau$ is outside the $l_i$ hence the region
in $E$ bounded by $\tau$ enclosed the whole region
${\cal E}$, see fig. \ref{nonsep}, b.
There is at least one other leaf $\tau'$ in
$\partial {\cal E}$. The same arguments we used for $\tau$
can be applied to $\tau'$. But it is impossible that
the $l_i$ are also nested with the $\tau'$, see fig. 
\ref{nonsep}, b.

This shows that the ideal points of ${\cal E}$ are
distinct except as mandated by consecutive rays. In addition
any line leaf in the boundary of ${\cal E}$ has distinct
ideal points and rays which are a bounded distance from
geodesic rays. It follows that the whole leaf is
a bounded distance from a geodesic in $E$.
This shows that 
${\cal E}$ is a spike region.
This finishes the proof of proposition \ref{band}.
\end{proof}

Finally in the case $\wls$ has Hausdorff leaf space
one can say much, much more about metric properties
of leaves of $\wlsf$:

\begin{proposition}{}{}
Suppose that $\Phi$ is an almost pseudo-Anosov flow transverse
to a foliation $\fol$ with hyperbolic leaves.
Suppose that $\wls$ has Hausdorff leaf space.
Then there is $k_0 > 0$
so that for any $F$ leaf of $\wls$, then the slice leaves 
of $\wlsf$ are uniform $k_0$-quasigeodesics.
\end{proposition}

\begin{proof}{}
If there is a leaf $F$ of $\fn$ and a slice leaf of $\wlsf$
with only one ideal point, then the proof of proposition
\ref{band} shows that there are leaves of $\wls$ non
separated from each other. This is impossible.

Suppose now that for any integer $i$, there are 
$x_i$ in $\mi$, $x_i$ in leaves $F_i$ of $\fn$ with
$x_i$ in line leaves $l_i$ of $\widetilde {\Lambda}^s_{F_i}$
with distance
from $x_i$ to $l^*_i$ in $F_i$ going to infinity.
Here $l^*_i$ is the geodesic in $F_i$ with same ideal
points as $l_i$.
Up to covering translations assume $x_i$ converges
to $x$.
Also assume all $x_i$ are in the same sector 
of $\wls$ defined by $x$. Since $l_i$
converges to $l$, the arguments
in lemmas \ref{encon} and \ref{cont} would show that
the ideal points of $l$ are the same. 
This was just disproved above.

Given that, the line leaves are within some global distance
$a_0$ 
of the respective
geodesics in their leaves.
It is well known that these facts imply that
the slice leaves of $\wlsf$ are uniform quasigeodesics.
For a proof of this well known fact see 
for example \cite{Fe-Mo}.
\end{proof}

\section{Continuous extension of leaves}

The purpose of this section is to prove the main theorem: \
the continuous extension property for leaves of foliations
which are almost transverse to quasigeodesic singular pseudo-Anosov
flows in atoroidal $3$-manifolds. As seen
before this implies that $M$ has negatively curved
fundamental group.

Suppose first that $\Phi$ is an almost pseudo-Anosov flow
which is transverse to a foliation $\fol$ 
with hyperbolic leaves
in a
general closed $3$-manifold $M$.
Given a leaf $F$ of $\fn$
we introduce geodesic ``laminations" in $F$  
coming from $\wlsf, \wluf$.
We only work with the stable foliation, similar results hold for the
unstable foliation.
Assume that a leaf $l$ of $\wlsf$ is not singular. If both ideal
points are the same let $l^*$ be empty. Otherwise let $l^*$ be the geodesic
with same ideal points as $l$. If $l$ is singular, then no line leaves
of $l$ 
have the same ideal point by proposition \ref{band}. For each line
leaf $e$ 
of $l$
let $e^*$ be the corresponding geodesic and $l^*$ their
union. Let now $\tau^s_F$ be the union of these geodesics
of $F$.
Leaves of $\wlsf$ do not have transverse intersections and
therefore the same happens for leaves of $\tau^s_F$.

Suppose that $\wlsf$ has non separated leaves 
$l, v$ which are not in the boundary of a spike region.
Then there are $l_i$ converging to $l \cup v$ (and maybe
other leaves as well),
but $l^*_i$ 
does not converge to $l^*$ or $v^*$. Notice none of the limit leaves
can have identified ideal points, because then they
would be in the interior of a spike region (proposition 
\ref{band}) and have a neighborhood which is product
foliated.
Let $\otsf$ be the closure
of $\tau^s_F$. Then $\otsf$ is a geodesic lamination in $F$.
Similarly define $\tau^u_F$, $\overline \tau^u_F$.
In a complementary region $U$ 
of $\otsf$ associated to non Hausdorffness, there
is one boundary component which is added (a leaf of $\otsf - \tau^s_F$)
and which is the limit of the $l^*_i$ as above.
All of the other boundary
leaves of the region
are associated to the non separated leaves
of $\wlsf$ and are in $\tau^s_F$.

\begin{lemma}{}{}
The new leaves in $\otsf$ (that is those
in $\otsf$ - $\tau^s_F$) come from non Hausdorff
pairs $(l,v)$ 
of $\wlsf$
as in the description above.
\end{lemma}

\begin{proof}{}
Let $e_i$ in $\tau^s_F$ converging to $e$ not in $\tau^s_F$.
Then choose $l_i$ line leaves in $\wlsf$ with
$e_i = l^*_i$. Given $u$ a point in $e$, there is
$u_i$ in $l^*_i$ very close to $u$. Then there are
$p_i$ in $l_i$ which are $2 \delta_0$ close to $u_i$. 
Up to subsequence
assume that $p_i$ converges to $p_0$ and
let $l$ be the line leaf of $\wlsf$ that the sequence
$l_i$ converges to. 
Since the $l^*_i$ converges to $e$ which is not in
$\otsf$ and $l^*$ is in $\otsf$, it follows that 
$l^*_i$ does not converge to $l^*$.
By lemma  \ref{encon} this is associated to 
a non Hausdorff situation:
$l_i$ converging to $l$ and other leaves as well and $l^*$ is
the added leaf associated to this non Hausdorfness. This finishes
the proof of the lemma.
\end{proof}

\begin{lemma}{}
The complementary regions of $\otsf$ are 
ideal polygons associated to singular leaves and non
Hausdorff behavior of $\wls$.
If $M$ is atoroidal 
then these regions
are finite sided ideal polygons.
\end{lemma}

\begin{proof}{}
Let $x$ be in a complementary region $U$ of $\otsf$.
Let $e$ be a leaf in the boundary $\partial U$.
Let $I$ be the interval of $\pin F - \partial e$ containing
other ideal points of $U$.
Suppose first that $e$ is an actual leaf of $\tau^s_F$,
which comes from a line leaf $l$ of $\wlsf$.
It may be that $l$ is contained in a singular leaf $z$ of
$\wlsf$ which
is singular on the $x$ side.  This means that $z$ has ideal
points in $I$.
In that case $x$ is in the complementary region obtained
by splitting $z$. This region must be $U$.
Otherwise $l$ is not singular on the side containing
$x$ and we may assume there
are $l_i$ leaves of $\wlsf$
with ideal points in the closure of $I$ in $\pin F$,
with $l_i$ converging to $l$. If
the ideal points of $l_i$ converge to that of $l$ then eventually
$l^*_i$ separates $x$ from $e$ and $x$ is
not in the complementary region $U$ $-$ impossible. 
Hence the ideal points of $l_i$ do not
converge to $\partial e$ and there is non Hausdorfness and
a complementary region in that side of $l$. Then $x$ needs to be
in this complementary region (which is $U$) and $e$ is a boundary
leaf of $U$ which comes from a line leaf of $\tau^s_F$. 

Suppose now that $e$ is an added leaf.
There are $l_i$ 
leaves of $\wlsf$ with $e_i = l^*_i$ converging to $e$ on the side
opposite to $x$, otherwise $x$ is not in $U$.
Then $l_i$ converges to more than one leaf of $\wlsf$ 
producing non Hausdorff behavior
and a complementary region with $e$ in its boundary.
The $x$ is in the region
associated to this non Hausdorff behavior, so
the complementary region must be $U$. 

If there is a complementary region of $\otsf$ with infinitely
many sides then it is associated to non Hausdorff behavior
and so there are leaves $l_i$ of $\wlsf$ converging to infinitely
many distinct leaves of $\wlsf$.
Then there is $L$ leaf of $\wls$ which is non separated from
infinitely many other leaves. Theorem \ref{theb} implies that
there is a ${\bf Z} \oplus {\bf Z}$ subgroup
of $\pi_1(M)$, contradiction.
This finishes the proof.
\end{proof}

We now turn to the continuous extension property.

\begin{theorem}{(Main theorem)}{}
Let $\fol$ be a foliation in $M^3$ closed, atoroidal.
Suppose 
that $\fol$ is almost transverse to a 
quasigeodesic, singular pseudo-Anosov flow $\Phi_1$
and transverse to an associated almost pseudo-Anosov flow $\Phi$.
Singular means $\Phi_1$ is not a topological Anosov flow.
Then for any leaf $F$ of $\fn$, the inclusion map
$\Psi: F \rightarrow \mi$ extends to a continuous map

$$\Psi: \ F \cup \pin F \ \rightarrow \ \mi \cup \si$$

\noindent 
The map
$\Psi$ restricted to $\pin F$, gives a continuous parametrization
of the limit set of $F$, which is then locally connected.
\label{exten}
\end{theorem}

\begin{proof}{}
The hypothesis imply that $\pi_1(M)$ is negatively curved.
Difficulties in the proof  of this result are that
$\wls, \wlu$ may have non Hausdorff leaf space \cite{Mo5,Fe6}
and so 
$\wlsf, \wluf$  can have non Hausdorff leaf space. This
implies
that the leaves of $\wlsf, \wluf$ 
cannot be uniform quasigeodesics in $F$.
In addition the leaves of $\wls, \wlu$ are not quasi-isometrically
embedded in $\mi$.
The proof is done in two steps: first we define an extension and then
we show that it is continuous.

The proof will fundamentally use the fact that $\Phi_1$ is
a quasigeodesic pseudo-Anosov flow.
It was proved in \cite{Fe-Mo} that this implies
that $\Phi$ is a quasigeodesic flow as well.
From now on we use the stable/unstable foliations
$\wls, \wlu$ of $\wwp$.
First we need to review some facts about quasigeodesic 
almost pseudo-Anosov flows. 
If $\gamma$ is an orbit of $\wwp$ then it is
a quasigeodesic and hence has unique distinct
ideal points $\gamma_-$ and $\gamma_+$ 
in $\si$ corresponding to
the positive and negative flow directions \cite{Th1,Gr,Gh-Ha,CDP}.
Hence given $x$ in $\mi$ define

$$\eta_+(x) \ = \ \gamma_+, \ \ \ \eta_-(x) = \ \gamma_-,
\ \ \ \eta_+(x) \ \not = \ \eta_-(x),$$

\noindent
where $\gamma$ is the $\wwp$ flowline through $x$.
If $L$ is a leaf of $\wls$ or $\wlu$ and $a$ is a limit point
of $L$ in $\si$, then there is an orbit $\gamma$ of $\wwp$ 
contained in $L$ with either $\gamma_- = a$ or $\gamma_+ = a$,
that is, any limit point of $L$ is a limit point of one
of its flow lines \cite{Fe2,Fe6}.
Also any such $L$ in $\wls$ is Gromov negatively curved
\cite{Gr,Gh-Ha,Fe2,Fe6} and
has an intrinsic ideal boundary $\partial L$ 
consisting of a single forward ideal point and distinct negative
ideal points for each flow line
\cite{Fe2,Fe6}. The set $L \cup \pin L$ is a natural
compactification of $L$ in the Gromov sense.
For instance if $L$ is a non singular leaf,
then $L \cup \pin L$ is a closed disk. In this
case the foliation by flow lines in $L$ is equivalent
to the foliation in $\hh$ by geodesics sharing a fixed
point in $\su$.

A very important fact for us is that  the inclusion

$$\kappa: L \rightarrow \mi \ \ \ \ 
{\rm extends \  to \ a \ continuous \
map} \ \ \ \  \kappa: L \cup \partial L \rightarrow \mi \cup \si,$$

\noindent
\cite{Fe2,Fe6}.
This all follows from the fact that $\Phi$ is quasigeodesic.
This works for any $L$ in $\wls$ or $\wlu$.
If $L$ is in $\wls$
there is a unique distinguished ideal point in $\si$ denoted
by $L_+$
which is the forward limit point in $\si$ of any flow line
in $L \subset \mi$.
If in addition $\Lambda^s$ is a quasi-isometric singular
foliation, then the extension $\kappa$ is always a homeomorphism
into its image, but this is not true if $\Lambda^s$ is
not quasi-isometric.
Similarly for $L$ in $\wlu$.

Throughout the proof we fix a unique identification
of $\mi \cup \si$ with the closed unit ball in
${\bf R}^3$.
The Euclidean metric in this ball induces the visual distance
in $\mi \cup \si$.
Then $diam(B)$ denotes the diameter in this distance for
any subset $B$ of $\mi \cup \si$.
A notation used throughout here is the following: 
if $A$ is a subset of a leaf $F$ of $\fn$, then $\overline A$
is its closure in $F \cup \pin F$. 

We now produce an extension $\Psi: \pin F \rightarrow \si$.

\vskip .1in
\noindent
{\underline {Case 1}} $-$ Suppose that $v$ in $\pin F$
 is not an ideal point
of a ray in $\wlsf$ or in $\wluf$.

Since $\pi_1(M)$ is negatively curved, then
complementary regions of $\otsf$ are finite sided ideal
polygons.
Hence there are $e_i$ in $\otsf$ 
so that $\{ e_i \cup \partial e_i \}, \ i \in {\bf N}$ define
a neighborhood basis of $v$ (in $F \cup \pin F$) 
and $\{ e_i \}$ forms a nested
sequence.
Here $\partial e_i$ are the ideal points of $e_i$ in $\pin F$.
We say that the $\{ e_i \}$ define a neighborhood basis at $v$.
Assume that no two $e_i$ share an ideal point $-$
possible because of hypothesis.
If $e_i$ is in $\otsf - \tau^s_F$ then it is the limit of leaves
in $\tau^s_F$ and by adjusting the sequence above we can
assume that $e_i$ is always in $\tau^s_F$.
Let $l_i$ in $\wlsf$ with $l^*_i = e_i$
and $L_i$ leaves of $\wls$ with $l_i \subset L_i$.

Similarly there are $c_i$ in $\otuf$ defining a neighborhood
basis of $v$.
Up to subsequence we may assume that $e_1, c_1, e_2, c_2$, etc..
are nested and none of them have any common ideal points
(in $F \cup \pin F$)
and $c_i$ is in $\tau^u_F$.
Let $b_i$ in $\wluf$
with $b^*_i = c_i$ and $B_i$ leaves of $\wlu$ with $b_i \subset B_i$.

At this point we need the following result:

\begin{lemma}{}{}
Let $L$ leaf of $\wls$, $B$ leaf of $\wlu$ and $F$ leaf of $\fn$
so that $F$ intersects both $L$ and $B$: $l = L \cap F,
\ b = F \cap B$. Suppose that $b$ and $l$ are disjoint
in $F$. Then $L$ does not intersects $B$ in $\mi$.
\end{lemma}

\begin{proof}{}
Suppose not. Recall that $\Theta(L), \Theta(B)$ are
finite pronged, non compact
trees and they intersect in a compact subtree.
The union is also a finite pronged tree.
In addition $\Theta(L \cap B)$ is connected.
The sets $\Theta(l), \Theta(b)$ are disjoint in this union.
Let $x$ be a boundary point of $\Theta(l)$ which is
either in $\Theta(L \cap B)$ or separates $\Theta(L \cap B)$
from $\Theta(l)$ in this union, see fig. \ref{inter}, a.

\begin{figure}
\centeredepsfbox{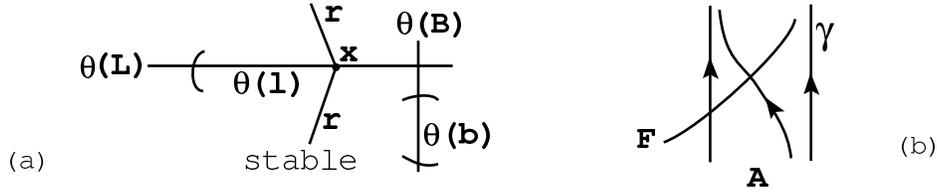}
\caption{
a. Obstruction to intersections of leaves,
b. The case of $F$ escaping up.}
\label{inter}
\end{figure}

Let $\gamma = x \times \rrrr$, an orbit of $\wwp$.
The first possibility is that $F$ escapes up
as $\Theta(F)$ approaches $x$.
Then $\gamma$ is a repelling orbit with respect to the
$\Theta(l)$ side, see fig. \ref{inter}, b
and $\gamma$ is in the boundary of a lift annulus
$A$.
This means that $\Theta(l)$ is a generalized unstable prong
from the point of view of $x$.
By proposition \ref{bounda} there is a stable slice
$r$ of $\oos(x)$ with $r$ contained in $\partial \Theta(F)$
and $F$ escapes up as $\Theta(F)$ approaches $r$, see
fig. \ref{inter}, a.
 The two sides of $r$ are the closest
generalized prongs to $\Theta(l)$ on either side of $\Theta(l)$.
This implies that $r$ separates $\Theta(b)$ from $\Theta(F)$
see fig. \ref{inter}, a.
Then $\Theta(b)$ cannot be contained in $\Theta(F)$, contradiction.

The second option is that $F$ escapes down as $\Theta(F)$ 
approaches $x$ along $\Theta(l)$.
Here there is a slice $r$ of $\oou(x)$ with $r$ contained
in $\partial \Theta(F)$ and the closest to $\Theta(l)$ 
on both sides of $\Theta(l)$. 
Either $\Theta(b) \subset r$ or $r$ 
separates $\Theta(b)$ from $\Theta(F)$. In any case
$\Theta(b)$ does not intersect $\Theta(F)$, again
a contradiction.
This finishes the proof of the lemma.
\end{proof}

\vskip .1in
\noindent
{\underline {Claim}} $-$ Both $L_i$ and $B_i$ escape in $\mi$.

Notice $e_i \cap c_j = \emptyset$ for any $i, j$.
If $l_i \cap b_j$ is non empty with $j > i$, then the nesting
property above implies that $b_{i+1}, b_{i+2}, ..., b_j$
all have to intersect. Since there is a global upper
bound on the number of prongs of leaves of $\wls, \wlu$, this
can happen for only finitely many times. Up to taking a 
further subsequence we may assume that all the $l_i, b_j$
are disjoint.

The lemma shows that 
$L_i \cap B_j$ is empty for any $i, j$, and they
form nested sequences of leaves in $\mi$.
Suppose that the sequence $\{ L_i \}$ does not escape
compact sets. Then there is $L$ in $\wls$ which is a limit
of $L_i$ (and possibly other leaves as well).
Let $\alpha$ be an orbit in $L$ which is not in
a lift annulus. Then $\wu(\alpha)$ is transverse
to $L$ in $\alpha$ and hence intersects $L_i$ for
$i$ big enough. Since the $L_i, B_j$ are nested
this would force $\wu(\alpha)$ to intersect
$B_j$ for $j$ big enough, contradiction.
It follows that both $L_i$ and $B_j$ escape compact sets
as $i, j \rightarrow \infty$.

Let $r$ be a geodesic ray in $F$ with ideal point $v$.
For each $i$, there is a subray of $r$ contained in
the component of $F - l_i$ which is in a small neighborhood
of $v$. Hence $\Psi(r)$ has a subray which is contained 
in the corresponding component $V_i$ of $\mi - L_i$.
These components $V_i$ form a nested sequence. The
ray $\Psi(r)$ can only limit in the limit set of $V_i$.
We need the following lemma which will be 
a key tool 
throughout the proof.

\begin{lemma}{(basic lemma)}{}
Let $\{ Z_i \}$ be a sequence of leaves or line leaves
or slices or any flow saturated sets in 
leaves of either in $\wls$ or
$\wlu$ (not all leaves $Z_i$
need to be in the same singular foliation).
If the sets $Z_i$ escape compact sets in $\mi$, then
up to taking a subsequence $\overline Z_i$ converges 
to a point in $\si$.
\label{narrow}
\end{lemma}

\begin{proof}{}
Let $Y_i$ be the leaf of $\wls$ or $\wlu$ which contains $Z_i$.
Up to subsequence assume $Y_i \in \wls$.
The statement is equivalent to $diam(Z_i)$ converges to $0$.
Otherwise up to subsequence we can assume $diam(Z_i) > a_0$
for some $a_0$ and all $i$ and hence no subsequence
can converge to a single point in $\si$.
Then there is $p_i$ in $Z_i$ with visual distance from
$p_i$ to $(Y_i)_+$ is bigger than $a_0/2$.
Notice that $(Y_i)_+$ is a point in $\overline Z_i$.
Let $\gamma_i$ the orbit of $\wwp$ through $p_i$.
If $(\gamma_i)_-$ is very close to $(\gamma_i)_+ = (Y_i)_+$
then the geodesic with these ideal points has
very small visual diameter. Since $\gamma_i$ is a global
bounded distance from this geodesic \cite{Gr,Gh-Ha,CDP},
the same
is true for $\gamma_i$ contradiction to the choice of $p_i$.
Hence the geodesic above intersects a fixed
compact set in $\mi$ and so does $\gamma_i$. This
contradicts the fact that $Z_i$ escape compact
sets in $\mi$ and finishes the proof.
\end{proof}

We claim that the limit sets of $V_i$ above shrink to a single
point in $\si$.
The limit sets form a weakly monotone decreasing sequence, because
the $L_i$ are nested and so are the $V_i$.
If the limit set does not have diameter going to zero,
then there are points in the limit set of $L_i$ which
are at least $ 2 \delta_1$ apart for some fixed
$\delta_1 > 0$.
By the previous lemma the $L_i$ cannot escape 
compact sets in $\mi$, contradiction.
Since the limit sets of $V_i$ shrinks to a point in $\si$, let
$\Psi(v)$ be this point. Clearly $\Psi(r)$ limits
to this point and so does $\Psi(r')$ for any other geodesic
ray $r'$ in $F$ with ideal point $v$.

\vskip .15in
\noindent
{\underline {Case 2}} $-$ Suppose that $v$ is an ideal point
of a leaf of $\wlsf$ or $\wluf$.

Let $l$ be a ray in say $\wlsf$ which limits on $v$
and $r$ a geodesic ray on $F$ with ideal point $v$.
Then $l$ is contained in $L$ leaf of $\wls$. 
Either $\Theta(l)$ escapes in $\Theta(L)$ or limits
to a point $x$ in $\Theta(L)$.

Consider the first case. Then in the intrinsic geometry of
$L$, the ray $l$ converges to the positive ideal point
of $L$, hence in $\mi \cup \si$, the image $\Psi(l)$ converges
to $L_+$.
In the other option let $\beta = x \times \rrrr$,
an orbit of $\wwp$.
As $l$ escapes in $F$ then in $L$ it either escapes up or
down. If it escapes down then it converges to the
negative ideal point of $\beta$ in $L \cup \pin L$
and hence $\Psi(l)$ converges to $\beta_-$. Otherwise
$l$ escapes up in $L$ as $\Theta(l)$ approaches
$x$. In this case $\beta$ is in the boundary of
a lift annulus and $l$ converges to the positive ideal
point in $L \cup \pin L$ and so $\Psi(l)$ converges to $L_+$ 
again.

The remaining case is that $\Theta(l)$ escapes in $\Theta(L)$.
Then as seen in $L \cup \pin L$ the ray $l$ converges
to the positive ideal point $p$ of all flow lines in $L$.
Hence $\Psi(l)$ converges to $\kappa(p) = L_+$.
Let $\Psi(v)$ be the limit point in any case.
Similarly if $l$ is a ray of $\wluf$.

Every point in $r$ it is $2 \delta_0$ close
to a point in $l$ in $F$, hence the limit of $\Psi(r)$ in
$\mi \cup \si$ is the same as that of
$l$. If $l'$ is another ray 
of $\wlsf$ or $\wluf$ converging to $v$,
then it will have
points boundedly close to $r$ which escape in 
$l'$ and therefore $\Psi(l')$ 
has the same ideal point
in $\si$. 
Therefore $\Psi(v)$ is well defined.

This finishes the construction of the extension
of $\Psi$ to $\pin F$.

\vskip .1in
\noindent
{\underline {Proof of continuity of the extension}} $-$

\vskip .1in
\noindent
{\underline {Case 1}} $-$  $v$ is not an ideal point of 
a ray in $\wlsf$ or $\wluf$.

Let $r$ be a geodesic ray in $F$ with ideal point $v$.
Recall the extension construction. 
There are $l_i$ in $\wlsf$ shrinking to $v$
in $F \cup \pin F$ and similarly $b_i$ in $\wluf$, 
assumed to be nested with the $l_i$.
Let $\{ l^*_i \}$ define a neighborhood basis of $v$ in
$F \cup \pin F$.
Let $L_i$ in $\wls$ with 
$l_i \subset L_i$, and $b_i \subset B_i \in \wlu$ 
as in the construction
case 1. Then as seen in the construction, the 
$L_i, B_i$ escape in $\mi$.
Let $U_i$ be the component of $F - l_i$ containing a subray of
$r$ and $V_i$ the component of $\mi - L_i$ containing $U_i$.
Notice that $\Psi(U_i) \subset V_i$.
Let now $z$ in $\overline U_i$ with the closure taken in 
$F \cup \pin F$ and $\overline V_i$ the closure of $V_i$
in $\mi \cup \si$.
Then $\overline U_i$ is a neighborhood of $v$ in $F \cup  \pin F$.
If $z$ is in $\Psi(\overline U_i)$ then using either of the
constructions in the extension part shows that $z$ is
a limit of points in $\Psi(U_i) \subset V_i$.
As seen in the construction arguments
the diameter of $\overline V_i$ in the visual distance
is converging to $0$. Hence we obtain continuity of $\Psi$ at
$v$.
This finishes the proof in this case.

\vskip .1in
\noindent
{\underline {Case 2}} $-$ $v$ is an ideal point of 
a ray of $\wlsf$ or $\wluf$.

This case is considerably more complicated, with several possibilities.

\vskip .1in
\noindent
{\underline {Case 2.1}} $-$
$v$ is an  ideal point of $\wlsf$ but not of $\wluf$
(or vice versa).

Suppose the first option occurs.
There is $l$ ray in $\wlsf$ with ideal 
point $v$. 
We may assume that $l$ is not in a leaf of $\wlsf$ with
same ideal points. Otherwise we can choose $l$ to be one
of the boundary leaves of the corresponding spike region.
Since $v$ is not an ideal point of $\wluf$,
there are $g_i$ line leaves
in $\wluf$ defining a basis neighborhood system at $v$.
Let $g_i$ be contained in $G_i$ leaves of $\wlu$.
Let $L$ in $\wls$ containing $l$.
If $G_i$ escapes in $\mi$ as $i \rightarrow \infty$, then as
seen in case 1, we are done.
Let then $G_i$ converge to the finite set of leaves

$${\cal V} \ =  \ H_1 \cup H_2 .... \cup H_m \ \ \ {\rm leaves
\ of } \ \ \wlu$$

\noindent
We can assume that $G_i \cap l \not = \emptyset$ for all $i$,
$G_i$ is non singular and the sequence $\{ G_i \}, i \in {\bf N}$
is nested.

\vskip .1in
\noindent
{\underline {Case 2.1.1}} $-$ Suppose that $L$ intersects
${\cal V}$, say $L \cap H_1 \not = \emptyset$.

Then $l$ escapes down as $\Theta(l)$ 
approaches $\Theta(L \cap H_1)$.
Otherwise $L \cap H_1$ is in the boundary of a lift annulus
$A$ and $l$ has a subray contained in this lift annulus.
But then $A$ is also contained in the unstable leaf $\wu(L \cap H_1)$
and so $G_i$ cannot intersect $l$, contradiction.
As $l$ escapes down in $L$, then
the ideal point of $\Psi(l)$ is
$(L \cap H_1)_-$ which is equal to $(H_1)_-$, the negative
ideal point of $H_1$. 

Since the values of $\Psi(p)$ for $p$ in $\pin F$ are obtained
as limits of values in $\Psi(F)$, then we only need to show that
if $z_k$ is in $F$ and $z_k$ converges to $p$ as $k \rightarrow \infty$,
then $\Psi(z_k)$ converges to $\Psi(p)$. Suppose this is not
the case.


By taking a subray if necessary, we may assume that $l$ 
does not intersect a lift annulus and hence it is transverse
to the unstable foliation $\wluf$ in $F$.
Parametrize the leaves of $\wluf$
intersected by $l$ as $\{ g_t, \  t \in \rrrr_+ \}$,
contained in $G_t \in \wlu$
(by an abuse of notation think of the $G_i$ as
a discrete subcollection of the $G_t, t \in \rrrr_+$). Let 

$${\cal U}  \ = \ \bigcup_{t>0} G_t$$

\noindent
No $g_t$ (or leaf of $\wluf$) 
has ideal point $v$ in $\pin F$.
As $\{ g_i \}, i \in {\bf N}$ converges to $v$ in $F \cup \pin F$, then
$g_t$ escapes compact sets in $F$ as
$t \rightarrow \infty$ and the ideal points of $g_t$ converge
to $v$ on either side of $v$. 
Up to subsequence assume that all of the elements
of the sequence $\{ z_k \}$ 
are either entirely contained
in ${\cal U}$ or disjoint from ${\cal U}$.

\vskip .1in
\noindent
{\underline {Situation 1}} $-$ Suppose that $z_k$ is not in ${\cal U}$ 
for any $k$.

Since $z_k$ is very close to $v$ in the compactification
$F \cup \pin F$
and $g_t$ converges to $v$ in $F \cup \pin F$ when $t \rightarrow \infty$,
then there are $t, s$
with $z_k$ between $g_t$ and $g_s$ (in $F$). Notice $z_k$ is not in
any of them. Now there is a unique time $t_k$ so that exactly at that
time $\Psi(z_k)$ switches from being in one side of $G_t$  in $\mi$
to the other
(equivalently compare the $z$ and $g_t$ in $F$).
In particular, either there is a line leaf $L_{t_k}$
of $G_{t_k}$ which separates $\Psi(z_k)$ from all the other 
$G_t$, see fig. \ref{sepa}, a, or there is a leaf $L_{t_k}$
non separated from $G_{t_k}$ with $\Psi(z_k)$ either in $L_{t_k}$
or $L_{t_k}$ separates $\Psi(z_k)$ from all $G_t$,
see fig. \ref{sepa}, b. This can be seen
in the leaf space of $\wlu$, which is a non Hausdorff tree 
\cite{Fe8,Ga-Ka,Ro-St}.

\begin{figure}
\centeredepsfbox{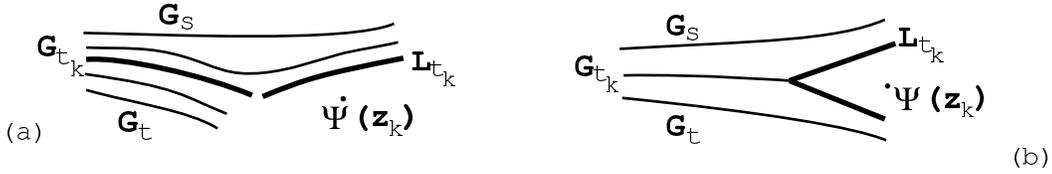}
\caption{
a. Line leaf separating points, b. Non separated
leaf separating points.}
\label{sepa}
\end{figure}

\vskip .1in
\noindent
{\underline {Claim}} $-$ In the Gromov-Hausdorff topology of
closed sets of 
$\mi \cup \si$, the sets $\overline L_{t_k}$ 
converge to $(H_1)_-$ as
$k \rightarrow \infty$.

If $L_{t_k}$ is a line leaf of 
$G_{t_k}$, then $(L_{t_k})_- = (G_{t_k})_-$. If
$L_{t_k}$ is not separated from $G_{t_k}$ then also
$(L_{t_k})_- = (G_{t_k})_-$. 
This is because there are $E_i$ leaves of $\wlu$ with
$E_i$ converging to $L_{t_k} \cup G_{t_k}$.
So there are $x_i, \ y_i$ in $E_i$ with 
$x_i \rightarrow x, \ y_i \rightarrow y$ and $x \in L_{t_k}, \
y \in G_{t_k}$. Then

$$\eta_-(x_i) \ \rightarrow \eta_-(x) \ = \ \eta_-(L_{t_k}), \ \ \ 
\eta_-(y_i) \ \rightarrow \eta_-(y) \ = \ \eta_-(G_{t_k})
\ \ \ \ {\rm and} \ \ \eta_-(x_i) \ = \ \eta_-(y_i).$$

\noindent
The last equality occurs
 because $x_i, y_i$ are in the same unstable
leaf $E_i$.
Therefore $(L_{t_k})_-$
converges to $(H_1)_-$ when $k \rightarrow \infty$. 
Suppose that $\overline L_{t_k}$ does not 
converge to $(H_1)_-$ in $\mi \cup \si$.
Since

$$(L_{t_k})_- \ \ {\rm converges \ to} \ \ (H_1)_-,$$

\noindent
then lemma \ref{narrow} shows that $L_{t_k}$ does not
escape compact sets in $\mi$. So it limits to some $u$ in
$\mi$ and up to subsequence  we may assume there 
$u_k$ in $L_{t_k}$ with $u_k$ converging
to $u$.
The first possibility is that the $L_{t_k}$ are subsets of 
%
%
the leaves $G_{t_k}$. This implies that $\wwp_{\rrrr}(u)$
is in the limit of the sequence of leaves $G_{t_k}$ (in $\mi$), so
it is contained in ${\cal V}$.
The second possibility is $L_{t_k}$ non separated from $G_{t_k}$ so
$L_{t_k}$ is between $G_{t_{k-1}}$ and
$G_{t_{k+1}}$ hence $u$ is again in the limit of the $G_t$ so
$u$ is in ${\cal V}$.
The leaves $H_j$ in ${\cal V}$ are non singular in the
side the $G_t$ are limiting on, 
so there is a neighborhood of $u$ on that side of $H_j$ which has
no singularities hence the $u_k$ will be in ${\cal U}$ for
$k$ big enough. This contradicts the hypothesis
in this case.

This shows that 
$\overline L_{t_k}$ converges to $(H_1)_-$ 
in $\mi \cup \si$.
Also $L_{t_k}$ 
either contains $\Psi(z_k)$
or separates it from a base point in $\mi$. It follows that
$\Psi(z_k)$ converges to $(H_1)_-$, which is what we wanted to prove.
This finishes the analysis in situation 1.

\vskip .1in
\noindent
{\underline {Situation 2}} $-$ For all $k$ assume that $\Psi(z_k)$ is
in ${\cal U}$.

Let $t_k$ with $\Psi(z_k)$ in $G_{t_k}$,
hence $z_k$ is in $G_{t_k} \cap F = g_{t_k}$. 
Then $(\Psi(z_k))_- = (G_{t_k})_-$ converges to $(H_1)_-$.
We want to show that $\Psi(z_k)$ converges to $(H_1)_-$.
Otherwise there is $q$  in $\si$
different from $(H_1)_-$ and a subsequence, still denoted by
$\Psi(z_k)$ so that
$\Psi(z_k)$ converges
to $q$.
As in the claim of situation 1 above, 
$\wwp_{\rrrr}(z_k)$ does not escape compact sets
in $\mi$ and there is $z$ in $\mi$, so that up to another
subsequence, we may 
assume that $\wwp_{\rrrr}(\Psi(z_k))$ converges to $\wwp_{\rrrr}(z)$.
Since $\Psi(z_k)$ is in $G_{t_k}$ then $z$ is
in ${\cal V}$, say $z$ is in $H_j$.
Let $p = \Theta(z)$.
At this point notice that $F$ does not intersect any 
leaf $H_i$ in ${\cal V}$.
If it did, say in $w$ then $F$ intersects the nearby leaves
$G_t$ (for any $t$ big enough) near $w$. This would
imply $F \cap G_t = g_t$ 
does not escape compact sets in $F$,
contradiction.
Therefore 
$p = \Theta(z)$ is in $\partial \Theta(F)$.
Let $x_k$ in $g_{t_k} \cap l$.
Then
$\Theta(x_k)$ converges to a point in $\Theta(H_1 \cap L)$.
There are segments $b_k$ in 
$F \cap G_{t_k} = g_{t_k}$ from $x_k$
to $z_k$. Then $\Theta(b_k)$ converges to a ray in
$\Theta(H_1)$
and a ray in $\Theta(H_j) \subset \oou(p)$ and possibly
other unstable leaves. Then there is a ray in $\oou(p)$
contained in $\partial \Theta(F)$. This implies that
$F$ escapes down as $\Theta(F)$ approaches this ray of $\Theta(H_j)$.
Hence $\Psi(z_k)$ is getting closer to $z_-$ which is $(H_j)_-$,
which is also equal to $(H_1)_-$.
This is what we
wanted to prove anyway.

This finishes the proof of case 2.1.1,
that is, when $L$ intersects ${\cal V}$.

\begin{lemma}{}{}
Let $A$ in $\wlu$, $B$ in $\wls$ satisfying:
there are
$R_i$ leaves of $\wlu$ intersecting $B$ 
with $R_i$ converging to $A$
and $R_i \cap B$ escaping compact sets in $B$.
Then $A_-$ is equal to $B_+$.
\label{conperf}
\end{lemma}

\begin{proof}{}
Since $R_i$ converges to $A$ then $(R_i)_-$ converges
to $A_-$.
Also $R_i$ intersects $B$ so
$(R_i)_- = (R_i \cap B)_-$. As $R_i \cap B$ escapes compact
sets in $B$ then in the intrinsic geometry of $B$,
the $R_i \cap B$ converges
to the positive ideal point of $B$. This implies that
$(R_i \cap B)_-$ converges to $B_+$. This implies the result.
\end{proof}

\noindent
{\underline {Case 2.1.2}} $-$ $L$ does not intersect 
${\cal V}$.

Then $\Theta(l)$ escapes in $\Theta(L)$ and so $\Psi(l)$ converges to
$L_+$. By the previous lemma,
this is also equal to $(H_1)_-$. From this point on,
the proof is the same
as in case 2.1.1.
This finishes the proof of case 2.1.

\vskip .2in
\noindent
{\underline {Case 2.2}} $-$ $v$ is an ideal point of both $\wlsf$ 
and $\wluf$.

\vskip .1in
\noindent
{\underline {Case 2.2.1}} $-$ For any ray $l$ of 
$\wlsf$ and $e$ of $\wluf$
with $l_{\infty} = e_{\infty} = v$, then $l$ does not intersect
$e$.

Let $l', e'$ be rays as above. 
We may assume that $l', e'$ do not have any singularities.
Parametrize the leaves of $\wlsf$
intersecting $e'$ as $\{ l_t, t \geq 0 \}$ 
where $l_t \cap e'$ converges to $v$
in $F \cup \pin F$
as $t$ converges to infinity.

Since $l'$ limits on $v$ and is disjoint from $e'$, then $l'$ 
is on a side
defined by $e'$. We will prove continuity of
$\Psi$ at $v$ from the other side
of $e'$. 
The point $p_t = l_t \cap e'$ disconnects $l_t$. For simplicity
we only consider those $l_t$ with $l_t \subset L_t \in \wls$ and
$L_t$ non singular. Let $l^1_t$ be the component of $(l_t - p_t)$ in 
the $e'$ side containing $l'$ 
union with $p_t$. Let $l^2_t$ be the other component
of $(l_t - p_t)$ union with $p_t$, see fig. \ref{onsid}, a.

The $l^1_t$ are rays (here we use $L_t$ non singular - but this
is just a technicality) and $(l^1_t)_{\infty}$ are not equal $v$ by
hypothesis. They cannot escape compact sets of $F$ since
$l'$ with ideal point $v$ is on that side of $e'$.
Hence as $t$ converges
to infinity $l^1_t$ converges to a leaf $l$ of $\wlsf$
with a ray (also denoted by $l$) with ideal point $v$ and
maybe some other leaves as well. The leaf
$l$ either shares a subray with $l'$ or separates $l'$ from $e'$
Let $e' \subset E$ leaf of $\wlu$ and $l \subset L$, leaf of $\wls$.

\vskip .1in
\noindent
{\underline {Case 2.2.1.1}} $-$ $l^2_t$ escapes in $F$
as $t \rightarrow \infty$.

Let $b_t$ be the ideal point of $l^2_t$. Then $b_t \not = v$.
Let $L^2_t$ be the union of $\wwp_{\rrrr}(p_t)$ and the component
of $L_t - \wwp_{\rrrr}(p_t)$ containing $l^2_t$. If $L^2_t$ escapes
in $\mi$, then the arguments in case 1 show 
continuity of 
$\Psi$ at $v$ in
the side of $e'$ not containing $l'$.

Now assume that $L^2_t$ converges to $R_1 \cup ... \cup R_m$
leaves of $\wls$
with union ${\cal R}$.
Notice $F$ may intersect some of these leaves or not.
If $\Theta(\Psi(p_t))$ does not escape in $\Theta(E')$,
then one of the $R_i$, call it $R_1$, is a leaf
intersecting $E'$.
As seen in the arguments for case 2.1.1, $F$ escapes
up in this direction so $\Psi(p_t)$ converges to
$(R_1)_+$.
If $\Theta(\Psi(p_t))$ escapes in $\Theta(E')$, then
lemma \ref{conperf} shows that $\Psi(p_t)$ also
converges to $(R_1)_+$.
This is equal to $(R_j)_+$ for any $j$.

Suppose there are $t_k \rightarrow \infty$ and 
$z_k$ in $l^2_{t_k}$ with $\Psi(z_k)$ not converging
to $(R_1)_+$. 
Here there is no need to assume that $L_{t_k}$ is non singular.
Then there $q$ in $\si$, \ $q \not = (R_1)_+$
and a subsequence still denoted by $\Psi(z_k)$,
so that $\Psi(z_k)$ converges to $q$.
As before there is $z$ in $\mi$ so that up to another 
subsequence, still denoted by $\wwp_{\rrrr}(z_k)$,
then $\wwp_{\rrrr}(z_k)$  converges
to $\wwp_{\rrrr}(z)$ and hence $z$ is in ${\cal R}$, say in
$R_i$.
Then $\wwp_{\rrrr}(z_k)$ are near 
$\wwp_{\rrrr}(z)$ and since a ray of $\Theta(R_i)$ is in
$\partial \Theta(F)$, then this is stable boundary. 
So $F$ escapes up as $\Theta(F)$ approaches $\Theta(z)$ and hence 
$\Psi(z_k)$ converges to $(R_i)_+$. 
This is equal to $(R_1)_+$.
The arguments of Case 2.1.1, situation 1 then show
continuity of $\Psi$ at $v$ on this side
of $e'$. This finishes the analysis of case 2.2.1.1.

\begin{figure}
\centeredepsfbox{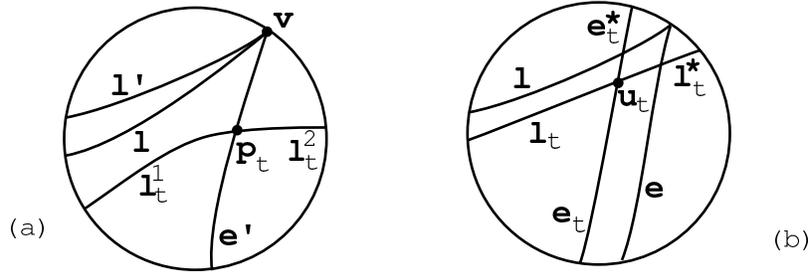}
\caption{
a. Convergence on one side,
b. Case 2.2.1.2 - intersection of leaves.}
\label{onsid}
\end{figure}

\vskip .1in
\noindent
{\underline {Case 2.2.1.2}} $-$ The \ $l^2_t$ \ limit to $r$ in $F$
as $t \rightarrow \infty$.

Choose the leaf $r$ with a ray which has ideal point $v$.
Then the leaves $r, l$ are not separated from each other in the leaf
space of $\wlsf$. 
Since $r, l$ have ideal point $v$ and there is no leaf
of $\wlsf$ non separated from $r, l$ and between them,
proposition \ref{nh} shows that the region bounded by these
rays of $r, l$ with ideal point $v$ projects in $M$ to a set
asymptotic to a Reeb annulus.
It follows that in $F$ this region is a bounded 
distance from
a geodesic ray with ideal point $v$. 
Now we restart
the process with the ray $r$ of $\wlsf$ instead of $e'$ of $\wluf$.
Let $\{ b_t, t \geq 0 \}$ be a parametrization of the 
leaves of $\wluf$ 
through the corresponding points $x_t$ of $r$.
If the  components of $(b_t - x_t)$ on the side opposite
of $e'$ escapes compact sets in $F$, then the analysis
of case 2.2.1.1 shows continuity of $\Psi$ at $v$ in that
side of $r$. Since $r$ and $e'$ are a bounded distance
from each other in $F$, this shows continuity of $\Psi$ at
$v$ on that side of $e'$.

Otherwise this process keeps being repeated.
Let $A_0 = L$, $A_1$ be the leaf of $\wls$ containing $r$.
If the process above does not stop, we keep producing $A_i$ in $\wls$,
so that they all disjoint and $A_i$ is non separated from
$A_{i+1}$. By theorem \ref{theb} 
up to covering translations there are only finitely many
leaves of $\wls$ which are not separated from some other leaf of
$\wls$.
There is then $m > n$ and $h$ covering translation with
$h(A_n) = A_m$. Let $f$ be the generator of the joint
stabilizer of $A_0, A_1$. 
This is non trivial by theorem \ref{theb}.
Then $f$ preserves all the prongs
of $A_1$ and therefore leaves invariant all the $A_i$.
Hence $h^{-1} f h(A_n) = A_n$ and so $h^{-1} f h = f^a$ for
some integer $a$. This implies there is a ${\bf Z} \oplus {\bf Z}$
in $\pi_1(M)$, see detailed arguments in \cite{Fe8}.
This is a contradiction.

There is then a last leaf $l_y$ (of $\wlsf$ or $\wluf$)
obtained from this process. The arguments of case 2.2.1.1 show
continuity of $\Psi$ at $v$ on the other side of $l_y$.
The region between $e'$ and $l_y$ is composed of a finite
union of regions between non separated rays of $\wlsf$ or
$\wluf$. They are all a bounded distance from a geodesic
ray with ideal point $v$, so the whole region also
satisfies this property.  It follows that this region
can only limit in $\Psi(v)$ as well and this proves
continuity of $\Psi$ at $v$ in that side of $e'$.

An entirely similar analysis shows continuity of $\Psi$ at
$v$ from the side of $l'$ not containing $e'$.

\vskip .1in

What remains to be analysed is the region of $F$ 
{\underline {between}} the rays $l'$ and $e'$.
Consider first the case that there is pair of non separated
leaves in the chain from $l'$ to $e'$.
Then as seen before the region between $l'$ 
and  $e'$
is a bounded distance (the bound is not uniform) from a geodesic
ray with ideal point $v$. 
This is not the case a priori if
there is no non Hausdorfness involved. 
In this case the region between $l'$ and $e'$ may not have
bounded thickness in $F$ and hence it is unclear whether
its image under $\Psi$ can only limit in $\Psi(v)$.
We analyse this case now.

In this last case 
parametrize the leaves of $\wluf$ intersecting the ray 
$l$ of $\wlsf$ as \ $\{ e_t \ | \ t \geq 0 \}$.
Since $l_t$ converges to $l$, then for big enough
$t$, the leaves $l_t, e_t$ intersect $-$ let $u_t$ be
their intersection point, see fig. \ref{onsid}, b.
Now define $l^*_t$ to be
the component of $l_t - u_t$  intersecting $e$ and $e^*_t$ 
the component of
$e_t - u_t$ intersecting $l$. 
Since $e'$ is on that side of $l$, the $e_t$ cannot escape
and converge to a leaf $e$ of $\wluf$ with an ideal point
$v$. Let $e \subset E$  leaf of $\wlu$.

Recall that $L_t$ is the leaf of $\wls$ containing $l^*_t$
and similarly let $E_t$ be the leaf of $\wlu$ containing $e_t$.
Let $L^*_t$ be the component
of $L_t - \wwp_{\rrrr}(u_t)$ containing $l^*_t$
and similarly define $E^*_t$.
In this remaining case the $l^*_t$ escape in $F$ and
so do the $e^*_t$. Hence
$\mu_t \ = \ l^*_t \cup \{ u_t \} \cup e^*_t$
defines a shrinking neighborhood system of
$v$ in $F \cup \pin F$. Consider the set

$$B_t \ = \ L^*_t \cup \wwr(u_t) \cup E^*_t$$

We want to show that $\overline B_t$ converges
to $L_+$ in the topology of closed sets of $\mi \cup \si$.

First consider $L^*_t \cap E$ which intersects $F$
in $(l^*_t \cap e)$. If $L^*_t \cap E$ does not
escape compact sets in $E$ then it limits to an orbit
$\gamma$ contained in a leaf $H$ of $\wls$.
Then $L, H$ are not separated from each other.
But for $t$ big enough then $E_t$ is near enough
$E$ and will intersect $H$ as well. This contradicts
$E_t \cap L$ is not empty and $L, H$ non separated.
Hence $L^*_t \cap E$ escapes in $E$ and similarly
$E^*_t \cap L$ escapes in $L$. 
Hence $L, E$ form a perfect fit.
This implies that $L_+ = E_-$.
Also $\Psi(e)$ limits to $E_-$ and $\Psi(l)$ limits
to $L_+ = E_-$.

The set $\overline L^*_t$ contains 
$(L^*_t \cap E)_+$ and this converges to $E_-$
when $t \rightarrow \infty$.
This is because $(L^*_t \cap E)$ escapes in $E$.
If $\overline L^*_t$ does
not converge to $E_-$ in $\mi \cup \si$, then
we find $t_k \rightarrow \infty$ and $x_k \in
L^*_{t_k}$ with $x_k$ converging to $x$ not equal to $E_-$.
Since $(x_k)_+ = (L_{t_k})_+$ converges to $E_-$,
then up to subsequence assume
$\wwr(x_k)$ converges to $\wwr(z)$ for some $z$ in $\mi$.
Then $z$ is in a leaf $H$ of $\wls$ which is non
separated from $L$.


The leaf $H$ does not intersect $F$, because $l^*_t$
escapes in $F$ by hypothesis in this final situation.
It follows that $\Theta(H)$ has a ray contained
in $\partial \Theta(F)$ and so this is stable boundary
of $\Theta(F)$. Hence $F$ escapes up as $\Theta(F)$
approaches $\Theta(H)$ and consequently 
$\Psi(x_k)$ limits to $H_+ = L_+ - = E_-$
$-$ which is what we wanted anyway.
This shows that $\overline L^*_t$ converges
to $E_-$ in $\mi \cup \si$.

Analysing the sets $E^*_t$ in the same manner we obtain
that $\overline E^*_t$ converges to $L_+$ as $t \rightarrow \infty$
as well.
This implies that $\overline B_t$ converges to 
$L_+ = \Psi(v)$. Since $B_t \cap F = \mu_t$ and the
$\mu_t$ define a neighborhood basis of $v$ in $F \cup \pin F$,
this shows continuity of $\Psi$ at $v$. 
This finishes the proof of case 2.2.1.2 and hence of
case 2.2.1.

\vskip .2in
\noindent
{\underline {Case 2.2.2}} $-$ There are rays $l$ of $\wlsf$ and
$e$ of $\wluf$ starting at $u_0$ and having the ideal point $v$.

We will first prove continuity on the side of $e$ not containing
a subray of $l$. There will be an iteration of steps.
Before we start the analysis we want to get rid of some
problems as described now. Suppose that there are $\alpha_0, 
\beta_0$ leaves of $\wlsf$ (or leaves of $\wluf$)
which have non separated rays converging to $v$ in $\pin F$
and on that side of $e$.
Suppose there are infinitely many of these on that side of
$e$. Let them be $\alpha_i, \beta_i$ and $G_i$ in
$\wls$ containing $\alpha_i$. Each region $B$ between $\alpha_0$
and any $\alpha_i$ is a bounded distance from a geodesic
ray in $F$ with ideal point $v$. The image $\Psi(B)$ then
can only limit in $\Psi(v)$.
If the $G_i$ do not escape in $\mi$ then they converge
to a leaf $G$ of $\wls$.
Let $A$ be an unstable leaf intersecting $G$ tranversely.
For $i$ big enough then $A$ intersects $G_i$ transversely,
which is impossible, as it would intersect
$\alpha_i$ and $\beta_i$ and these are not separated.
Hence the 
the $G_i$ escapes in $\mi$. Then  as seen in case 1,
there is continuity of $\Psi$ at $v$ in that side of $\alpha_1$.

Another situation is when there are leaves $\alpha_i$ in
that side of $e$ with two rays with ideal point $v$.
Then they are in the interior of a spike region $B$ with
one boundary $g$ with ideal point $v$. If there are infinitely
many of these, where none of the $\alpha_i$ are
nested with each other, then
 let $G_i$ in $\wls$ containing $\alpha_i$.
As in the previous paragraph,
the $G_i$ have to escape in $\mi$ and we have continuity
in that side of $\alpha_1$.

Therefore we can assume there are only finitely many
occurrences of spike regions or non separated leaves
with ideal point on this side of $e$. 
If there is any of these let $e_0$ be the last ray
in that side coming from such occurrences. 
Otherwise let 
$e_0$ be the ray given $e$ by the hypothesis in this case.
For simplicity assume that $e_0$ is a ray in $\wluf$,
the other case being similar. Let $e_0 \subset E_0 \in \wlu$.

Parametrize the ray of $e_0$ as 
\ $\{ p_t \ | \ t \geq 0 \}$ \
with $p_t$ converging to $v$ as $t \rightarrow \infty$. Let $l_t$ be
the leaf of $\wlsf$ through $p_t$ and $L_t$ in $\wls$
with $l_t \subset L_t$. If $L_t$ escapes $\mi$ as
$t \rightarrow \infty$ then as seen before we have continuity 
of $\Psi$ at
$v$ in that side of $e_0$.
So now suppose that $L_t$ converges to $A_1 \cup .... A_m$,
leaves of $\wls$. 
This case is considerably more involved, with several 
possibilities.

\vskip .1in
\noindent
{\underline {Claim}} $-$ $\Psi(e_0)$ converges to $(A_i)_+$ 
(notice the  $(A_i)_+, 1 \leq i \leq m$ are all equal).

If $E_0$ intersects some $A_i$, say $A_1$, then
as seen in case 2.1.1,
$F$ escapes positively along $\Psi(e_0)$ as $\Theta(F)$
approaches $A_1$. This implies that
$\Psi(e_0)$ converges to $(E_0 \cap A_1) _+ = (A_1)_+$.
If $E_0$ does not intersect any 
$A_i$ then $\Psi(e_0)$ converges 
to $(E_0)_- = (A_1)_+$. This proves the claim.

Let $l^1_t$ be the component of $(l_t - p_t)$ in the side
of $e_0$ we are considering.
We are really interested
in the behavior for $t \rightarrow \infty$, so we may
assume $p_t$ is not singular and there is only one
such component.

Suppose first that no $l^1_t$ has a ray with
ideal point $v$ and that
$l^1_t$ escapes in $F$ as $t \rightarrow \infty$.
In this case it is easy to show continuity 
of $\Psi$ at $v$ and in this side of $e_0$:
Suppose there are $x_i$ in $l^1_{t_i}$ with $t_i \rightarrow \infty$
and $\Psi(x_i) \not \rightarrow (A_i)_-$.
Since $(x_i)_+$ converges to $(A_i)_+$ then up to subsequence
assume that $(x_i)_- \rightarrow b \not = (A_i)_+$.
Up to subsequence $\wwr(x_i) \rightarrow \wwr(x)$. Then $x$
is in some $A_i$ say $x \in A_2$. But $F$ escapes positively
as 
$\Theta(F)$ approaches $\Theta(A_2)$,
so $\Psi(x_i) \rightarrow (A_i)_+$,  as we wanted.
Then as in case 2.1.1 this implies continuity.

There are 2 other options: 1) There is no $t$ with $l^1_t$
with an ideal point $v$ and $l^1_t$ does not escape in $F$; and
2) There is $t$ with $l^1_t$ having ideal point $v$.
These two options interact and intercalate in appearance
as explained below:

\vskip .1in
\noindent
{\underline {Situation 1}} $-$  There is no $t$ with $l^1_t$
with ideal point $v$ and $l^1_t$ does not escape in $F$.

There could be several leaves of $\wlsf$ in the limit of $l^1_t$ 
as $t \rightarrow \infty$ 
but there is a single leaf, call it $g$ with ideal point $v$.
If there is more than one such leaf with ideal point $v$, then
there would have to be one with
two rays with ideal point $v$. This leaf would be in a spike region
and it is separated from any other leaf in $\wlsf$,
contradiction.
Let $g$ be contained in a leaf 
 $G$ of $\wls$.

\begin{figure}
\centeredepsfbox{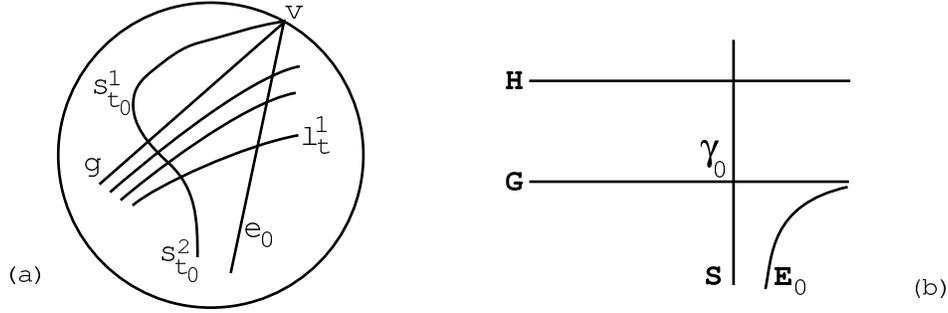}
\caption{
Some limits in $F$, b. The picture in $\mi$.}
\label{remo}
\end{figure}

Parametrize the ray $g$ as \ $\{ q_t \ | \ t \geq 0 \}$,
with $q_t \rightarrow v$ as $t \rightarrow \infty$.
Let $s_t$ be the unstable leaf of $\wluf$ through
$q_t$. 
Let $s^1_t$ be the component of $(s_t - q_t)$ on the
side of $g$ opposite to $e_0$ and $s^2_t$ the other
component. Then $s^2_t$ cannot have ideal point
$v$: for $t$ big enough it intersects $l^1_t$,
see fig. \ref{remo}, a.
Then $s^2_t$ converges to $e_0$.
By hypothesis there are no more occurrences of non separated
leaves of $\wlsf$
with ideal point $v$ on that side of $e_0$,
which implies that $s^1_t$ 
cannot limit to a leaf of $\wluf$ at $t \rightarrow \infty$
(it would distinct but non separated from $e_0$).
Hence the $s^1_t$ have to escape compact sets
in $F$. 
If $s^1_t$ does not have
an ideal point at $v$ for any $t$, then 
the previous analysis shows continuity of $\Psi$ at $v$ in
that side of $g$.
As in case 2.2.1.2
if $B$ is the region between $g$ and $e_0$ then
$\Psi(B)$ can only limit in $\Psi(v)$.

Hence assume 
there is some $t_0$ so that $s^1_{t_0}$ has ideal point $v$,
see fig. \ref{remo}, a.
Then for $t$ bigger than $t_0$ all ideal points of $s^1_t$
are $v$.
Let $s^1_{t_0}$ be contained in a leaf $S$
of $\wlu$ and $s_t$ contained in $S_t$ leaf of $\wlu$.
Since

$$l^1_t \ \rightarrow \ g, \ \ \ s^2_t \ \rightarrow \ e_0
\ \ \ \ {\rm when } \ \ t \rightarrow \infty,$$

$${\rm then} \ \ L_t \ \rightarrow \ G, \ \ \ 
S_t \ \rightarrow \ E_0, \ \ \ {\rm when} \ \ 
t \rightarrow \infty.$$

\noindent
It follows that $E_0, G$ form a perfect fit,
see fig. \ref{remo}, b.
Hence $(E_0)_- = G_+$.
If $\Theta(s^1_{t_0})$ is a ray in $\Theta(S)$ then
$\Psi(s^1_{t_0})$ converges to $S_-$.
But $\Psi(s^1_{t_0})$ also converges to 

$$\Psi(v)  \ = \ (E_0)_- \ = \ G_+ \ = \ (G \cap S)_+.$$

\noindent
Let $\gamma_0 = G \cap S$, an orbit of $\wwp$ in $G$.
The above equations imply that 

$$(\gamma_0)_+ \ = \ (G \cap S)_+  \ =  \
\Psi(v) \ = \ S_- \ = \ 
(\gamma_0)_-,$$

\noindent
which is a contradiction.
Hence $\Theta(s^1_t)$ is not a ray and has an endpoint
$x_1$ in $\Theta(S)$. Let $\gamma_1 = x_1 \times \rrrr$.
Let $H = \wls(\gamma_1)$. 
But $F$ does not intersect $H$.
If $F$ escapes down as $\Theta(F)$ approaches $x_1$, then
$\Psi(v) = (\gamma_1)_-$. 
But then

$$(\gamma_0)_- \ = \ (\gamma_1)_- \ = \ \Psi(v) 
\ = \ (\gamma_0)_+$$

\noindent
contradiction.
This implies that $F$ escapes up as $\Theta(F)$ approaches
$x_1$. 
Hence $\Theta(H)$ has a ray in $\partial \Theta(F)$.
Therefore $\Psi(s^1_t)$ limits to $(\gamma_1)_+$.
This implies that $(\gamma_0)_+ = (\gamma_1)_+$,
where $\gamma_0, \gamma_1$ are distinct orbits
of $\wwp$ in the {\underline {same}} unstable leaf $S$.
This is dealt with by the following theorem proved
in \cite{Fe5}:

\begin{theorem}{}{(\cite{Fe5})}
Let $\Phi$ be a quasigeodesic almost pseudo-Anosov flow in $M^3$ with
$\pi_1(M)$ negatively curved. Suppose there is
an unstable leaf $V$ of $\wlu$ and different orbits
$\beta_0, \beta_1$ in $V$ with $(\beta_0)_+ = (\beta_1)_+$.
Then $C_0 = \wls(\beta_0),  \ C_1 = \wls(\beta_1)$ are both periodic,
invariant under a nontrivial covering translation $f$,
and the periodic orbits in $C_0, C_1$ are connected
by an even chain of lozenges all intersecting $V$.
\label{ident}
\end{theorem}

\noindent
{\bf {Remark}} $-$
This result is case 2 of theorem 5.7 of \cite{Fe5}.
In that article the proof is done for quasigeodesic
Anosov flows in $M^3$ with $\pi_1(M)$ negatively
curved. The proof goes verbatin to the case of
pseudo-Anosov flows. The singularities make no
difference. By the blow up operation, the same holds
for almost pseudo-Anosov flows.
\vskip .1in

The theorem implies that $G, H$ are in the boundary of a chain
of adjacent lozenges all intersecting $S$. The first lozenge,
call it ${\cal C}$ has one stable side contained in $G$
and an unstable side $D_1$ which makes a perfect fit with
$G$. Suppose first $D_1$ is in the side of $S$ opposite
to $E_0$, see fig. \ref{spread}, a. The other unstable side
of ${\cal C}$ is a leaf $D_2$ which intersects $G$ on the
other side of $S$. Hence $G$ is some $S_c$ with $c > t_0$.
Then $S_c \cap F = s_c$ \ is a leaf of $\wluf$ and
$\Psi(s_c)$ has ideal point $\Psi(v)$. Notice that
$\Theta(s_c)$ (which is contained in $\Theta(F)$) escapes
in $\Theta(F)$ $-$ otherwise it would produce stable/unstable
boundary in $\Theta(F)$ before it hits $\Theta(H)$
and $\Theta(F)$ could not limit on $\Theta(H)$, impossible.
Hence $\Psi(s_c)$ limits to $(S_c)_-$ which is equal to $\Psi(v)$.
Then

$$(S_c \cap G)_- \ = \ (S_c)_- \ = \ \Psi(v) \ = \ G_+$$

\noindent 
which contradicts the orbit $S_c \cap G$ being a
quasigeodesic.


It follows that the perfect fits with $G$ occurs in the $E$ side
of $S$, see fig. \ref{spread}, b.
Here $\Theta(H), \Theta(D_1)$ are
contained in the boundary of $\Theta(F)$. 
We now look at the region $B$ in $F$ bounded by $s_{t_0} = S \cap F$
and $e_0 = E_0 \cap F$.

\vskip .1in
\noindent
{\underline {Claim 1}} $-$ 
The image $\Psi(B)$ can only
limit in $\Psi(v)$.

The region $\Psi(B)$ is contained in the region ${\cal E}$ of
$\mi$ which is bounded by $E, D_1$ (maybe other unstable leaves non
separated from $D_1$ as well), $H$ and $S$, see fig. \ref{spread}, b.
Notice that $F$ does not intersect $D_1$ or any leaf non separated from
$D_1$ which is beyond $D_1$. Otherwise $b_0 = (D_1 \cap F)$ is
contained in $B$ and non separated from $e_0$, so it would have
both ideal points $v$. Then it would be contained in the
interior of a spike region and could not be non separated
from another leaf $-$ impossible.
On the other hand since $\Theta(H)$ has a line leaf in
the stable boundary of 
$\Theta(F)$, then $\Theta(D_1)$ has a line
leaf in the unstable boundary of $\Theta(F)$.
Hence $F$ escapes down as $\Theta(F)$ approaches 
$\Theta(D_1)$.

Let $z_k$ in $B$ escaping in $F$ and hence
converging to $v$ in $\pin F$.
Suppose that $\Psi(z_k)$ does not converge to $\Psi(v)$.
Given that $z_k$ escapes $F$ and the structure of the region
${\cal E}$, it follows that up to subsequence
either $\wu(z_k)$ converges
to $D_1$ or $\ws(z_k)$ converges to $H$.
Suppose that $\ws(z_k)$ converges to $H$.
In that case $(z_k)_+$ converges to $H_+ = \Psi(v)$.
Then as seen before if $(z_k)_-$ does not converge
to $\Psi(v)$ we can assume up to subsequence $\wwr(z_k)$
converges to $\wwr(z)$. Then $z$ is in a leaf
non separated from $H$ and since $\Psi(z_k)$ has to
be in ${\cal E}$ then $z$ can only be in $H$.
As $F$ escapes up as $\Theta(F)$ approaches $\Theta(H)$ then
$\Psi(z_k)$ converges to $H_+ = \Psi(v)$.
The case $\wu(z_k)$ converges to $D_1$ leads to 
$\wwr(z_k)$ converging to $\wwr(z)$ with $z$ in
unstable leaf non separated from $D_1$. As $F$ escapes
down as $\Theta(F)$ approaches these unstable leaves,
then $\Psi(z_k)$ converges to $(D_1)_- = \Psi(v)$.
Since this works for any subsequence of $z_k$, then
$\Psi(z_k)$ has to converge to $\Psi(v)$ always.
This proves claim 1.

\begin{figure}
\centeredepsfbox{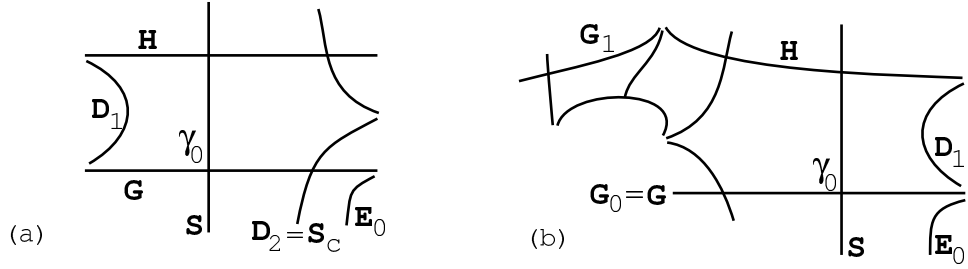}
\caption{
Perfect fit
with $G$ in the side opposite to $E_0$,
b. Perfect fit in the $E$ side.}
\label{spread}
\end{figure}

\vskip .1in
Let $G_0 = G$.  Notice that $G$ is periodic and connected to $H$ 
by an even chain of lozenges. 
We consider the ray $s_{t_0} = S \cap F$ which has
ideal point $v$. Parametrize it as \ $\{ z_t \ | \ t \geq 0 \}$.
Let $y_t$ be the leaf of $\wlsf$ through $z_t$ and $y^1_t$ the
component of $(y_t - z_t)$ in the side opposite to $e_0$.
The ray $s_{t_0}$ has the same behavior as the original ray
$e_0$. Hence we obtain continuity in that side of $s_{t_0}$
unless $y^1_t$ converges to a leaf $\mu$ of $\wlsf$ with
ideal point $v$. Let $G_1$ in $\wls$ with $\mu \subset G_1$.
Then $G_1$ is non separated from $H$, see fig. \ref{spread}, b \
and therefore connected to it by a chain of lozenges.
It follows that $G_1$ is connected to $G_0$ by a chain
of lozenges.
As in the proof of claim 1, the region $B_1$ of $F$
between $e_0$ and $(F \cap G_1)$ has image
$\Psi(B_1)$ which can limit only in $\Psi(v)$.

We restart the process with $g_1 = G_1 \cap F$ 
instead of $g$. The leaves of $\wluf$ through
points of $g_1$ already converge to the unstable
leaf $(D_3 \cap F)$ of $\wluf$ \ ($D_3$ is depicted in fig.
\ref{spread}, b).
The leaf $(D_3 \cap F)$ cannot be non separated from any
other leaf of $\wluf$ in that side of $(D_3 \cap F)$.
It follows that 
the unstable leaves intersected by $g_1$ escape
in $F$. The only case to be analysed is that some of
these unstable leaves have ideal point $v$. This brings
the process exactly to the situation of some $s^1_t$
of $\wluf$ having ideal point $v$ as described before
(it was $s^1_{t_0}$). So this would produce 
$H_1$ of $\wls$ with similar properties as $H$.
This process can now be iterated.
As in claim 1 the region of $F$ between $g_i$ and $g_{i+1}$
maps to $\mi$ to a region which can only limit in $\Psi(v)$.

We show that this process has to stop. Otherwise produce
$G_i$ leaves of $\wls$ which are all connected to $G_0$ by
a chain of lozenges. The $G_i$ are all non separated from
some other leaf of $\wls$, Hence there are $G_i, G_j$ which
project to the same stable leaf in $M$.
There is a covering translation $h$ taking $G_i$ to $G_j$.
If $f$ is a generator of the isotropy group of $G_0$ leaving
all sectors invariant, then it leaves invariant
all lozenges in any chain starting in $G_0$ so leaves
invariant all the $G_i$. As before this leads to $h^{-1} f h = f^n$
for some $n$ in ${\bf Z}$ and to a ${\bf Z} \oplus {\bf Z}$ in
$\pi_1(M)$. This is disallowed.
Therefore the process finishes after say $j$ steps and we
obtain continuity of $\Psi$ at $v$ in that side of $g_j = G_j \cap F$.
As seen above the region between $s_{t_0}$ and $g_j$ maps by $\Psi$
into a region that can only limit in $\Psi(v)$. 
This proves continuity of $\Psi$ at $v$ in that side of $e_0$.
This finishes the analysis of situation 1.

\vskip .15in
\noindent
{\underline {Situation 2}} $-$ There is $l^1_{t_0}$ with ideal
point $v$.

Recall the setup before the analysis of situation 1. 
Let $\{ u_t \ | \ t \geq 0 \}$  be the collection
of unstable leaves intersected by the ray $l^1_{t_0}$.
The analysis is 
extremely similar to the analysis of situation 1, which shows all
cases produce continuity in the first step except
when $u_t$ converges to a leaf $u$ of $\wluf$ 
with ideal point $v$.
Then consider the stable leaves intersecting $u$.
The analysis of situation 1 shows continuity unless there 
is stable leaf
with ideal point $v$.
From now on the analysis is exactly the same as in situation
1, with unstable replaced by stable and vice versa.

%
%

\vskip .15in
So far we proved continuity of $\Psi$ at $v$
from the side of $e_0$ opposite to $l$.
The same works for the other side of $l$, producing
$l_0$ with similar properties as $e_0$.
We now must consider the regions between $e_0$ and
$e$, between $e$ and $l$ and between $l$ and $l_0$.

First consider the region between $e$ and $e_0$,
which occurs only when they are distinct. This implies
that the ray $e_0$ is a bounded distance from a
geodesic ray in $F$ with ideal point $v$.
Let $\{ \mu_t \ | \ t \geq 0 \}$ be a parametrization
of the stable leaves of $\wlsf$ through $e$.
Let $\mu^1_t$ be the component of $(\mu^1_t - e)$ 
in the side of $e$ we are considering. If some 
$\mu^1_t$ has ideal point $v$ then both ideal
points of $\mu_t$ are $v$ and $\mu_t$ is inside
a spike region. The same is true for $e$ and
so $e$ is a bounded distance from a geodesic
ray in $F$ with ideal point $a$.
Hence the region between $e$ and $e_0$ is a bounded
distance from a geodesic ray and we are finished in this
case.


The remaining case to be analysed here is that $\mu^1_t$
has no ideal point $v$. Then $\mu^1_t$ does not escape
$F$ as $t \rightarrow \infty$, because $e_0$ is in that
side of $e$.
So $\mu^1_t$ converges to a leaf $\mu$ which has ideal
point $v$.
Now consider a parametrization $\{ \nu_t \ | \ t \geq 0 \}$ of
the unstable leaves intersected by $\mu$.
Then $\nu_t$ converges to the leaf $e$.
If it converges to some other leaf, then $e$ is
a bounded distance from a geodesic ray in $F$ and we are done.
Otherwise it must be that some $\nu_t$ has ideal
point $v$. Therefore we exactly in the setup analysed
in situation 1 above. 

This shows continuity of $\Psi$ for the region between $e$ and
$e_0$ and similarly for the region between $l$ and $l_0$.

Finally we analyse the region $B$ between $e$ and $l$.
First notice there is no singularity in the interior
of $B$. Otherwise there would be a line leaf in
$B$ and hence a leaf with both endpoints $v$. 
It would have to be part of a spike region and the
spike region does not have any singularities in its 
interior.

Parametrize the leaves of $\wluf$ through $l$ as
$\{ e_s \ | \ s \geq 0 \}$ and similarly
those of $\wlsf$ through $e$ as $\{ l_t \ | \ t \geq 0 \}$.
Let $L, L_t$ leaves of $\wls$ with $l \subset L, \ l_t
\subset L_t$ and similarly define $E, E_t$.
There are 2 possibilities:

\vskip .1in
\noindent
1) Product case $-$ Any $l_t$ intersects every $e_s$ and
vice versa. 

Equivalently $\wlsf, \wluf$ define  a product structure 
in the region $B$ bounded by $l_0$ and $e_0$. If the $L_t$ escapes
in $\mi$ as $t \rightarrow \infty$, then there is a stable
product region defined by a segment in $L_0$. 
But then theorem \ref{prod} implies that
$\Phi$ is topologically conjugate to a suspension,
contradiction.
It follows that the 
$L_t$ converge to $H_1 \cup ... H_m$ as $t \rightarrow \infty$.
Since the $l_t$ are stable leaves, it follows that
$F$ escapes up as $\Theta(F)$ appraches $\Theta(H_i)$.
This implies that $\Psi(e)$ limits to $(H_i)_+$ which
is then equal to $\Psi(v)$.
Similarly $E_s$ converges to $V_1 \cup ... V_n$ and
$F$ escapes down as $\Theta(F)$ approaches $\Theta(V_j)$. 
Hence 
$\Psi(l)$ limits to $(V_j)_- = \Psi(v)$.
If some $H_i$ intersects some
$V_j$, then

$$(V_j \cap H_i)_+ \ = \ (H_i)_+ 
\ = \ \Psi(v) \ = \ 
(V_j)_- \ = \ (V_j \cap H_i)_-,$$

\noindent
contradiction.
Let now $\{ z_k \}$ be a sequence in $B$ converging
to $v$. The product structure implies that up to subsequence
we may assume that either $\ws(z_k)$ converges to
$H_i$ or $\wu(z_k)$ converges to $V_j$. 
This is analysed carefully in Claim 1 above, which
shows that $\Psi(z_k)$ must converge to $\Psi(v)$.
This shows continuity of $\Psi$ when restricted to the
region $B$.

%
%
%
%
%

\vskip .1in
\noindent
2) Non product case.

There are $t, u > 0$ with $l_t \cap e_u = \emptyset$.
Consider one such $u$.
Let $a$ be the infimum of $t$ with $l_t \cap e_u = \emptyset$.
Now let $b$ be the infimum of $u$ with $l_a \cap e_u = \emptyset$.
Then $l_a \cap e_b = \emptyset$, but for any $0 \leq t < a$
and $0 \leq u < b$ one has 
$l_t \cap e_u \not = \emptyset$. 
Since $l_a \cap e_b = \emptyset$, then $L_a \cap E_b = \emptyset$.
It follows that $L_a, E_b$ form a perfect fit,
see fig. \ref{nene}, a.
If $\Theta(l_a)$ does not escape in $\Theta(L_a)$, then
there would be unstable boundary of $\Theta(F)$ in the
limit and that would keep $F$ from intersecting $E_b$,
contradiction. Hence $\Theta(l_a)$ escapes in $\Theta(L_a)$
and $\Theta(e_b)$ escapes in $\Theta(E_b)$.
Hence $\Psi(l_a)$ limits to $(L_a)_+$ and $\Psi(e_b)$ limits
to $(E_b)_-$. 
Also $l_a, e_b$ limit to $v$ in $\pin F$. 

Let $p_t = l_t \cap e$. If $\Theta(p_t)$ escapes in $\Theta(E)$,
then $\Psi(e)$ converges to $E_-$. Notice that $\Psi(e)$ converges
to $\Psi(v)$ so:

$$E_- \ = \ \Psi(v) \ = \ (L_a)_+ \ = \ (L_a \cap E)_+$$

\noindent
contradiction.
It follows that $\wwr(p_t)$ converges to
$\wwr(x)$ with $x$ in $E$.
Also $F$ has to escape up as $\Theta(F)$ approaches $\Theta(x)$
$-$ same as in Situation 1 above.
Hence $\Psi(e)$ limits to $x_+$. So

$$x_+ \ = \ \Psi(v) \ = \ (L_a)_+ \ = \ (p_a)_+$$

\noindent
Let $X = \wws(x)$.
Then $x, p_a$ are in 2 distinct orbits of $E$ with
the same positive ideal point.
Therefore theorem \ref{ident} implies that
$L_a,  X$ are connected by an even chain of lozenges
all intersecting $E$.
Let ${\cal C}$ be the first lozenge. It has a stable
side in $L_a$ and one unstable side, call it $D_1$ which
makes a perfect fit with $L_a$. Suppose first that $D_1$ is
in the component of $\mi - E$ opposite to $E_b$.
Then the other unstable side of ${\cal C}$, call it $D_2$
has to intersect $L_a$ in the other side of $E$.
Then $D_2$ must be some $E_t$,
let it be $E_{b'}$,
 see
fig. \ref{nene}, a. 
Then $\Theta(e_{b'})$ has to escape in $\Theta(E_{b'})$ 
or else one produces stable boundary to $\Theta(F)$ and
$\Theta(F)$ cannot limit to $\Theta(x)$ contradiction.
Hence $\Psi(e_{b'})$ converges to $\Psi(v)$ and 
also to $(E_{b'})_-$.
But then 

$$(E_{b'} \cap L_a)_- \ = \ 
(E_{b'})_- \ = \ \Psi(v) \ = \ (E_b)_- \ = \ (L_a)_+$$

\noindent
again a contradiction.

\begin{figure}
\centeredepsfbox{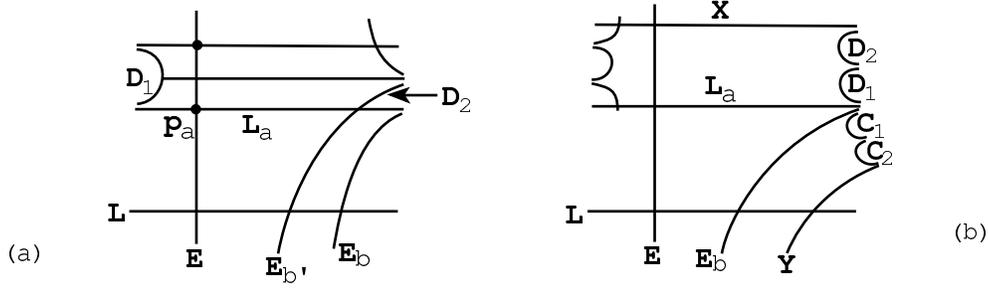}
\caption{
a. Reaching before, b. Reaching at the exact time.}
\label{nene}
\end{figure}

This implies that $D_1$ is on the side of $E$ containing
$E_b$, see fig. \ref{nene}, b.

If there are only 2 lozenges in the chain from $L_1$ to
$X$, then  $D_1$ also makes a perfect fit with
$X$. Otherwise there are $D_2, ..., D_i$ all non separated
from $D_1$ and so that $D_i$ makes a perfect fit 
with $X$ and the $D_j$ are all in the boundary of the
chain of lozenges.
As seen in claim 1 above, $F$ cannot intersect
any $D_j$ ($1 \leq j \leq i$), but all
$\Theta(D_j)$ are contained in the unstable boundary
of $\Theta(F)$. Also $F$ escapes down as $\Theta(F)$
limits to $\Theta(D_j)$.
The set $\Theta(X)$ also has a line leaf which is
a stable boundary of $\Theta(F)$ and $F$ escapes
up when $\Theta(F)$ approaches $\Theta(X)$.

The same discussion applies to $L$, so there is 
$y$ in $L$, $Y = \widetilde W^u(y)$ with $\Theta(Y)$
having a line leaf in the unstable boundary
of $\Theta(F)$ and $F$ escapes down accordingly.
There are $C_1, ..., C_n$ leaves in $\wlu$, 
all non separated from each other and in the
boundary of the lozenges in the chain from $E_b$ to
$Y$ so that $C_1$ makes a perfect fit with $E_b$ and
$C_n$ makes a perfect fit with $Y$, see fig. \ref{nene}, b.
Finally $\Theta(C_j)$ has a line leaf in the
stable boundary of $\Theta(F)$ and $F$ escapes up
accordingly.

Let ${\cal E}$ be the region in $\mi$ bounded
by

$$E, \ L, \ X, \ Y, \ C_1, ..., C_n, \ D_1, ..., D_i$$

\noindent
Then ${\cal E} \cap F$ is exactly the region $B$
bounded by the rays $e$ and $l$.
Let $z_k$ in $B$ escaping to $v$.
Then the region ${\cal E}$ shows that
up to subsequence one of the following must occur:

\vskip .08in
1) $\ws(z_k)$ converges to either $X$ or $C_1$. 
The analysis of claim 1 above shows that
$\Psi(z_k)$ converges to either $X_+$ or $(C_1)_+$
both of which are equal to $\Psi(v)$.

2) $\wu(z_k)$ converges to either $Y$ or $D_1$.
Here $\Psi(z_k)$ converges to either $Y_-$ or
$(D_1)_-$ both of which are equal to $\Psi(v)$.

In any case this shows continuity of $\Psi$ in the
region $B$.
This finishes the non product case.

\vskip .05in
This finishes the proof of theorem \ref{exten},
the continuous extension theorem.
\end{proof}


%
%
%
%

%
%
%
%

\section{Foliations and Kleinian groups}

There are many similarities between foliations in hyperbolic
3-manifolds and Kleinian groups.
We refer to \cite{Mi,Can,Mar} for basic definitions concerning
degenerate and non degenerate Kleinian groups, in particular
singly and doubly degenerate groups.

If the foliation is $\rrrr$-covered then the limit set
of any leaf in $\mi$ is the whole sphere \cite{Fe3}. 
This corresponds to doubly degenerate surface
Kleinian groups \cite{Th1,Mi,Can,Mar,Bon}.
There is always a pseudo-Anosov flow which is transverse
to the foliation \cite{Fe7,Cal1}. If the flow is quasigeodesic then
the results of this article imply that the foliation
has the continuous extension property.

If the foliation has one sided branching, say branching down,
then limit sets of leaves can only have domain of discontinuity
``above" \cite{Fe3}. 
Let $F$ in $\fn$ and $\Lambda_F$ its limit set. If $p$ is not
in $\Lambda_F$, the $p$ is said to be {\em above} $F$ if
there is a neighborhood $V$ of $p$ in $\mi \cup \si$, so
that $V \cap \mi$ is on the positive side of $F$.
This corresponds to simply degenerate surface Kleinian
groups \cite{Th1,Mi,Can}. There are examples of foliations with one sided
branching transverse to 
suspension pseudo-Anosov flows
provided by Meigniez \cite{Me}. Suspension flows are always
quasigeodesic flows \cite{Ze}. 
The results of this article show the continuous extension property
for such foliations.
Under these conditions, the limit
sets are locally connected, the continuous extension provides
parametrizations of these limit sets.

Finally if there is branching in both directions, 
then there can be domain of discontinuity above
and below leaves. This corresponds to non degenerate
Kleinian groups \cite{Th1,Mi,Can}. These occur for example
in the case of finite depth foliations, where
the depth 0 leaves are not virtual fibers \cite{Fe6}.

There are many interesting questions:

\vskip .1in
\noindent
{\bf {Question 1}} $-$ Given a 
foliation $\fol$, is it  $\rrrr$-covered if and only if 
for every $F \in \fn$ then the limit set $\Gamma_F$ is $\si$?

The forward direction is true. The backwards direction
is true if there is a compact leaf \cite{Go-Sh}. In addition if there is one
leaf with limit set the whole sphere then all leaves have
limit set the whole sphere \cite{Fe3} $-$ whether $\fol$ is
$\rrrr$-covered or not.

\vskip .1in
\noindent
{\bf {Question 2}} $-$ Given $\fol$ an $\rrrr$-covered foliation, is there 
a quasigeodesic transverse pseudo-Anosov flow?

This is true in the case of slitherings or uniform foliations
as defined by Thurston \cite{Th5}. Examples are fibrations,
$\rrrr$-covered Anosov flows and many others.
There is always a transverse pseudo-Anosov flow, the
question is whether it is quasigeodesic.

\vskip .1in
\noindent
{\bf {Question 3}} $-$ Is there domain of discontinuity of $\Lambda_F$ only
above $F$ if and only if $\fol$ has one sided branching
in the negative direction?

This occurs for the examples constructed by Meigniez \cite{Me}.

\vskip .1in
\noindent
{\bf {Question 4}} $-$ Are the pseudo-Anosov flows constructed
by Calegari \cite{Cal2} and transverse to one sided branching foliations
quasigeodesic?

\vskip .1in
\noindent
{\bf {Question 5}} $-$ If $\fol$ has 2 sided branching is there
always domain of discontinuity above and below?
Is there a quasigeodesic pseudo-Anosov flow almost transverse to $\fol$?

{\footnotesize
{
\setlength{\baselineskip}{0.01cm}

\noindent
Florida State University

\noindent
Tallahassee, FL 32306-4510, USA

}
}

\end{document}